\documentclass{elsarticle}

\usepackage{lineno,hyperref}
\usepackage{amsmath}
\usepackage{amsfonts,amssymb}
\usepackage{esint}
\usepackage{subfigure}
\usepackage{wasysym}
\usepackage{amsthm}
\usepackage{multirow}
\usepackage{algorithm}
\usepackage{algorithmicx}
\usepackage{algpseudocode}
\usepackage{natbib}
\usepackage{color}
\setcitestyle{numbers}
\usepackage{graphicx}
\usepackage{xcolor}
\usepackage{array}
\usepackage{geometry}
\newtheorem{remark}{Remark}

\modulolinenumbers[5]

\begin{document}

\begin{frontmatter}

\title{A multi-resolution limiter for the Runge-Kutta discontinuous Galerkin method}


\author[add1]{Hua Shen}
\ead{huashen@uestc.edu.cn}
\author[add2]{Bangwei She\corref{mycorrespondingauthor}}
\cortext[mycorrespondingauthor]{Corresponding author}
\ead{bangweishe@cnu.edu.cn}

\address[add1]{School of Mathematical Sciences, University of Electronic Science and Technology of China, Chengdu, Sichuan 611731, China}
\address[add2]{Academy for Multidisciplinary studies, Capital Normal University, West 3rd Ring North Road 105, 100048 Beijing, China}

\begin{abstract}
 We propose a novel multi-resolution (MR) limiter for the Runge-Kutta discontinuous Galerkin (RKDG) method for solving hyperbolic conservation laws on a general unstructured mesh. 
Unlike classical limiters, which detects only solution discontinuities to dichotomize cells into good or troubled, the proposed MR limiter also takes into account the derivative discontinuities to divide cells into several groups.  The method operates by performing a successive comparison of the local DG polynomial's derivatives, from high-order to low-order, against a baseline constructed from neighboring cell averages. 
  If a $k$th-order derivative of the DG polynomial is larger than the baseline, 
  then we reduce the order to $k-1$ and set the corresponding $k$th-order terms to be 0; 
  Otherwise, the remaining $k$th-order DG polynomial is used to represent the final solution.
  Only if all the derivatives are larger than the baseline, a TVD slope limiter is used to reconstruct the solution. 
In this manner, the limiter dynamically selects an optimal polynomial suited to the local solution smoothness without problem-dependent parameter to tune.   Notably, it also possesses a scale-invariance property that is absent in most classical limiters. A series of numerical examples demonstrate the accuracy and robustness of the proposed MR limiter.
\end{abstract}

\end{frontmatter}


\section{Introduction}\label{sec:intro}
The discontinuous Galerkin (DG) method was initially proposed by Reed and Hill \cite{Reed1973DG}
for solving the neutron transport, which is governed by steady-state linear hyperbolic equations.
One of the main advantages of the DG method is that it can easily achieve high-order accuracy
on unstructured meshes due to its compactness and finite element feature. 
However, when solving nonlinear time-dependent hyperbolic conservation laws, high-order linear DG schemes
inevitably generate numerical oscillations near discontinuities.
Cockburn and Shu and their coworkers \cite{Cockburn1989RKDGII, Cockburn1989RKDGIII, Cockburn1990RKDGIV, Cockburn1998RKDGV} 
systematically developed the DG method by adopting 
the efficient total variation diminishing (TVD) Runge-Kutta schemes \cite{Shu1988EfficientENO,Shu1989EfficientENOII}
for time discretization and constructing a total variation bounded (TVB) limiter \cite{Shu1987TVB} for shock capturing.
Their work has significantly promoted the application of the DG method for solving nonlinear hyperbolic conservation laws.

The TVB limiter \cite{Cockburn1989RKDGII, Cockburn1989RKDGIII, Cockburn1990RKDGIV, Cockburn1998RKDGV} first
detects the cells with a large variation in the solution and then modifies the solution by a TVD limiter.
In this way, the TVB limiter can suppress spurious oscillations near discontinuities,
but it contains a free parameter $M$ which is difficult to choose for specific problems.
In view of this, some improved limiters were proposed, such as the moment-based limiters \cite{Biswas1994, Burbeau2001}.
However, these limiters always have a certain probability of 
mistakenly replacing the high-order DG solutions in smooth regions with a low-order reconstruction. 
Generally speaking, there are two ways to address this issue.
The first way is using a good limiter to reconstruct the solutions in troubled cells,
so the detection of a good cell as a troubled one does not significantly affect the accuracy. 
Qiu and Shu \cite{Qiu2004WENOLimiter,Qiu2005WENOLimiter,Qiu2005ComparisonIndicator} performed 
WENO reconstructions with the same order as the DG solution in the troubled cells 
detected by a troubled-cell indicator. 
They systematically compared different troubled-cell indicators
and concluded that no indicator can defeat others in every problem \cite{Qiu2005ComparisonIndicator}.
Overall, the minmod-based TVB indicator \cite{Cockburn1989RKDGII} 
with a suitably chosen constant, the indicator based on Harten's subcell resolution idea \cite{Harten1989Subcell}, 
and the KXRCF shock detector proposed by Krivodonova \emph{et al.} \cite{Krivodonova2004KXRCF} performed better 
than other indicators existing at that time.
Later on, Zhu \emph{et al.} \cite{Zhu2008WENOLimiterTri,Zhu2012WENOLimiter3D} implemented WENO limiters 
for the DG schemes on two- and three-dimensional unstructured meshes.
Although WENO reconstructions can preserve the high-order accuracy even when some good cells are mistakenly flagged as troubled ones,
the high-order WENO reconstructions on unstructured meshes become complex due to the large stencils. 
To address this issue, Zhong and Shu \cite{Zhong2013SimpleDG_limiter} proposed a simple limiter 
that adopts the classical WENO weights to convexly combine the entire DG polynomials at the local cell and its neighboring cells.
In this way, the compactness of the DG method can be maintained.
Zhu \emph{et al.} \cite{Zhu2013SimpleWENOLimiterTri} further extended the simple WENO limiter \cite{Zhong2013SimpleDG_limiter} to unstructured meshes.
Dumbser \emph{et al.} \cite{Dumbser2014subcell} proposed a posteriori subcell limiter 
which first detects troubled cells based on the so-called MOOD paradigm \cite{Clain2011MOOD}
and then recomputes the discrete solutions in the troubled cells 
by using a finite volume scheme based on subcells.
Du \emph{et al.} \cite{Du2022ImprovedWENOLimiterTri} improved the efficiency of the simple WENO limiter \cite{Zhong2013SimpleDG_limiter}
by reducing the number of polynomials transformed to the characteristic fields for each direction.
Recently, Wei and Xia \cite{Wei2024HybridLimiter} proposed a jump indicator coupled with a hybrid limiter 
which is compact and less dissipative than the TVD limiter.
Furthermore, Wei \emph{et al.} \cite{Wei2025JumpFilter} proposed a compact filter 
which can be seamlessly integrated into the hybrid limiter 
and is compatible with the other damping techniques 
such as the oscillation-free DG (OFDG) \cite{Lu2021OFDG,Liu2022OFDG} 
and the oscillation-eliminating DG (OEDG) \cite{Peng2025OEDG,Ding2025OEDG}.
Another way to improve the limiter is to use a better troubled-cell indicator to detect troubled cells more accurately,
so that we can use a simple slope limiter to reconstruct the solutions in troubled cells.
The KXRCF shock detector proposed by Krivodonova \emph{et al.} \cite{Krivodonova2004KXRCF}
is a widely used troubled-cell indicator and generally performs well for many cases 
without tuning a problem-dependent parameter. 
Fu and Shu \cite{Fu2017NewLimiter} proposed a compact troubled-cell indicator for the RKDG method, 
which can identify shocks with a parameter only depending on the degree of the DG polynomial.
We refer to this limiter simply as the FS limiter.
Vuik and Ryan \cite{Vuik2016automated} treated the troubled-cell detection process
as an outlier-detection problem and utilized Tukey's boxplot approach to automatically 
tune parameters for various troubled-cell indicators.
Ray and Hesthaven \cite{Ray2018ANN,Ray2019ANN} trained an artificial neural network as a troubled-cell indicator.
Zhu \emph{et al.} \cite{Zhu2021K_means} used the K-means clustering algorithm to detect troubled cells 
without tuning problem-dependent parameters. 

The limiters mentioned above focus on suppressing the spurious oscillations 
induced by the high-order DG schemes, but they cannot guarantee the boundedness of the solutions,
such as the positivity of the density and the internal energy. 
An elegant bound-preserving limiter for the high-order DG schemes was proposed 
by Zhang and Shu \cite{Zhang2010MaximumDG,Zhang2010PositivityDG,Zhang2011MaximumReview,Zhang2012MaximumDGTri,Zhang2017PositivityDGNS}.
A shock-capturing limiter can be used together with the bound-preserving limiter
to guarantee a strong stability for the DG method.

Inspired by our previous work on the ENO-MR finite difference schemes \cite{Shen2025ENO_MR}, 
we propose a multi-resolution (MR) shock-capturing limiter for the RKDG method 
for solving hyperbolic conservation laws.
The new limiter is based on the fact that a solution to nonlinear hyperbolic conservation laws may contain different types of discontinuities 
which can be approximated by polynomials of different degrees.
Therefore, it is not the best choice to categorize the cells as either `good' or `troubled' cells.
There should be some cells falling between these two categories.
This is probably the main reason why traditional limiters are difficult to perform well 
for different problems without tuning problem-dependent parameters.
In view of this, the proposed MR limiter tags the cells with a certain order of accuracy 
by detecting the discontinuities of the solution as well as its derivatives using a simple strategy.
In this way, the limiter can choose an optimal polynomial to approximate the solutions with different smoothness levels, 
so it can efficiently capture different types of discontinuities while preserving high-order accuracy for smooth solutions
on both quadrilateral and triangular meshes without problem-dependent parameters to tune.
Additionally, the limiter has the scale-invariant property, which is lacking in most of the classical limiters.

\section{A brief description of the Runge-Kutta discontinuous Galerkin method}\label{SEC:RKDG}
We consider the following time-dependent hyperbolic conservation law,
\begin{equation}\label{Eq:HCL}
  \begin{cases}
    \frac{\partial u}{\partial t}+\nabla\cdot f(u)=0,\\
    u(\mathbf{x},0)=u_0(\mathbf{x}),
  \end{cases}
\end{equation}
where $\mathbf{x}=(x_1,...,x_D)\in\mathbb{R}^D$, $u:\mathbb{R}^{D+1}\to \mathbb{R}^M$,
and $f:\mathbb{R}^M\to \mathbb{R}^{MD}$.
The computational domain is divided into $N$ non-overlapping elements $\{K_j\}_{j=1}^N$.
The DG solution $u_h$ is defined
in the discrete space $V_h^k=\{v:v\in P^k(K_j), 1\le j\le N\}$, 
where $P^k(K_j)$ denotes the set of polynomials of degree at most $k$ defined on $K_j$.
Then the discrete solution can be written as 
\begin{equation}\label{Eq:Solution}
  u_h(\mathbf{x},t)=\sum_{l=0}^{k} u_j^l(t)v_j^l(\mathbf{x}), \quad\mathbf{x}\in K_j,
\end{equation}
where $\{v_j^l\}_{l=0}^k$ are the basis functions for $P^k(K_j)$
and $\{u_j^l\}_{l=0}^k$ are the degrees of freedom to be solved.
Substituting the discrete solution into the hyperbolic conservation law (\ref{Eq:HCL})
and enforcing the projection of the residual onto the test functions $v_j^m$ $(1\le m\le k)$ equal to 0,
we obtain the semi-discrete DG scheme as
\begin{subequations}\label{Eq:DG_scheme}
  \begin{equation}\label{Eq:DG_scheme_time_evolution}
    \sum_{l=0}^{k}a_{ml}\frac{du_j^l(t)}{dt}=\int_{K_j} f(u_h(\mathbf{x},t))\cdot\nabla v_j^m(\mathbf{x})dV
  -\int_{\partial K_j} f(u_h(\mathbf{x},t))\cdot n v_j^m(\mathbf{x})ds,
\end{equation}
\begin{equation}\label{Eq:DG_scheme_initial_condition}
  \sum_{l=0}^{k}a_{ml}u_j^l(0)=\int_{K_j}u_0(\mathbf{x})v_j^m(\mathbf{x})dV,
\end{equation}
\end{subequations}
where $dV=dx_1...dx_D$ and $a_{ml}=\int_{K_j} v_j^m(\mathbf{x})v_j^l(\mathbf{x})dV$.
We usually adopt a set of locally orthogonal basis on ${K_j}$ such that $a_{ml}=\delta_{ml}$
which significantly simplifies the computations. 
The integral terms in (\ref{Eq:DG_scheme}) are computed by numerical quadratures 
which are exact for polynomials of degree up to $2k$ for the element integral 
and up to $2k+1$ for the boundary integral.
Since the solution at the element boundaries $\partial K_j$ is discontinuous,
the flux term $f(u_h(\mathbf{x},t))\cdot n$ in the boundary integral 
is approximated by a suitable numerical flux.
In this study, we adopt the widely used Lax-Friedrichs flux.

Denote $\mathbf{u}_j=(u_j^0,u_j^1,...,u_j^k)$,
then the initial value problem (\ref{Eq:DG_scheme}) is solved by the third-order 
strong-stability-preserving Runge-Kutta (SSP-RK3) scheme \cite{Shu1988EfficientENO} as 
\begin{subequations}\label{Eq:TVD3rdRK}
  \begin{equation}\label{SubEq:RKStage1}
    \mathbf{u}_j^{(1)}=\mathbf{u}_j^n+\Delta t\mathcal{L}_h(\mathbf{u}^n),
  \end{equation}
  \begin{equation}\label{SubEq:RKStage2}
    \mathbf{u}_j^{(2)}=\frac{3}{4}\mathbf{u}_j^n+\frac{1}{4}\left(\mathbf{u}_j^{(1)}+\Delta t\mathcal{L}_h(\mathbf{u}^{(1)})\right),
  \end{equation}
  \begin{equation}\label{SubEq:RKStage3}
    \mathbf{u}_j^{n+1}=\frac{1}{3}\mathbf{u}_j^n+\frac{2}{3}\left(\mathbf{u}_j^{(2)}+\Delta t\mathcal{L}_h(\mathbf{u}^{(2)})\right),
  \end{equation}
  \end{subequations}
where $\Delta t$ is the time step and $\mathcal{L}_h$ represents the spatial operator.
\section{A multi-resolution limiter for the RKDG method}\label{SEC:MR_limiter}
\subsection{The description of the algorithm}\label{SUBSEC:MR_limiter}

 \begin{figure}
  \centering
  \includegraphics[width=5.5 cm]{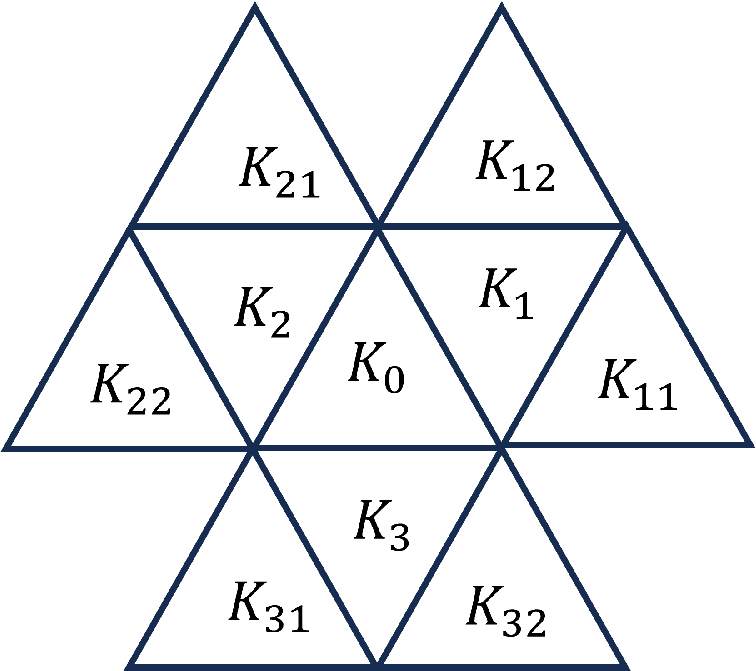}
  
  \caption{The indicator stencil.}
 \label{FIG:Stencil}
 \end{figure}
Without loss of generality, we describe the limiter for 2D unstructured triangular meshes.
For the convenience of explanation, we relabel the target cell as $K_0$ 
and define the associated indicator stencil as $S:=\{K_0,K_1,K_{11},K_{12},K_2,K_{21},K_{22},K_3,K_{31},K_{32}\}$, that is the set of all neighboring elements of $K_0$ and their neighbors, see Fig. \ref{FIG:Stencil}.
Then, we denote the DG solutions in these cells as $u_h(x,y)=p_j^k(x,y)$ 
and the corresponding cell averages as $\bar{p}_j=\frac{1}{|K_j|}\int_{K_j} p_j^k(x,y) dxdy$, where\\
$j\in\{0,1,11,12,2,21,22,3,31,32\}$.
Based on these notations, the MR limiter for the target cell $K_0$ is described as follows:
\begin{itemize}
  \item [Step 1.] Define three substencils as $S_l:=\{K_0,K_l,K_{l1},K_{l2}\}$, $l=1,2,3$. 
  The smoothness indicator corresponding to $S_l$ is defined as 
  \begin{equation}\label{Eq:IS_l}
    IS_l=\max(|\bar{p}_0-\bar{p}_l|,|\bar{p}_0-\bar{p}_{l1}|,|\bar{p}_0-\bar{p}_{l2}|), \quad l=1,2,3.
  \end{equation}
  Based on $IS_l$, we calculate a baseline smoothness indicator as
 \begin{equation}\label{Eq:IS_0}
    IS^0=\min(IS_1,IS_2,IS_3).
  \end{equation}
  For the simplicity of the boundary implementation, we do not consider the cells outside the computational domain when $K_0$ is close to a boundary.
  For example, if $S_3$ is outside the computational domain, we just let $IS^0=\min(IS_1,IS_2)$.
  \item [Step 2.] Calculate the smoothness indicator of the DG solution $p_0^k(x,y)$ in the target cell as 
  \begin{equation}\label{Eq:ISH}
    IS^k=|K_0|^{\frac{k}{D}}\sum_{l=0}^k \frac{1}{(k-l)!l!}\left|\frac{\partial^kp_0^k(x,y)}{\partial x^{k-l}\partial y^{l}}\right|,
  \end{equation}
  which measures the scaled $k$th-order derivative of the solution $p_0^k(x,y)$. 
The factor $|K_0|^{\frac{k}{D}}$ is used to eliminate the influence of mesh size,
where $|K_0|$ is the area of the element and $D=2$ is the spatial dimension of the problem. 
\item[Step 3.] We compare $IS^k$ with the baseline smoothness indicator $IS^0$.
\begin{itemize}
  \item [Step 3.1.] If $IS^k\le C_k\times IS^0$, we directly use $p_0^k(x,y)$ as the final solution. 
  Here, $C_k$ is a constant. 
  \item [Step 3.2.] If $IS^k>C_k\times IS^0$, we set the coefficients of all the basis with nonzero $k$th-order derivatives to be 0, 
  and then reduce the order of the polynomial $k$ by one. 
  \begin{itemize}
    \item [$\bullet$] If $k\ge 1$, we return to the beginning of Step 2.
    \item [$\bullet$] If $k=0$, we reconstruct a first order polynomial $\tilde{p}_0^1(x,y)=a_0+a_1x+a_2y$ with
  \begin{equation}
    a_0=\bar{p}_0,\quad a_l=\text{minmod}(a_{l,1},a_{l,2},a_{l,3}), l=1,2,
  \end{equation}
  where the classical minmod function is defined as
   \begin{equation}
    \text{minmod}(\xi_1,...,\xi_n)=\begin{cases}
      \text{sgn}(\xi_1)\min(|\xi_1|,...,|\xi_n|),\text{if } \text{sgn}(\xi_1)=...=\text{sgn}(\xi_n),\\
      0,\quad \text{otherwise}.
    \end{cases}
   \end{equation}
  The slopes $\{a_{1,1},a_{2,1}\}$, $\{a_{1,2},a_{2,2}\}$, and $\{a_{1,3},a_{2,3}\}$ 
  are respectively reconstructed based on the stencils 
  $S_1^*:=\{K_0,K_1,K_2\}$, $S_2^*:=\{K_0,K_2,K_3\}$, and $S_1^*:=\{K_0,K_3,K_1\}$. 
  More specifically, by enforcing
  \begin{equation}
    \frac{1}{|K_j|}\int_{K_j} \tilde{p}_0^1(x,y) dxdy=\bar{p}_j,\quad K_j\in S_l^*,
  \end{equation}
  we can trivially obtain $\{a_{1,l},a_{2,l}\}$ ($l=1,2,3$).
  \end{itemize}
\end{itemize}
\end{itemize} 

\begin{remark}
  Although the stencil for the baseline indicator is not compact, the reconstruction stencil retains the compactness of the DG method.
  The MR-limiter contains a hyperparameter $C_k$, but plenty of numerical examples 
  demonstrate that the limiter performs very well for different problems 
  and different types of meshes with a fixed $C_k$.
\end{remark}
\begin{remark}
  The above limiter can be naturally degenerated to 1D 
and can be straightforwardly extended to arbitrary polygonal meshes in 2D
and polyhedral meshes in 3D.
  For the 1D case, the baseline smoothness indicator for $K_j$ is calculated as
  $IS^0=\min(IS_L,IS_R)$ with $IS_L=\max(|\bar{p}_j-\bar{p}_{j-1}|,|\bar{p}_{j}-\bar{p}_{j-2}|)$
  and $IS_R=\max(|\bar{p}_j-\bar{p}_{j+1}|,|\bar{p}_{j}-\bar{p}_{j+2}|)$,
  and the smoothness indicator of the DG solution $p_j^k(x)$ in the target cell
  is calculated as $IS^k=\frac{h^k_j}{k!}\left|\frac{\partial^kp_j^k(x)}{\partial x^k}\right|$
  where $h_j$ is the mesh size of $K_j$.
  If $IS^k>C_k\times IS^0$ for all $k\ge 1$, the minmod slope limiter is used 
  to reconstruct the solution as $\tilde{p}_j^1=a_0+a_1x$ 
  where $a_0=\bar{p}_j$ and 
  $a_1=\text{minmod}\left(\frac{\bar{p}_j-\bar{p}_{j-1}}{x_j-x_{j-1}},\frac{\bar{p}_{j+1}-\bar{p}_{j}}{x_{j+1}-x_{j}}\right)$.
  Here, $x_j$ is the cell center of $K_j$.
\end{remark}
\subsection{The properties of the multi-resolution limiter}\label{SUBSEC:MR_limiter_property}

Firstly, the computational cost of the MR limiter is relatively low,
because the baseline smoothness indicator $IS^0$ depends on 
the difference between several neighboring cell averages 
and $IS^k$ only depends on the derivatives of the local DG polynomial in the target cell, which are easy to calculate.
We note that the indicator only depends on the cell averages of neighboring cells, which are invariant during the entire reconstruction.
It means that reconstructing the solutions in neighboring cells does not change the indicator for $K_0$.
This basic property is lacking in many classical indicators, which may lead to a paradox.
For example, the classical KXRCF indicator \cite{Krivodonova2004KXRCF} is defined by
\begin{equation}\label{Eq:KXRCF_Indicator}
  I_{K_0}^{KXRCF}:=\frac{\left|\int_{\partial K_0^-}(p_0(x,y)-p_l(x,y))ds\right|}{h^{\frac{k+1}{2}}|\partial K_0^-|\|p_0(x,y)\| },
\end{equation}
where $\partial K_0^-$ denotes the inflow part of the boundaries of $K_0$,
and $p_l(x,y)$ are the DG polynomials in the neighboring cells of $K_0$ sharing the edge(s) in $\partial K_0^-$.
 As we see, $I_{K_0}^{KXRCF}$ measures the jump of the solution at inflow boundaries,
so the indicator changes accordingly
when the solutions in the neighboring cells are reconstructed. 
Then a `good' cell may become a troubled one after implementing the limiter for reconstructions.
The same paradox might appear in the FS indicator \cite{Fu2017NewLimiter} which is defined as follows
\begin{equation}\label{Eq:FS_Indicator}
  I_{K_0}^{FS}:=\frac{\sum_{j=1}^{3}\left|\frac{1}{|K_0|}\int_{K_0}p_j(x,y)dxdy-\bar{p}_0\right|}{\max_{j\in{0,1,2,3}}\{|\bar{p}_j|\}}.
\end{equation}

Secondly, the MR limiter can detect the discontinuities of the solution and its derivatives.
On one hand, $IS_l$ is defined by the difference of cell averages between neighboring cells,
so we have $IS_l=\mathcal{O}(h)$ when the solution is smooth in the sub-stencil $S_l$ $(l=1,2,3)$
where $h$ is the mesh size. 
When the solution is discontinuous in the substencil $S_l$ $(l=1,2,3)$, we have $IS_l=\mathcal{O}(1)$. 
However, the baseline smoothness indicator is defined as $IS^0=\min(IS_1,IS_2,IS_3)$,
so we always have $IS^0=\mathcal{O}(h)$ if 
the solution is smooth at least in one of the three sub-stencils, which is true for most problems.
On the other hand, when the $k$th-order derivative of $u$ is smooth in $K_0$, $IS^k=\mathcal{O}(h^k)\le C_k\times IS^0$ 
if $h$ is smaller than a certain value. 
However, when the $k$th-order derivative of $u$ is not smooth in $K_0$,
we have $IS^k=\mathcal{O}(1)>C_k\times IS^0$. 
Therefore, the MR limiter described in Section \ref{SUBSEC:MR_limiter}
can effectively detect the discontinuities of the solution and its derivatives by simply comparing $IS^k$ with $C_k\times IS^0$.
When the solution is sufficiently smooth, the limiter directly uses the high-order DG polynomial to approximate the solution.
When the solution itself is discontinuous, the limiter uses a shock-capturing slope limiter 
to suppress numerical oscillations.
When the solution is weakly discontinuous (its derivatives are not smooth),
the limiter uses a truncated DG polynomial to approximate the solution.
In short, the limiter can adaptively choose a proper polynomial to approximate solutions 
with different smooth levels.

Finally, we show that the proposed limiter has the scale-invariant property.
If we perform an affine transformation on the solution $u=p(x,y)$ as 
$\hat{u}=\lambda u+\sigma$, where $\lambda$, $\sigma$ are constants and $\lambda\neq 0$.
From (\ref{Eq:IS_l}), (\ref{Eq:IS_0}), and (\ref{Eq:ISH}), it is straightforward 
to obtain the following equations
\begin{equation}
  IS^0(\lambda u+\sigma)=|\lambda| IS^0(u),\quad IS^k(\lambda u+\sigma)=|\lambda| IS^k(u),
\end{equation}
which indicates that the inequality relationships between $IS^0$ and $C_k\times IS^k$ are invariant
when the scale of the solution changes.
As a result, the MR limiter described in Section \ref{SUBSEC:MR_limiter} 
is not affected by the scale of the solution.

\begin{remark}
  The scale-invariant property is important for the generalization performance of the limiter.
  The OEDG with the damping technique has this property \cite{Peng2025OEDG,Ding2025OEDG}.
  However, many widely used limiters do not possess this property.
  For example, from Eqs. (\ref{Eq:KXRCF_Indicator}) and (\ref{Eq:FS_Indicator}),
  it is trivial to conclude that 
  $I_{K_0}^{KXRCF}(u)\neq I_{K_0}^{KXRCF}(\lambda u+\sigma)$ and $I_{K_0}^{FS}(u)\neq I_{K_0}^{FS}(\lambda u+\sigma)$ if $\sigma\neq 0$,
  because the term $\|p_0(x,y)\|$ in the denominator of $I_{K_0}^{KXRCF}$ 
 and the term $\max_{j\in{0,1,2,3}}\{|\bar{p}_j|\}$ in the denominator of $I_{K_0}^{FS}$ depend on the value of $\sigma$.
 In addition, the term $h^{\frac{k+1}{2}}$ in the denominator of $I_{K_0}^{KXRCF}$
 makes the KXRCF indicator also depend on the scale of the mesh size.
\end{remark}

\section{Numerical examples}\label{SEC:_NumExam}
In this section, we use some 1D and 2D benchmarks to test the performance of 
the RKDG method coupled with the proposed MR limiter
which is applied after every RK stage.
For the compressible Euler equations, the density is used as the indicator variable
and the minmod slope limiter is applied to characteristic variables.
For 2D cases, both quadrilateral and triangular cells are included.
To further improve the efficiency, we adopt the rotated characteristic decomposition technique \cite{Shen2020RotatedDecompositon} 
for the slope limiter, which requires only one-time characteristic decomposition in the direction of the gradient of the density.
In addition, the positivity-preserving limiter proposed 
by Zhang and Shu \cite{Zhang2011MaximumReview,Zhang2012MaximumDGTri} is adopted 
to preserve the positivity of the density and the pressure for the compressible Euler equations.
The CFL numbers are respectively set to 0.3, 0.15, 0.1, 0.06, 0.05, and 0.04 for 
$k=1,2,3,4,5,$ and $6$ except for the convergence tests. 
The hyperparameter $C_k$ is set to 3 for all the simulations except the first test case 
which is used to assess the effect of $C_k$.
All the numerical results are plotted with only cell averages.
\subsection{Numerical examples for the 1D linear advection equation}
We consider the 1D linear advection equation
\begin{equation}\label{Eq:1DAdvEq}
 \frac{\partial u}{\partial t}+\frac{\partial u}{\partial x}=0.
\end{equation}

\begin{figure}[htbp]
 \centering
 \includegraphics[width=8 cm]{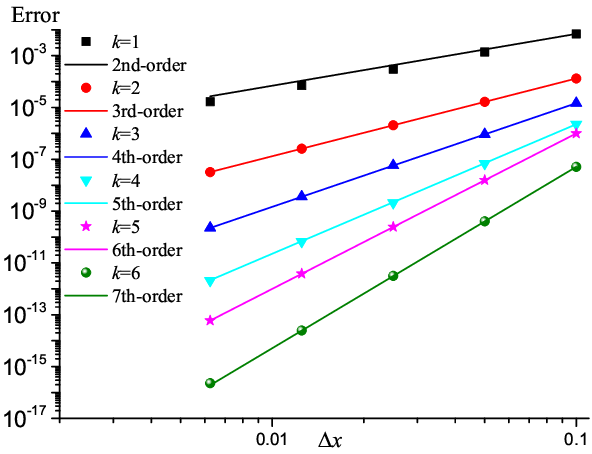}
 \caption{The maximum numerical errors 
at the cell centers versus the mesh size at $t=2$.}
\label{FIG:Linear_convergence}
\end{figure}

\paragraph{Example 4.1.1} We test the convergence rate by setting 
the initial condition as $u_0(x)=\sin(\pi x)$, $x\in[-1,1]$.
Periodic boundary conditions are applied to the left and right boundaries.
For the cases of $k=1$ and $2$, we respectively set $\Delta t=0.3h$ and $\Delta t=0.15h$.
For the cases of $k\ge 3$, the spatial order of the DG operator does not match with
the temporal order of the third-order RK scheme, 
so we respectively set $\Delta t=0.25h^{\frac{4}{3}}$, $\Delta t=0.3h^{\frac{5}{3}}$,
$\Delta t=0.5h^2$, and $\Delta t=0.4h^{\frac{7}{3}}$
for $k=3,4,5,$ and $6$ to guarantee the expected spatial convergence order.
When $C_k\ge 2.91$ and $h\le 0.1$, 
all the RKDG schemes with the MR limiter do not detect any troubled cells 
during the computations of this smooth problem,
so we observe the ideal convergence order for all the cases 
as shown by Fig. \ref{FIG:Linear_convergence}.
For $k=4, 5$ and $6$, the numerical errors already reach the level of truncation errors of double precision,
so we use quadruple precision to eliminate the effect of truncation errors.
As we analyzed in Section \ref{SUBSEC:MR_limiter_property}, $IS^k=\mathcal{O}(h^k)$ and $IS^0=\mathcal{O}(h)$ for smooth solutions,
so the relation $IS^k\le C_k\times IS^0$ can be hold for a much smaller $C_k$ when $k\ge2$.
As a result, the RKDG schemes with the MR limiter can achieve the optimal convergence order 
with a much smaller $C_k$.
To be more specific, the RKDG schemes do not detect any troubled cells on the test meshes 
when we set $C_k\ge1.43$ for $k=2$, $C_k\ge0.09$ for $k=3$, $C_k\ge 0.013$ for $k=4$, $C_k=2.3\times10^{-4}$ for $k=5$,
and $C_k\ge 4\times10^{-5}$ for $k=6$.
We note that the largest mesh size in the convergence tests is 0.1.
The value of $C_k$ can be even smaller to preserve the accuracy when we refine the mesh for $k\ge2$.
Therefore, $C_k=3$ is a safe choice to guarantee the optimal convergence rates for smooth solutions.

\paragraph{Example 4.1.2} We test the accuracy and robustness of the proposed multi-resolution limiter for 
solutions with different smooth levels.
The initial condition is 
a combination of Gaussians, a square wave, a sharp
triangle wave, and a half ellipse arranged from left to right \cite{Jiang_Shu1996WENO}.
Periodic boundary conditions are adopted.

\begin{figure}[htbp]
  \centering
  \subfigure[$k=1,C_k=3$]{
  \includegraphics[width=5.5 cm]{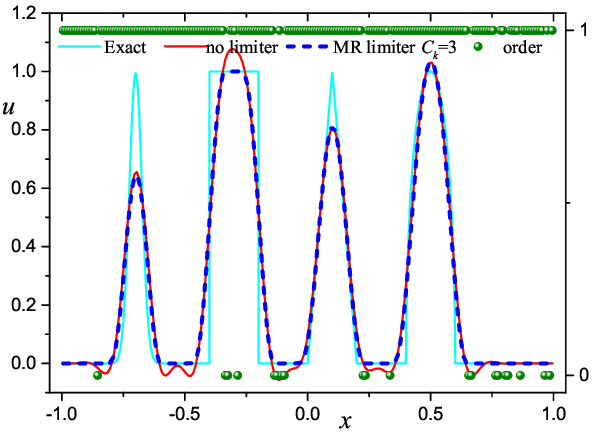}}
  \subfigure[$k=2,C_k=3$]{
  \includegraphics[width=5.5 cm]{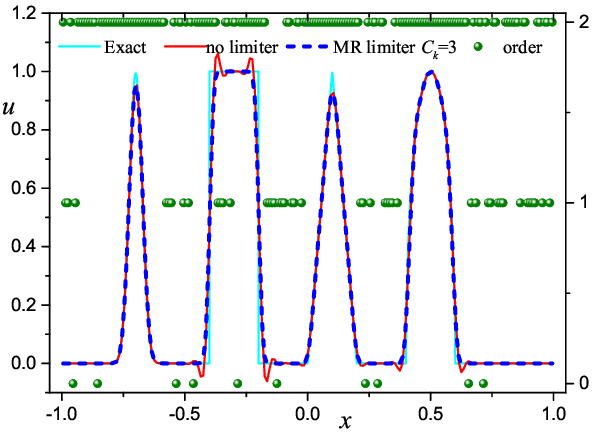}}
  \subfigure[$k=3,C_k=3$]{
  \includegraphics[width=5.5 cm]{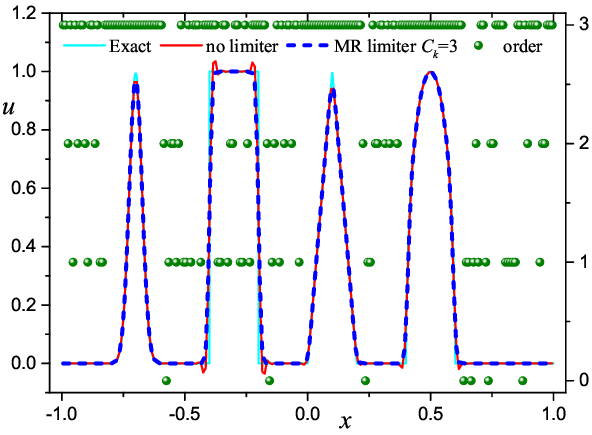}}
  \subfigure[$k=1,C_k=3.5$]{
  \includegraphics[width=5.5 cm]{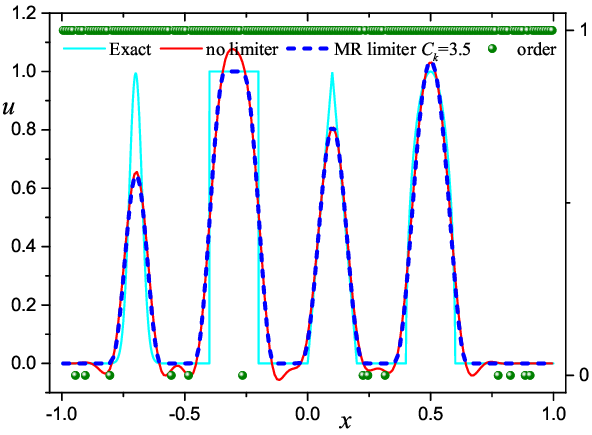}}
  \subfigure[$k=2,C_k=3.5$]{
  \includegraphics[width=5.5 cm]{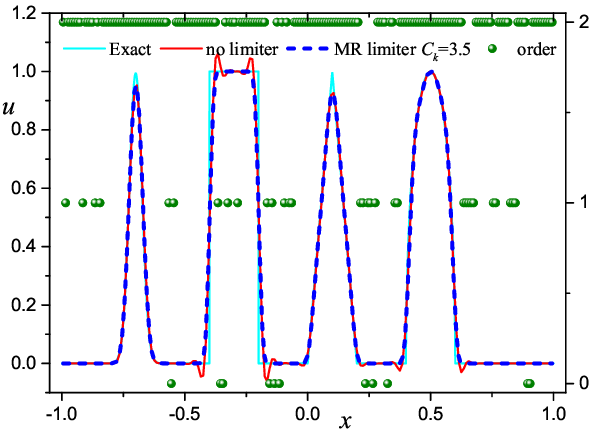}}
  \subfigure[$k=3,C_k=3.5$]{
  \includegraphics[width=5.5 cm]{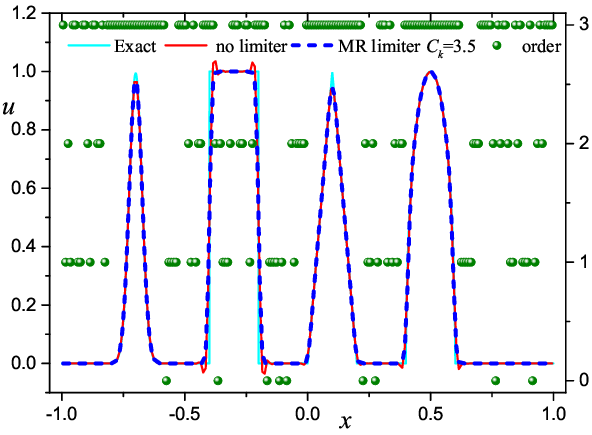}}
  \subfigure[$k=4,C_k=3$]{
  \includegraphics[width=5.5 cm]{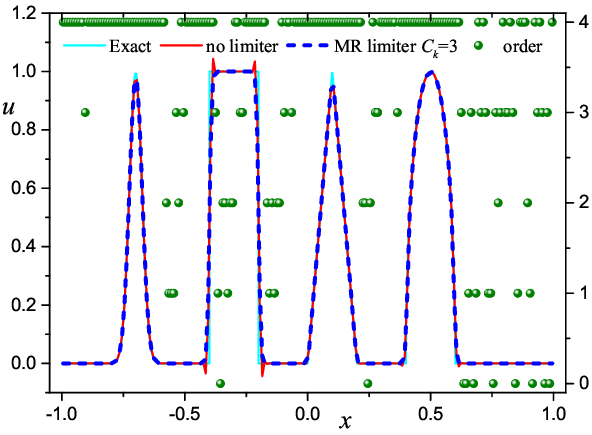}}
  \subfigure[$k=5,C_k=3$]{
  \includegraphics[width=5.5 cm]{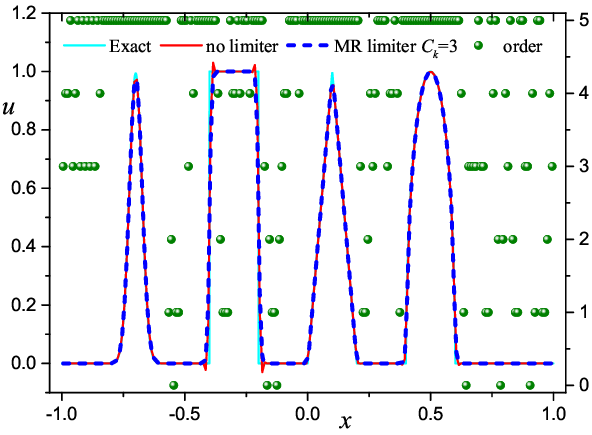}}
  \subfigure[$k=6,C_k=3$]{
  \includegraphics[width=5.5 cm]{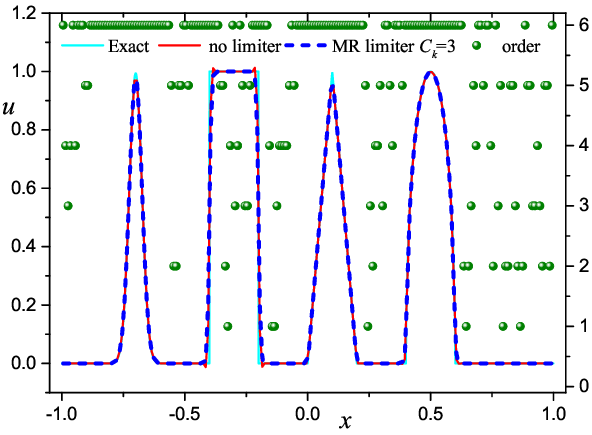}}
  \subfigure[$k=4,C_k=3.5$]{
  \includegraphics[width=5.5 cm]{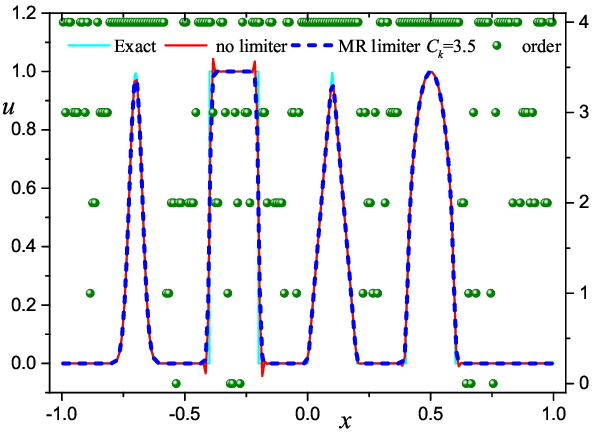}}
  \subfigure[$k=5,C_k=3.5$]{
  \includegraphics[width=5.5 cm]{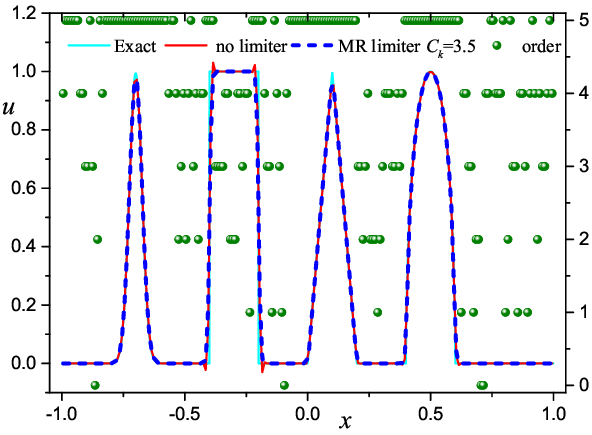}}
  \subfigure[$k=6,C_k=3.5$]{
  \includegraphics[width=5.5 cm]{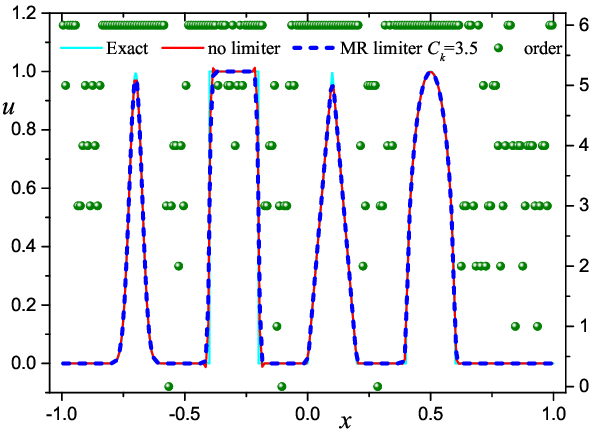}}
  \caption{Distributions of $u$ and polynomial order for 
  linear advection of a combination of Gaussians, a square wave, a sharp
triangle wave, and a half ellipse at $t=20$ computed by RKDG schemes using 200 cells.}
 \label{FIG:Linear_Advection}
 \end{figure}

Fig. \ref{FIG:Linear_Advection} shows the distributions of $u$ and polynomial orders at $t=20$ 
computed by RKDG schemes with and without the MR limiter using $200$ cells. 
Order=$0$ means using the minmod slope limiter to reconstruct the solution in the local cell.
We observe that the original DG schemes with no limiter generate prominent
overshoots and undershoots near discontinuities, especially for $k\le5$.
However, all the RKDG schemes adopting the MR limiter with both $C_k=3$ and $3.5$
can effectively suppress the numerical oscillations near discontinuities while preserving almost the same accuracy
of the original RKDG schemes with no limiter in smooth regions.
In all of the following numerical tests, we uniformly set $C_k=3$.

\begin{figure}[htbp]
  \centering
  \includegraphics[width=8 cm]{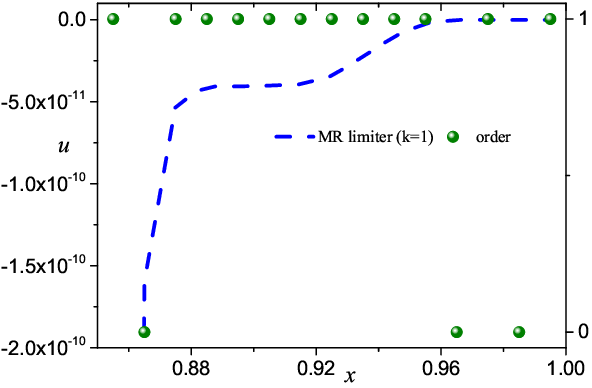}
  \caption{Zoom-in view of the distribution of $u$ for linear advection of a combination of Gaussians, a square wave, a sharp
triangle wave, and a half ellipse with different scale factors at $t=20$ computed by RKDG schemes with the MR limiter using 200 cells ($k=1$).}
 \label{FIG:Linear_Advection_Zoom_In}
 \end{figure}

\begin{remark}
  We note that some cells in the constant solution regions reduce the polynomial order due to the truncation errors of the computer.
  When the solution is constant within the indicator stencil $S$,
  both $IS^0$ and $IS^k$ are theoretically equal to 0.
However, they are actually very small numbers due to the accumulation of truncation errors.
Fig. \ref{FIG:Linear_Advection_Zoom_In} shows the zoom-in view of the distribution of $u$ 
at $t=20$ computed by RKDG schemes with the MR limiter and $k=1$.
It shows that the solution is indeed discontinuous at the scale of $10^{-10}$.
Therefore, $IS^k\le C_k\times IS^0$ is not always guaranteed in `constant' solution regions,
and the MR limiter may select a low-order polynomial.
This precisely demonstrates that the MR limiter is scale-invariant
and can detect discontinuities of very small scales.
Nevertheless, low-order polynomials do not affect the computational accuracy for constant solutions.
We just can not get a nice distribution of troubled cells when there are large portions of constant solutions.
However, we have to emphasize that computational accuracy is more important 
than a nice distribution of troubled cells. 
If we seek nice distributions of troubled cells,
we can set a small threshold to filter out the troubled cells in constant solution regions,
but the limiter will inevitably lose the scale-invariant property.
\end{remark}

\begin{figure}[htbp]
  \centering
  \subfigure[$k=1,\lambda=10^{-10}$]{
  \includegraphics[width=5.5 cm]{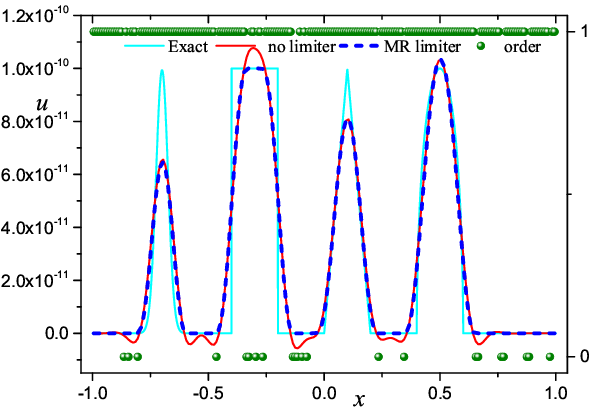}}
  \subfigure[$k=2,\lambda=10^{-8}$]{
  \includegraphics[width=5.5 cm]{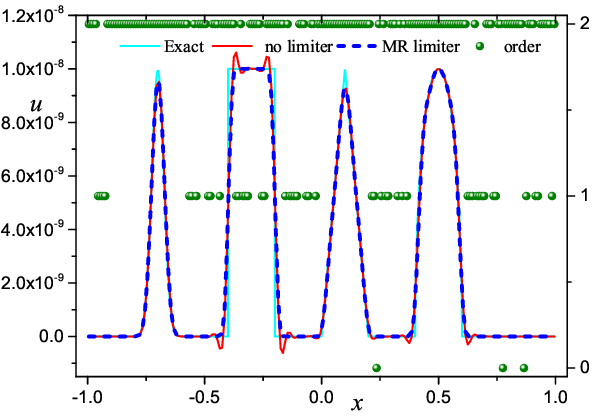}}
  \subfigure[$k=3,\lambda=10^{-6}$]{
  \includegraphics[width=5.5 cm]{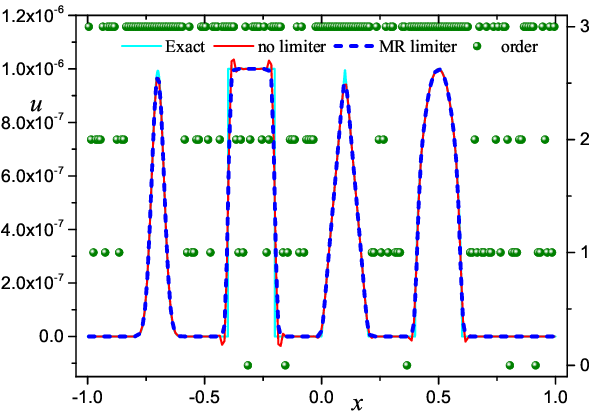}}
  \subfigure[$k=4,\lambda=10^{6}$]{
  \includegraphics[width=5.5 cm]{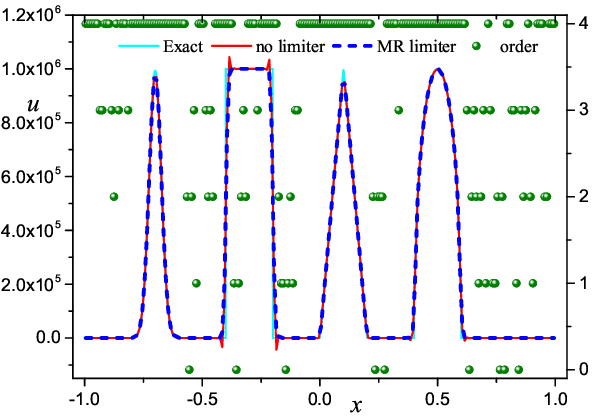}}
  \subfigure[$k=5,\lambda=10^{8}$]{
  \includegraphics[width=5.5 cm]{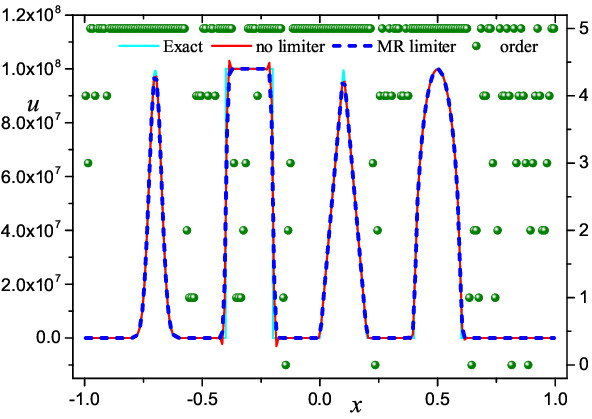}}
  \subfigure[$k=6,\lambda=10^{10}$]{
  \includegraphics[width=5.5 cm]{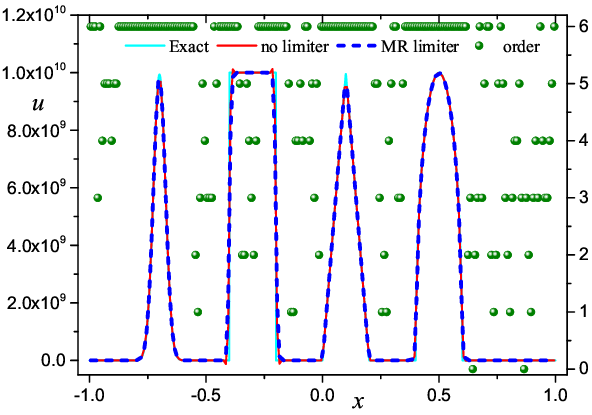}}
  \caption{Distributions of $u$ and polynomial order for linear advection of a combination of Gaussians, a square wave, a sharp
triangle wave, and a half ellipse with different scale factors at $t=20$ computed by RKDG schemes using 200 cells.}
 \label{FIG:Linear_Advection_Scale_Invariant}
 \end{figure}

To further test the scale-invariant property of the MR limiter,
we amplify the solution by a factor of $\lambda\in[10^{-10},10^{10}]$.
Fig. \ref{FIG:Linear_Advection_Scale_Invariant} shows that DG schemes 
with the MR limiter retain the same accuracy and robustness
regardless of the scales of the solutions.

\subsection{Numerical examples for the 1D compressible Euler equations}
The 1D compressible Euler equations read
\begin{equation}
 \frac{\partial \mathbf{U}}{\partial t}+\frac{\partial \mathbf{F}}{\partial x}=0,
\end{equation}
where $\mathbf{U}=[\rho,\rho u,\rho e]^T$ and $\mathbf{F}=[{\rho u,\rho u^2+p,(\rho e+p)u}]^T$
with $\rho$, $u$, $p$, and $e$ respectively representing the density, velocity, pressure, and specific total energy.
The specific total energy is calculated as $e=\frac{p}{\rho(\gamma-1)}+\frac{1}{2} u^2$,
where $\gamma$ is the specific heat ratio.
Without particular instructions, $\gamma$ is set to 1.4.

\begin{figure}[htbp]
  \centering
  \subfigure[$k=1$]{
  \includegraphics[width=5.5 cm]{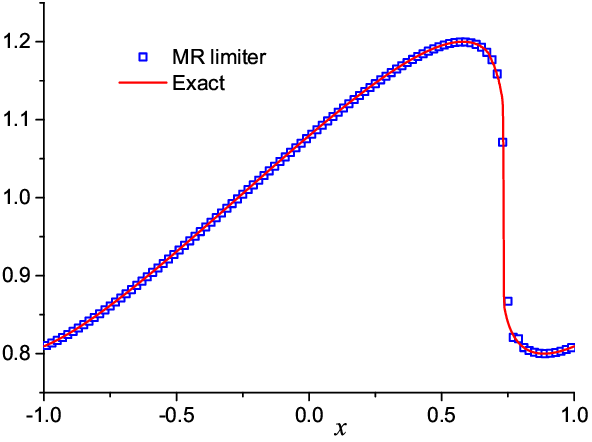}}
  \subfigure[$k=2$]{
  \includegraphics[width=5.5 cm]{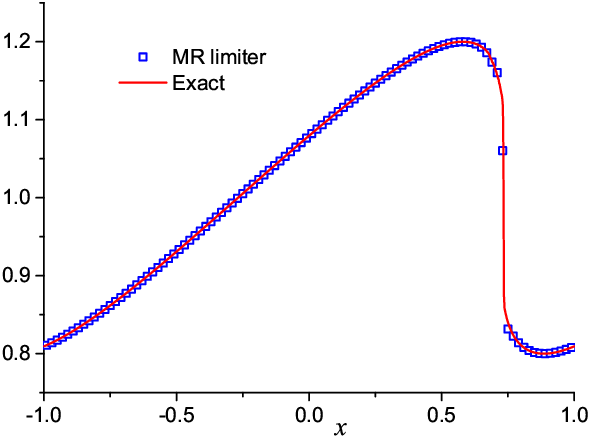}}
  \subfigure[$k=3$]{
  \includegraphics[width=5.5 cm]{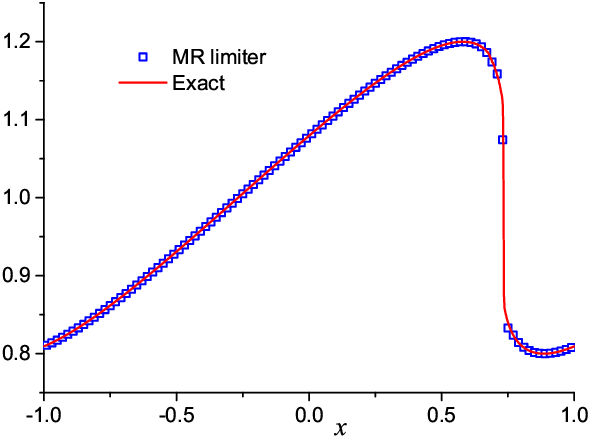}}
  \subfigure[$k=4$]{
  \includegraphics[width=5.5 cm]{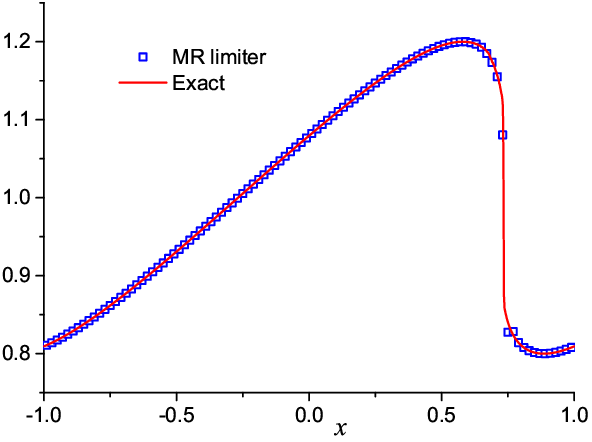}}
  \subfigure[$k=5$]{
  \includegraphics[width=5.5 cm]{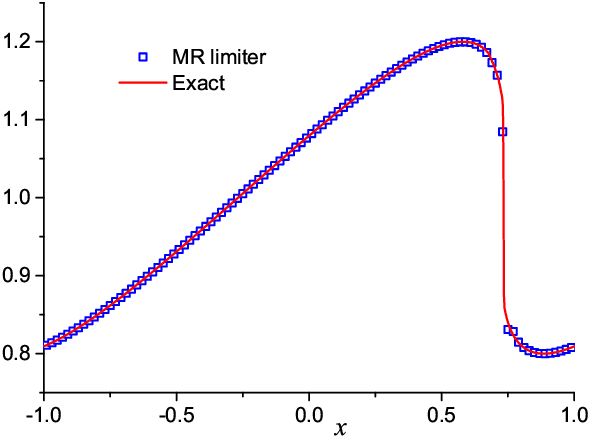}}
  \subfigure[$k=6$]{
  \includegraphics[width=5.5 cm]{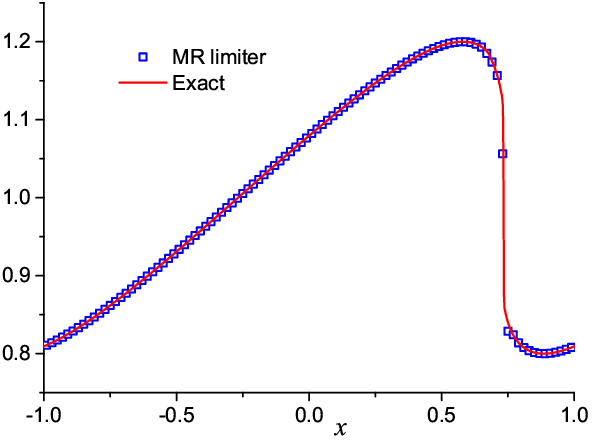}}
  \subfigure[$k=1$]{
  \includegraphics[width=5.5 cm]{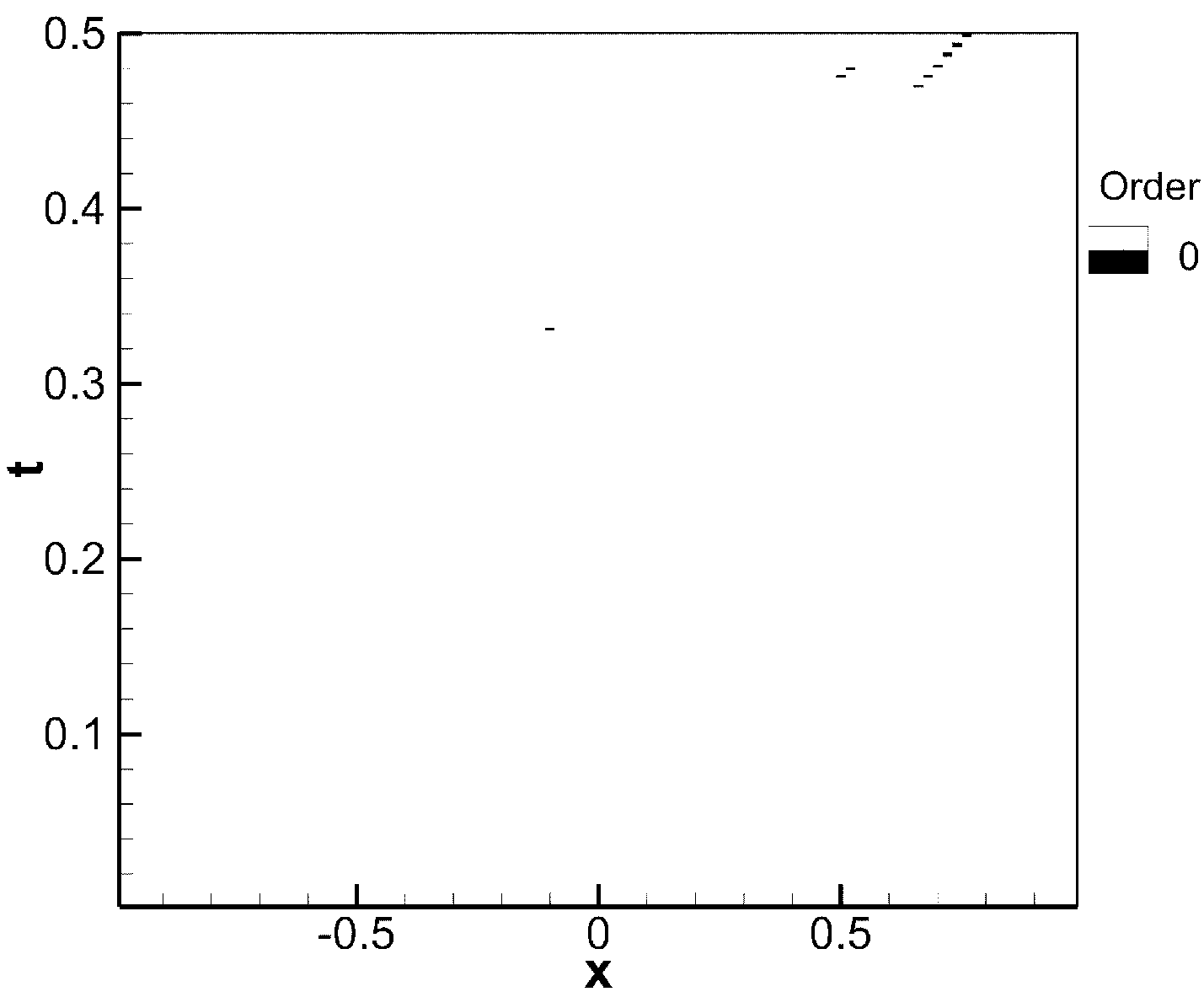}}
  \subfigure[$k=2$]{
  \includegraphics[width=5.5 cm]{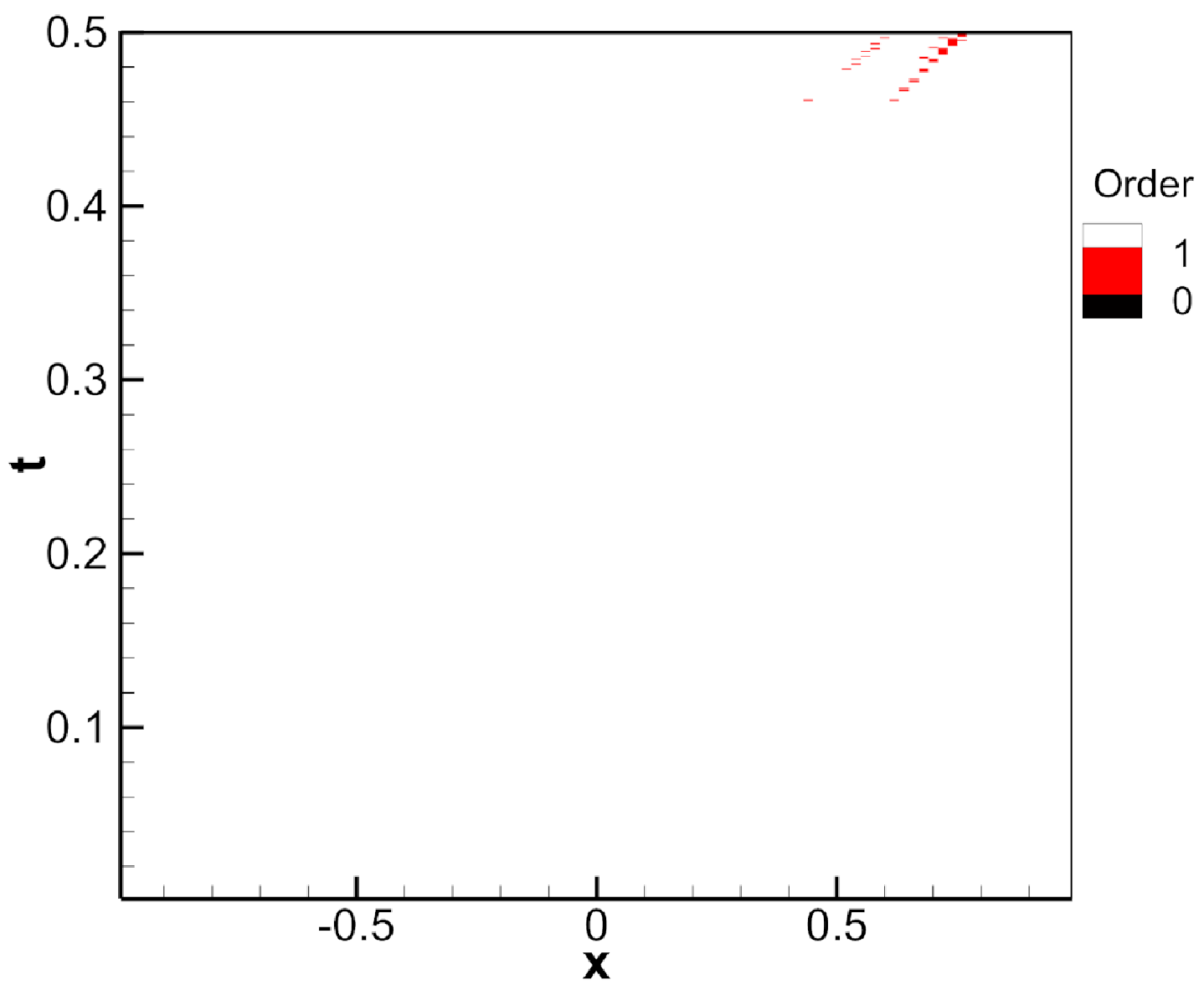}}
  \subfigure[$k=3$]{
  \includegraphics[width=5.5 cm]{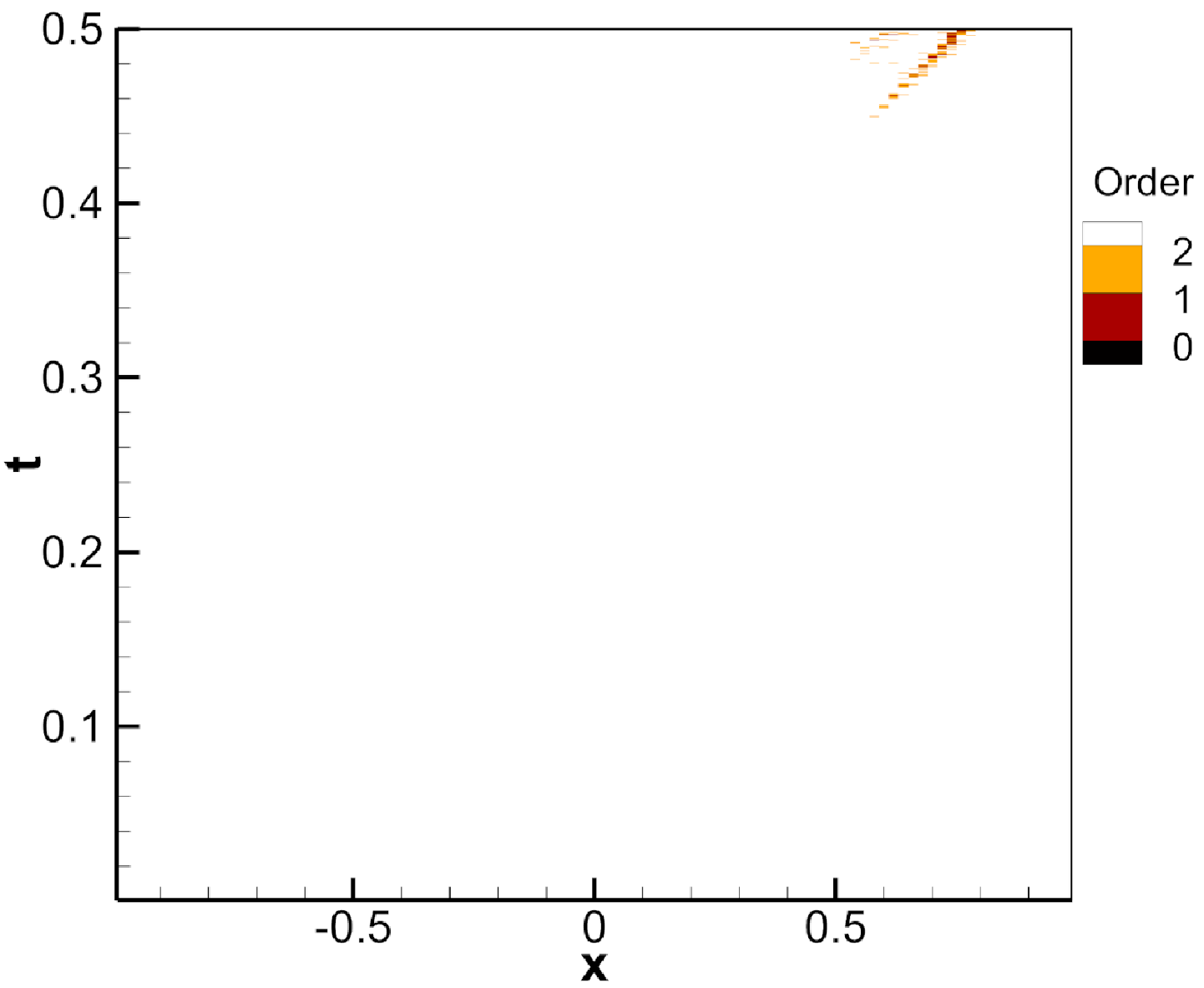}}
  \subfigure[$k=4$]{
  \includegraphics[width=5.5 cm]{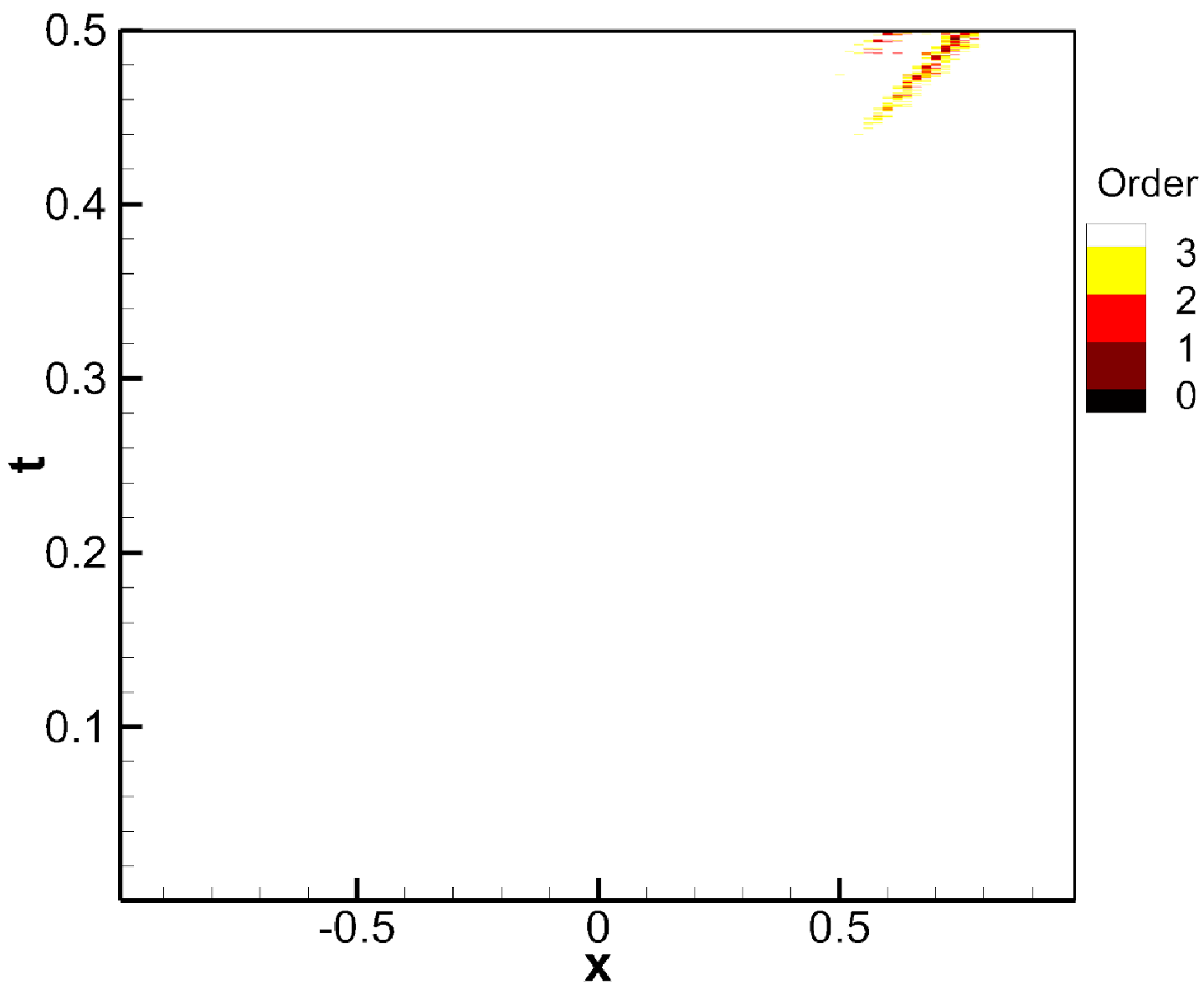}}
  \subfigure[$k=5$]{
  \includegraphics[width=5.5 cm]{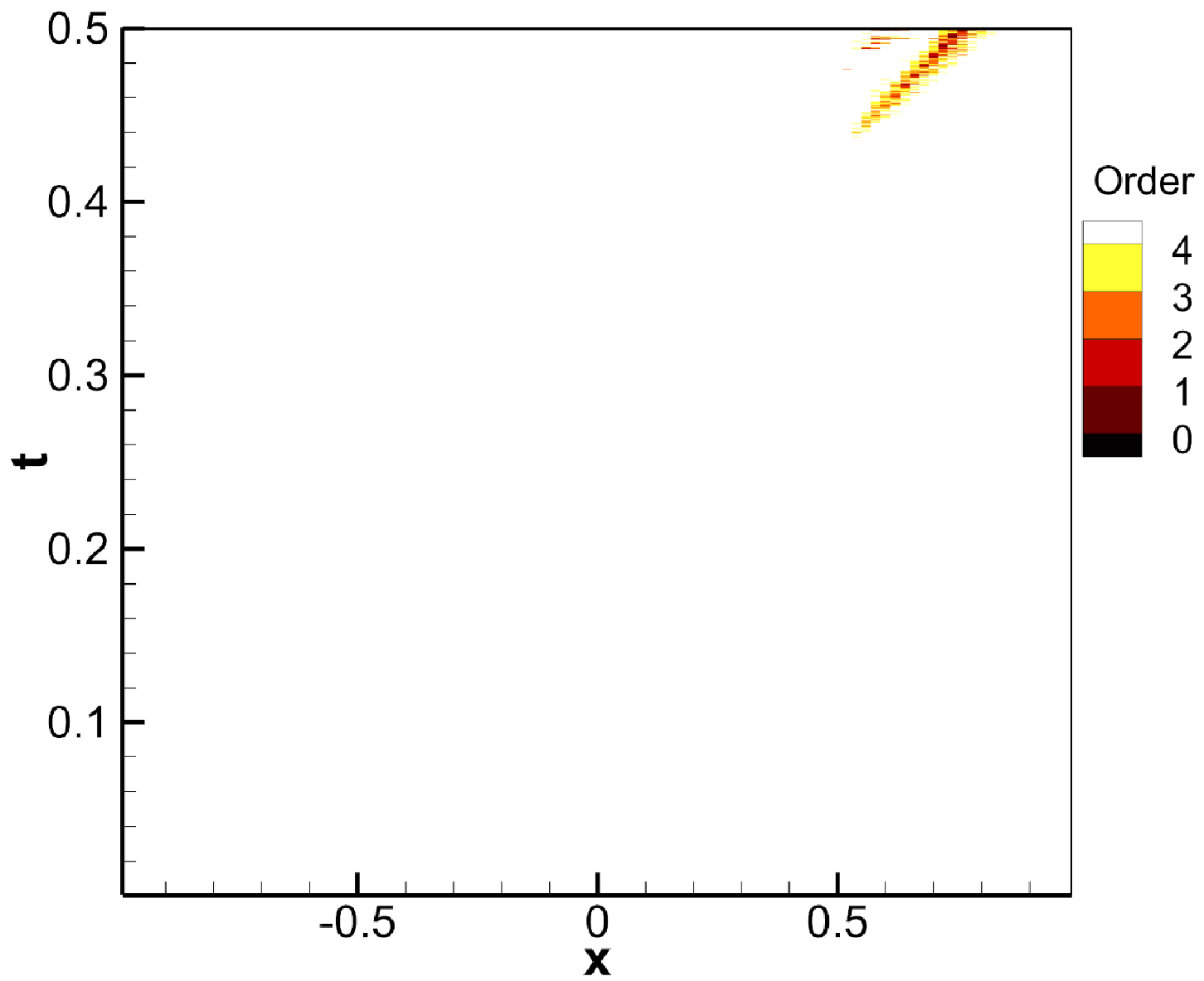}}
  \subfigure[$k=6$]{
  \includegraphics[width=5.5 cm]{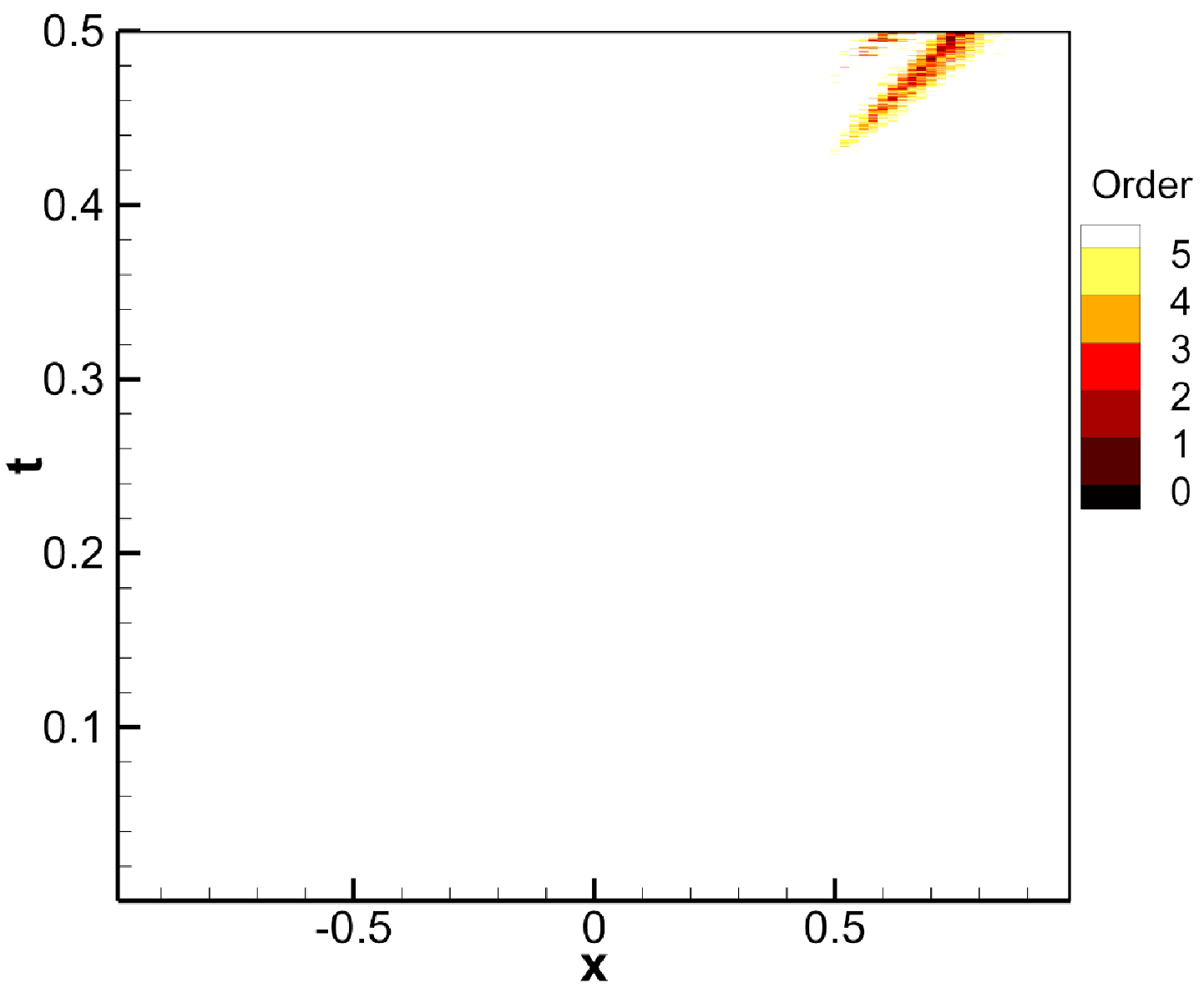}}
  \caption{Density profiles at $t=0.5$ and the time history of the polynomial order
  for the first case of the 1D compressible Euler equations computed by RKDG schemes with the MR limiter using 100 cells. 
  The first and second rows: density profiles; The third and fourth rows: polynomial order, 
  the white parts represent the original $k$th-order DG polynomial.}
 \label{FIG:Burgers}
 \end{figure}

\paragraph{Example 4.2.1} This is a non-physical problem taken from \cite{Fu2017NewLimiter}.
The initial conditions are 
\begin{equation*}
  \rho(x,0)=1+0.2\sin(\pi x),\quad u(x,0)=\sqrt{\gamma}\rho(x,0),\quad p(x,0)=\rho(x,0)^{\gamma},
\end{equation*}
where $\gamma=3$. Periodic boundary conditions are implemented on both sides of the computational domain $[-1,1]$.
The exact solution for $\rho(x,t)$ is the solution to the following scalar Burgers equation
\begin{equation*}
  \frac{\partial\rho}{\partial t}+\frac{\partial(\sqrt{\gamma}\rho^2)}{\partial x}=0,
\end{equation*}
with the initial condition $\rho(x,0)=1+0.2\sin(\pi x)$,
and the exact solutions for the velocity and the pressure are 
$u(x,t)=\sqrt{\gamma}\rho(x,t)$ and $p(x,t)=\rho(x,t)^{\gamma}$.
The solution will form a discontinuity at $t=\frac{5\sqrt{3}}{6\pi}\approx 0.46$ due to the nonlinearity.
We use 100 cells to compute this problem to $t=0.5$.
Fig. \ref{FIG:Burgers} shows the density profiles and the time history of the polynomial order
computed by the RKDG schemes with the MR limiter.
We observe that no troubled cell is detected by the MR limiter before the formation of the discontinuity,
and the polynomial order only decreases near the discontinuities.
In addition, when $k>1$, there are indeed some polynomials falling between the highest order and the lowest order.

\begin{figure}[htbp]
  \centering
  \subfigure[$k=1$]{
  \includegraphics[width=5.5 cm]{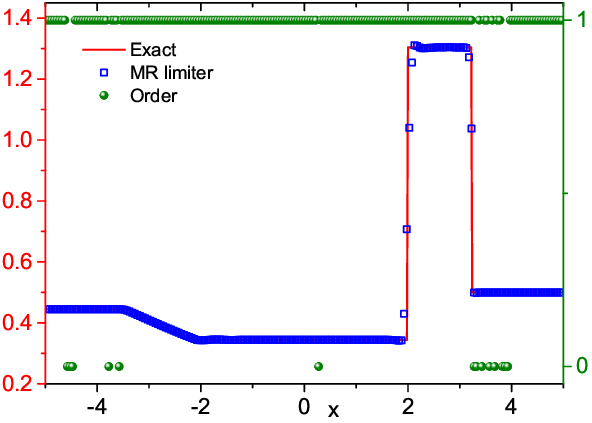}}
  \subfigure[$k=2$]{
  \includegraphics[width=5.5 cm]{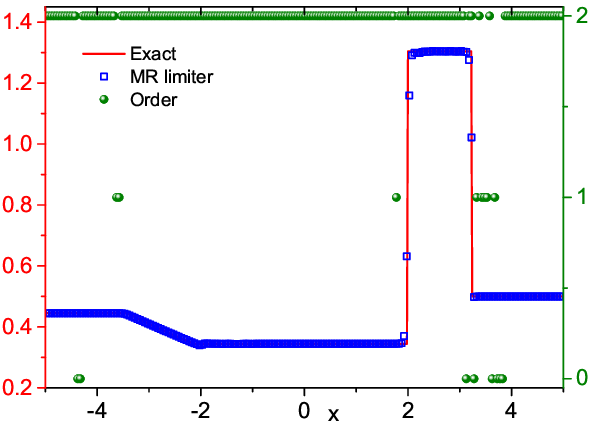}}
  \subfigure[$k=3$]{
  \includegraphics[width=5.5 cm]{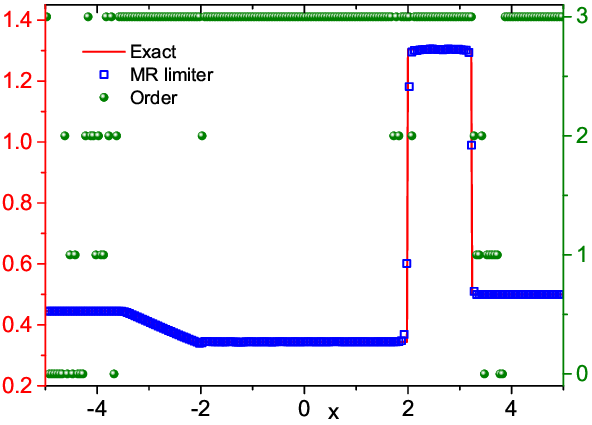}}
  \subfigure[$k=4$]{
  \includegraphics[width=5.5 cm]{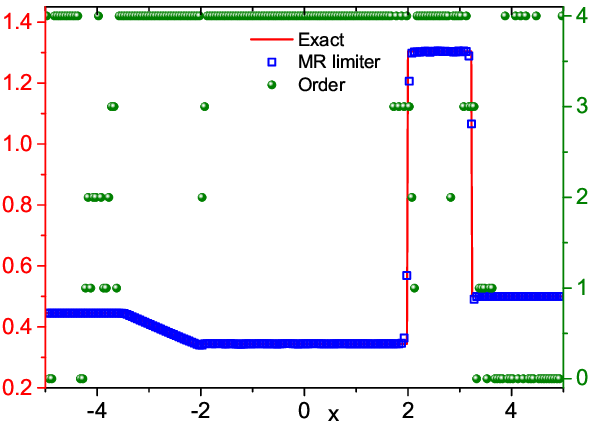}}
  \subfigure[$k=5$]{
  \includegraphics[width=5.5 cm]{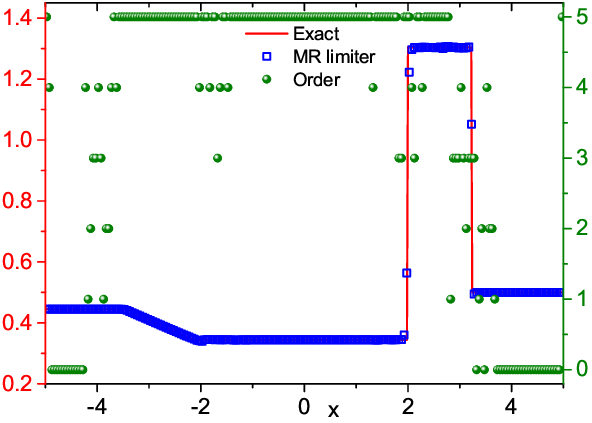}}
  \subfigure[$k=6$]{
  \includegraphics[width=5.5 cm]{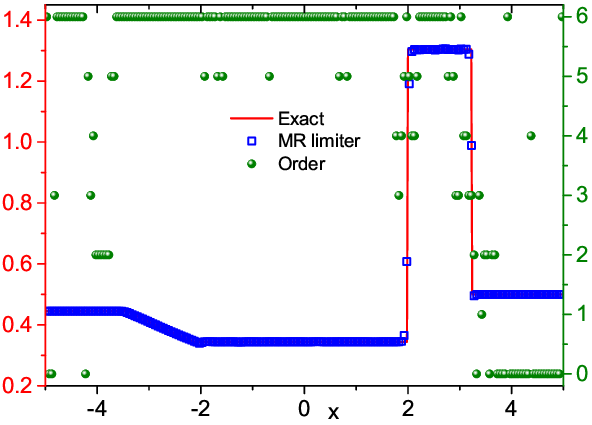}}
  \subfigure[$k=1$]{
  \includegraphics[width=5.5 cm]{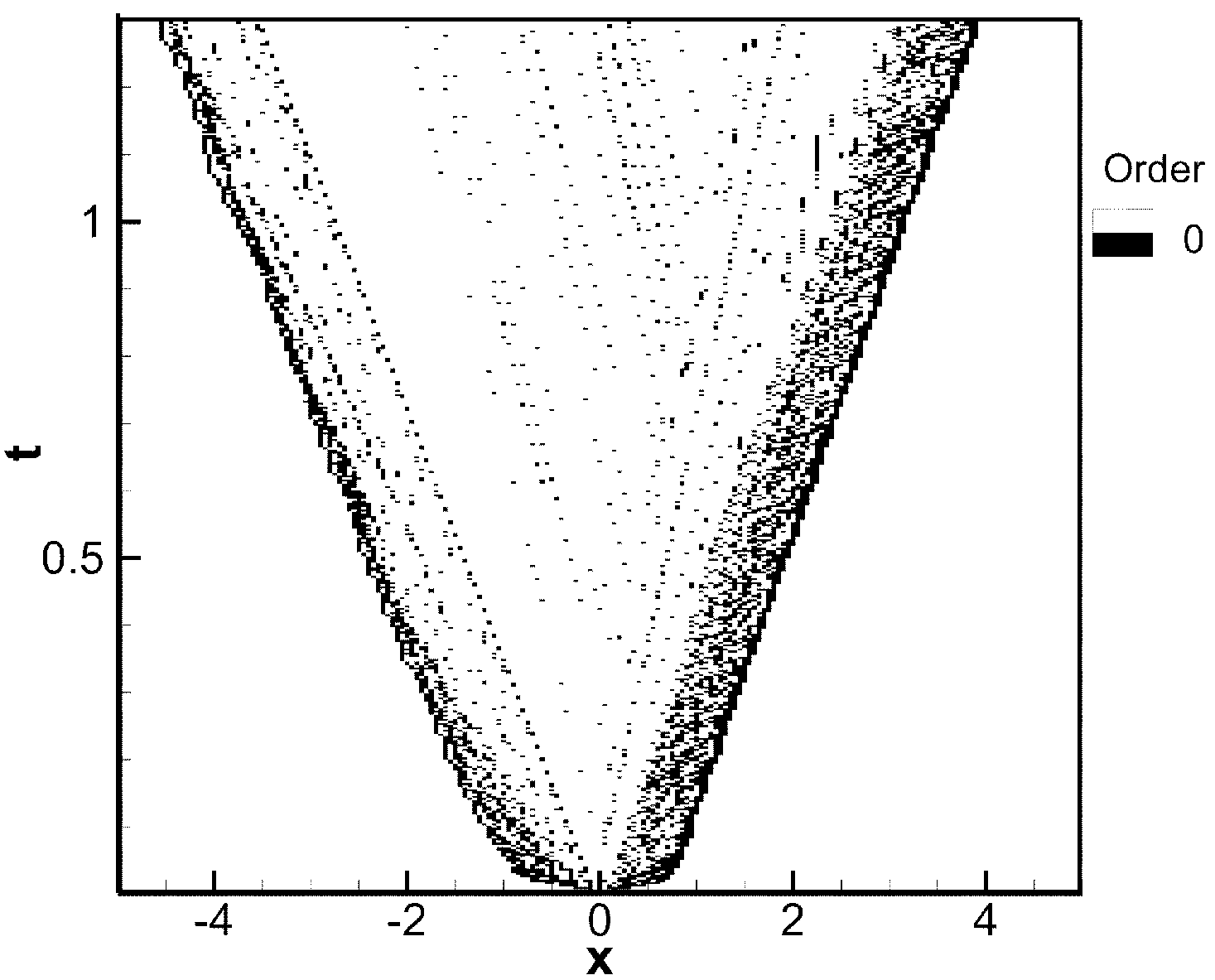}}
  \subfigure[$k=2$]{
  \includegraphics[width=5.5 cm]{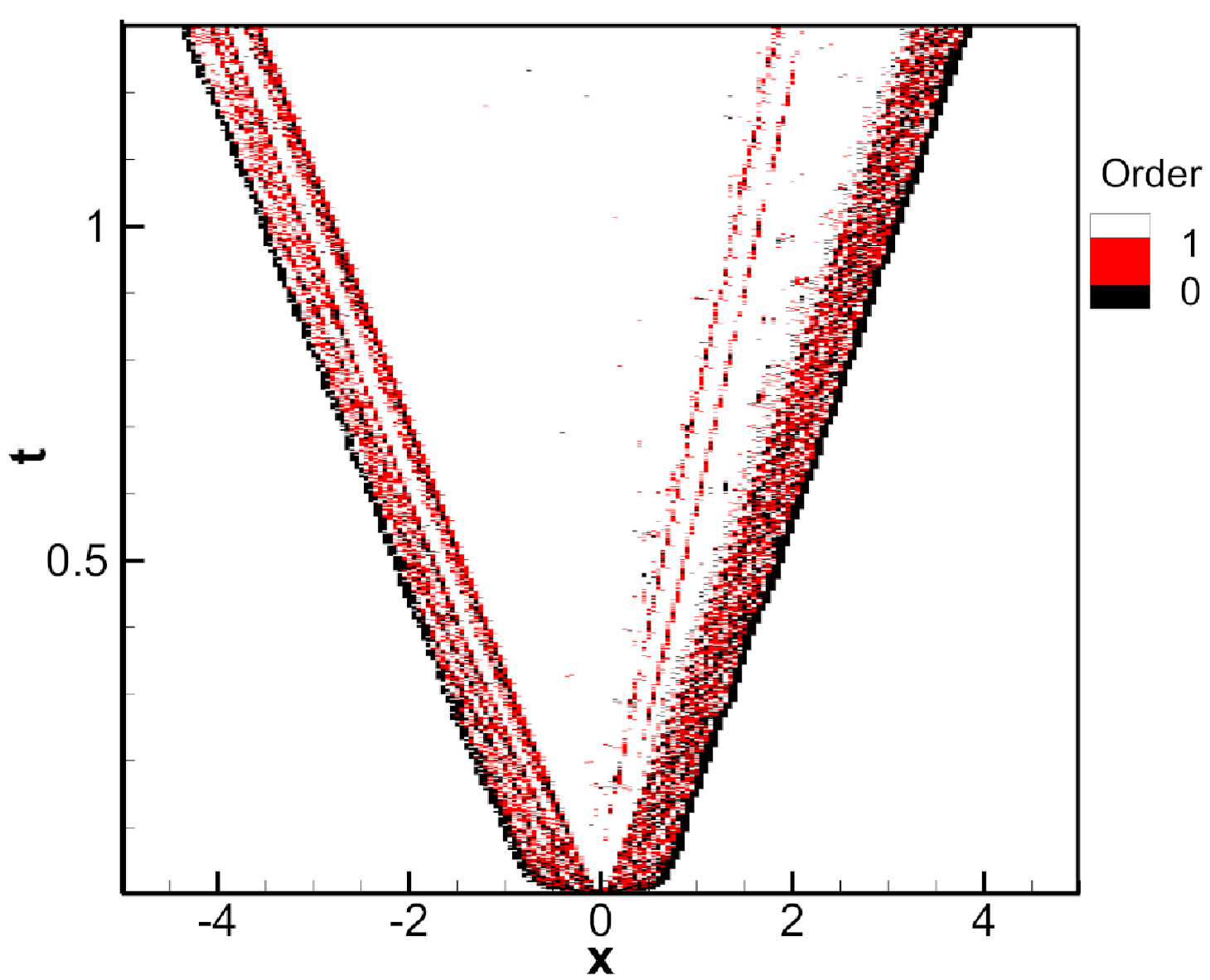}}
  \subfigure[$k=3$]{
  \includegraphics[width=5.5 cm]{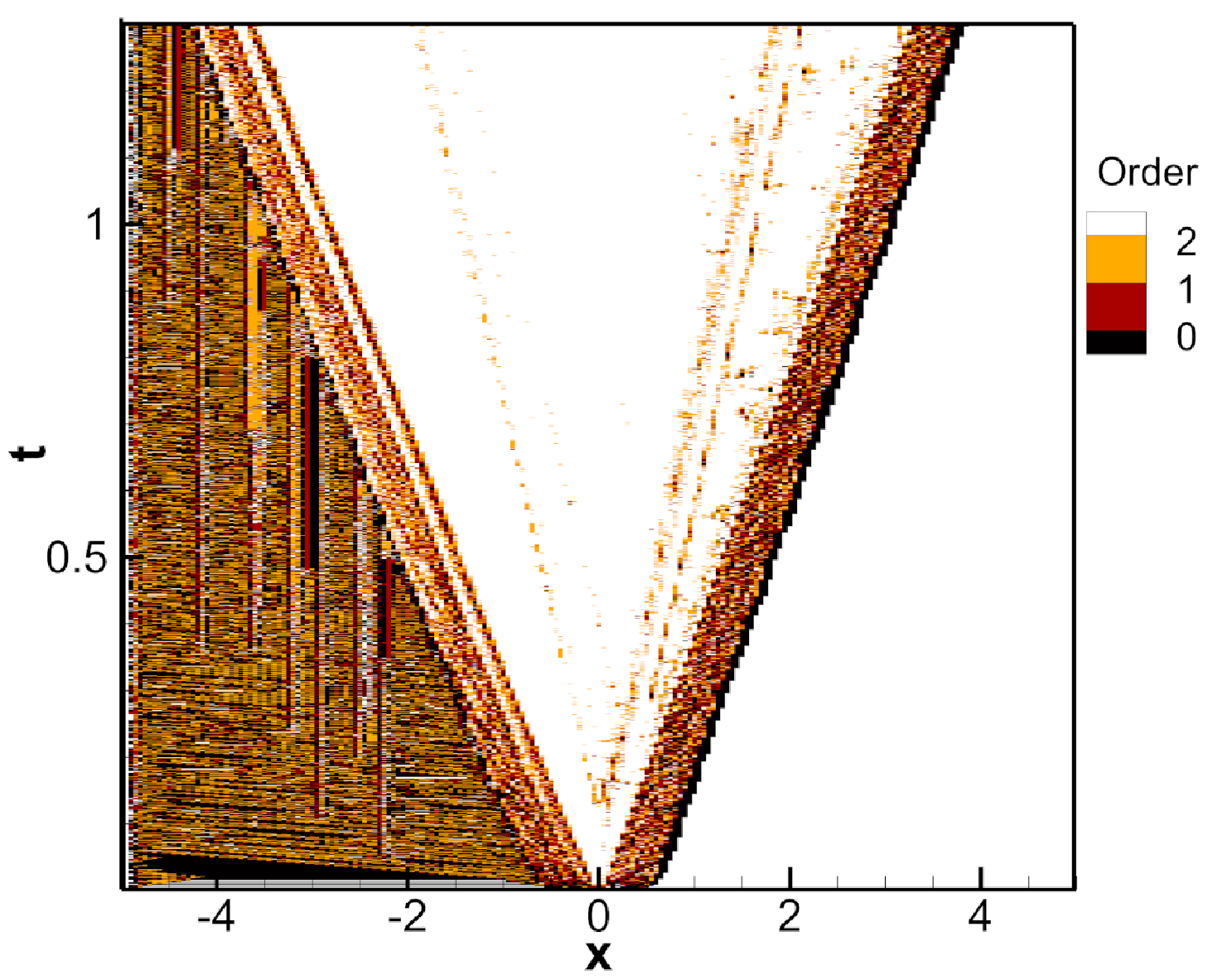}}
  \subfigure[$k=4$]{
  \includegraphics[width=5.5 cm]{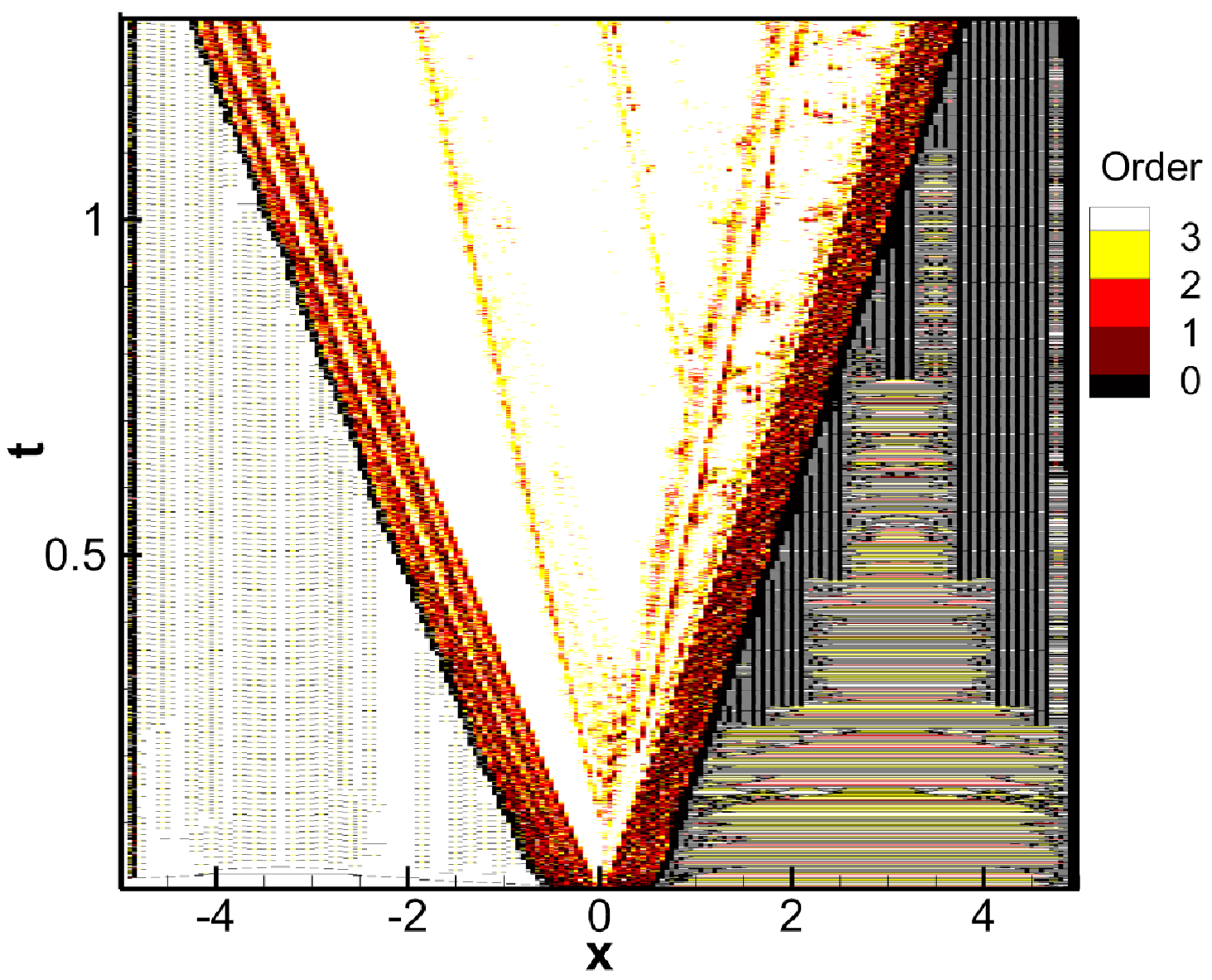}}
  \subfigure[$k=5$]{
  \includegraphics[width=5.5 cm]{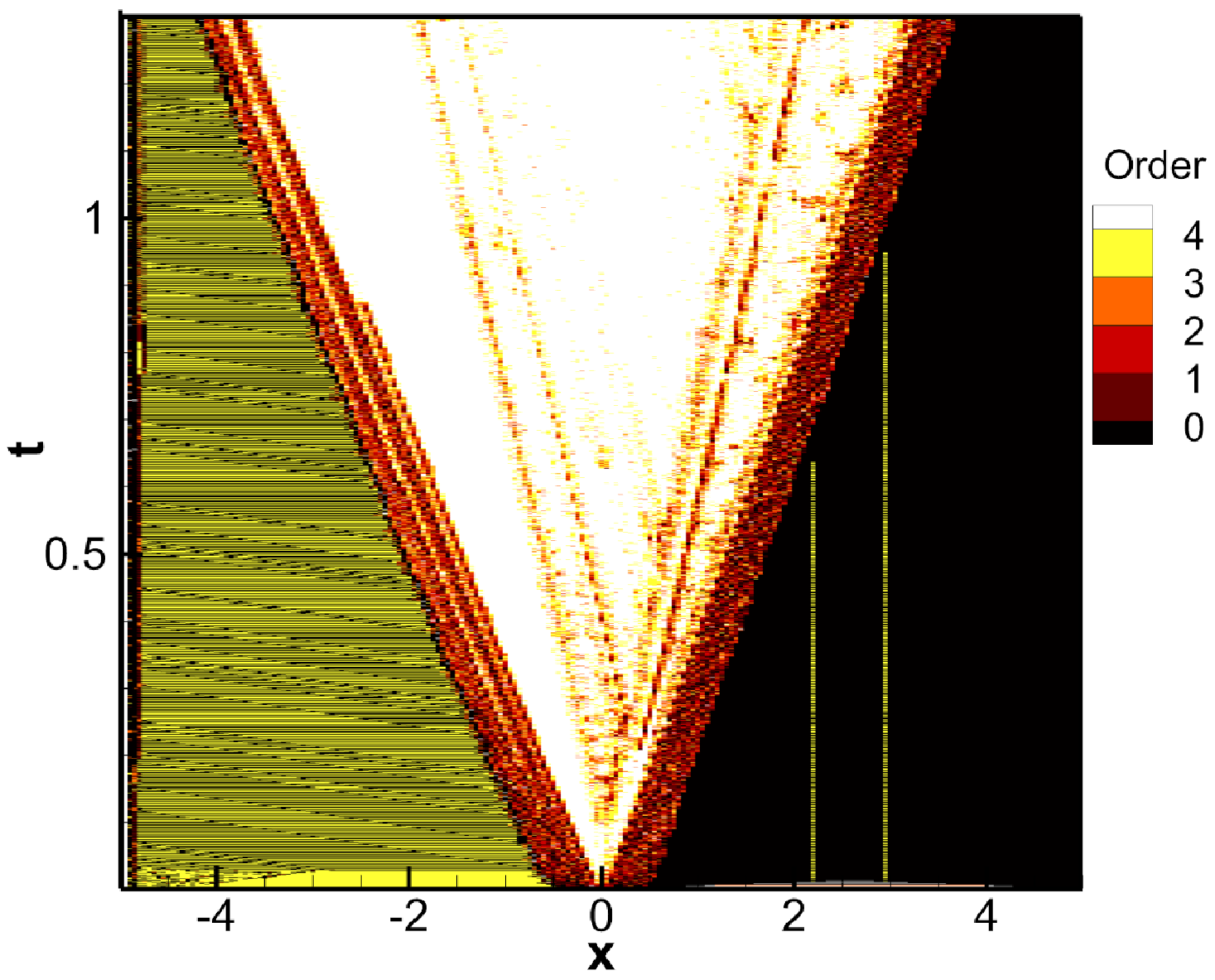}}
  \subfigure[$k=6$]{
  \includegraphics[width=5.5 cm]{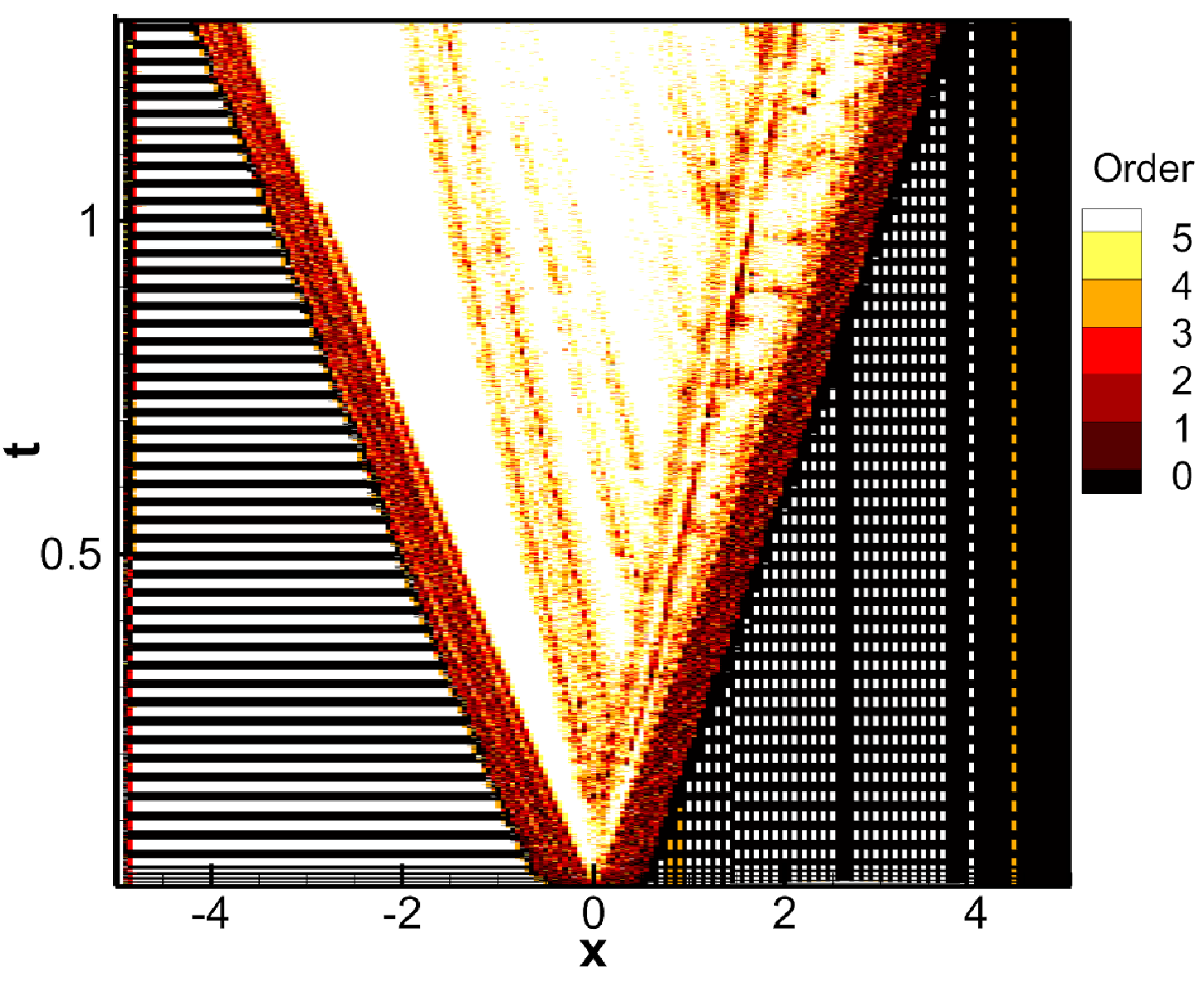}}
  \caption{The Lax shock tube problem computed by RKDG schemes with the MR limiter using 200 cells. 
  The first and second rows: distributions of the density and the polynomial order at $t=1.3$;
   The third and fourth rows: the time history of the polynomial order, 
  the white parts represent the original $k$th-order DG polynomial. }
 \label{FIG:Sod}
 \end{figure}

 \paragraph{Example 4.2.2} The Lax shock tube problem. 
 The computational domain is [-5,5], and the initial condition is given by
\begin{equation*}
 (\rho,u,p)=\begin{cases}
              (0.445,0.698,3.528), & \mbox{if } x<0, \\
              (0.5,0,0.571), & \mbox{otherwise}.
            \end{cases}
\end{equation*}
We compute this problem using RKDG schemes with the MR limiter using 200 cells.
Fig. \ref{FIG:Sod} shows the distributions of the density and the polynomial order at $t=1.3$
and the time history of the polynomial order.
We observe that the order of the polynomial decreases near the shock
and retains the highest order for the expansion fans. 
Some troubled cells are detected in the constant density regions due to truncation errors.

\paragraph{Example 4.2.3} The double rarefaction wave problem.
The computational domain is [-1,1], and the initial condition is given by
\begin{equation*}
 (\rho,u,p)=\begin{cases}
              (7,-1,0.2), & \mbox{if } x<0, \\
              (7,1,0.2), & \mbox{otherwise}.
            \end{cases}
\end{equation*}
When the simulation starts, the gases move fast in the opposite direction,
and a near-vacuum region forms at the center due to the strong rarefaction waves.
Fig. \ref{FIG:Sod} shows the distributions of the density and the polynomial order at $t=0.6$
and the time history of the polynomial order computed by the RKDG schemes with the MR limiter using 200 cells.
The MR limiter retains non-oscillatory while detecting much less troubled cells within the rarefaction fans
than the KXRCF indicator and the FS indicator, as shown by Fig. 3.7 in \cite{Fu2017NewLimiter}.
For the cases of $k>1$, although some troubled cells are detected, the polynomial order does not reduce to 0.
In fact, the flow field of rarefaction waves is continuous,
so it is suitable to approximate it by high-order polynomials.
The troubled cells detected outside the expansion region are caused by truncation errors.

\begin{figure}[htbp]
  \centering
  \subfigure[$k=1$]{
  \includegraphics[width=5.5 cm]{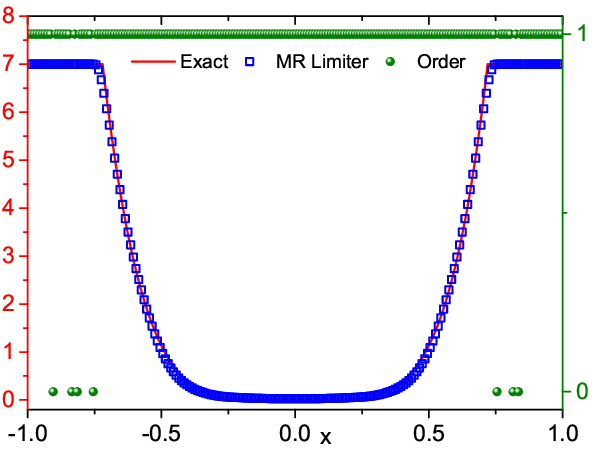}}
  \subfigure[$k=2$]{
  \includegraphics[width=5.5 cm]{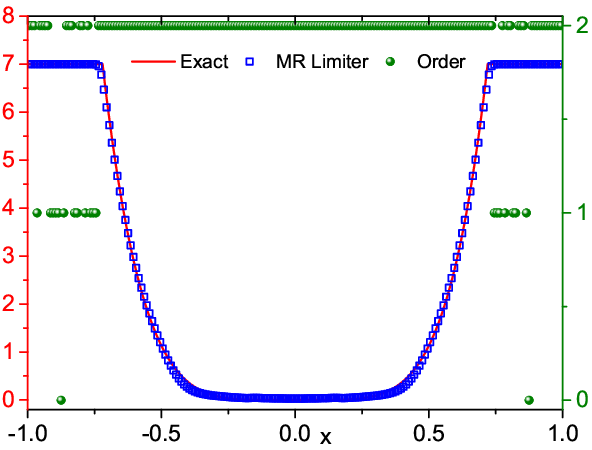}}
  \subfigure[$k=3$]{
  \includegraphics[width=5.5 cm]{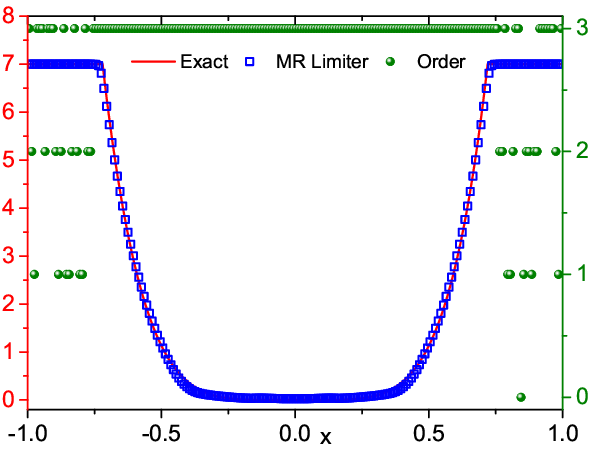}}
  \subfigure[$k=4$]{
  \includegraphics[width=5.5 cm]{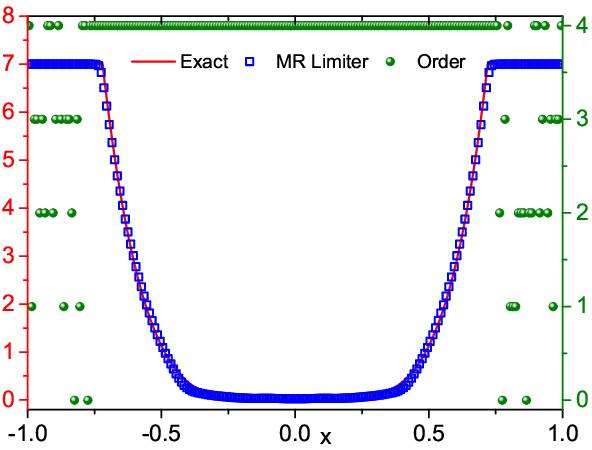}}
  \subfigure[$k=5$]{
  \includegraphics[width=5.5 cm]{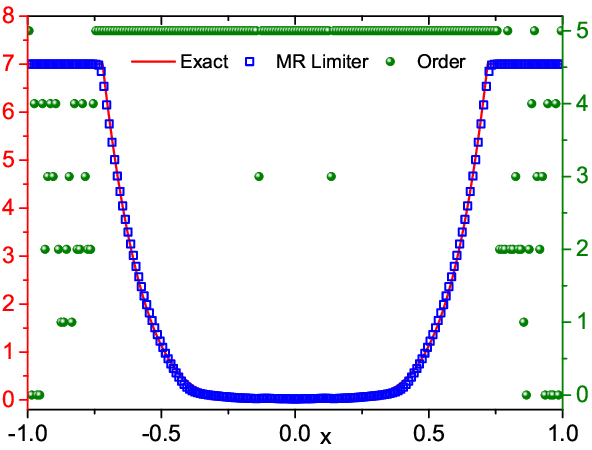}}
  \subfigure[$k=6$]{
  \includegraphics[width=5.5 cm]{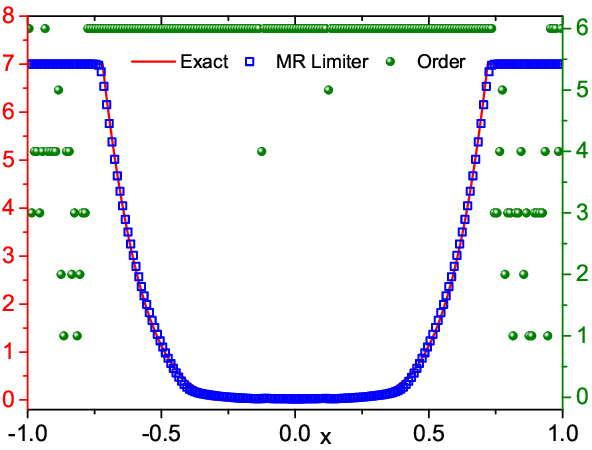}}
  \subfigure[$k=1$]{
  \includegraphics[width=5.5 cm]{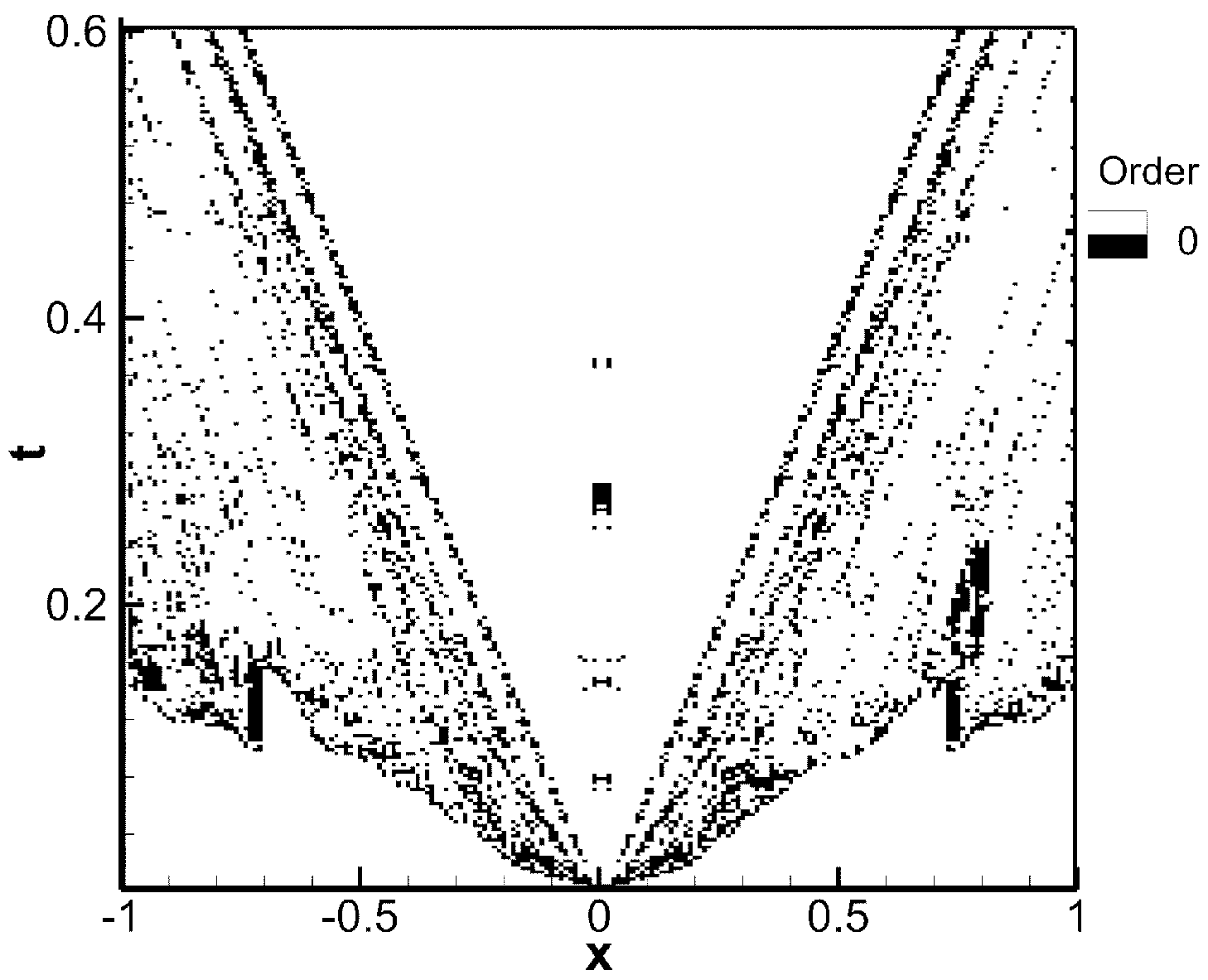}}
  \subfigure[$k=2$]{
  \includegraphics[width=5.5 cm]{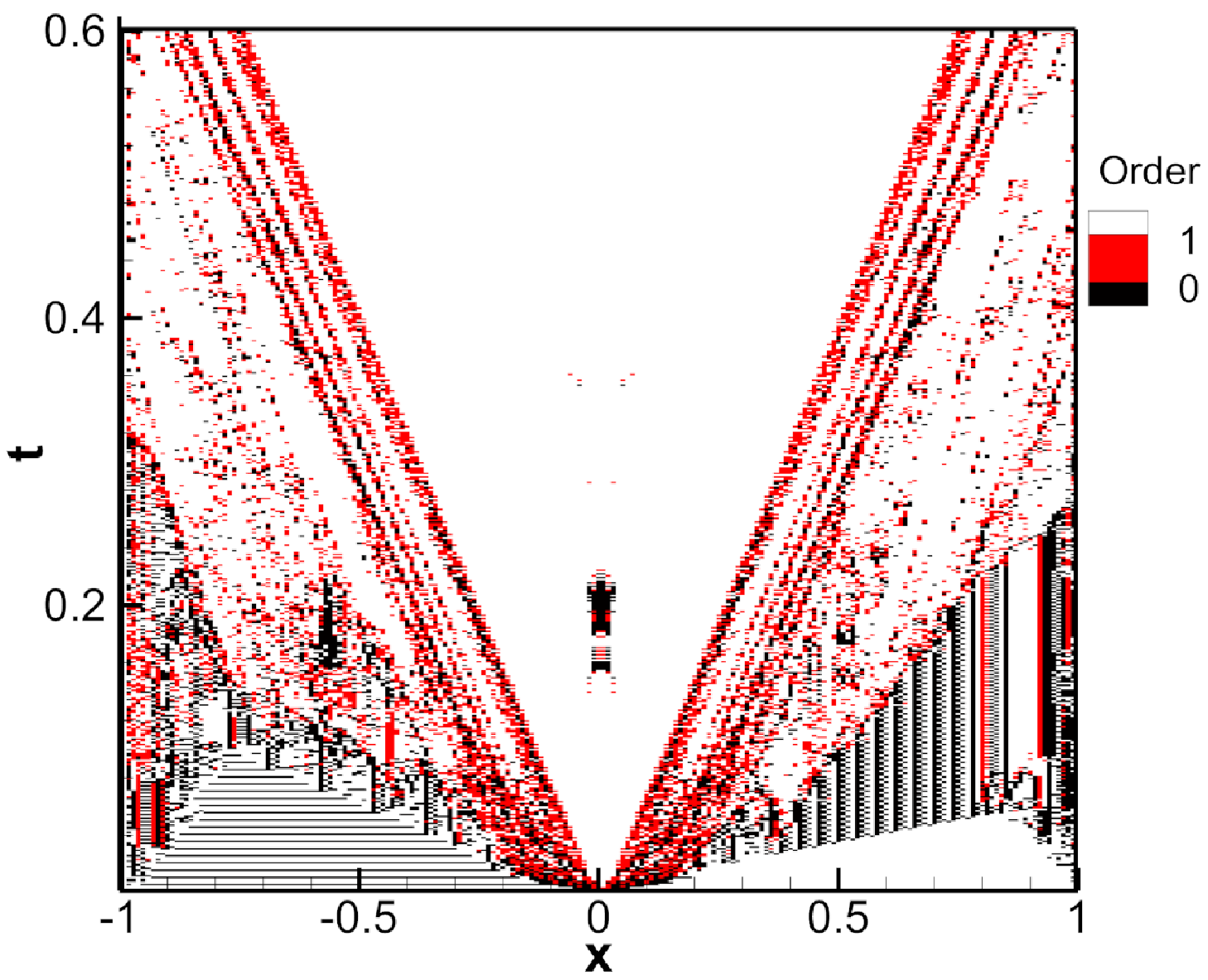}}
  \subfigure[$k=3$]{
  \includegraphics[width=5.5 cm]{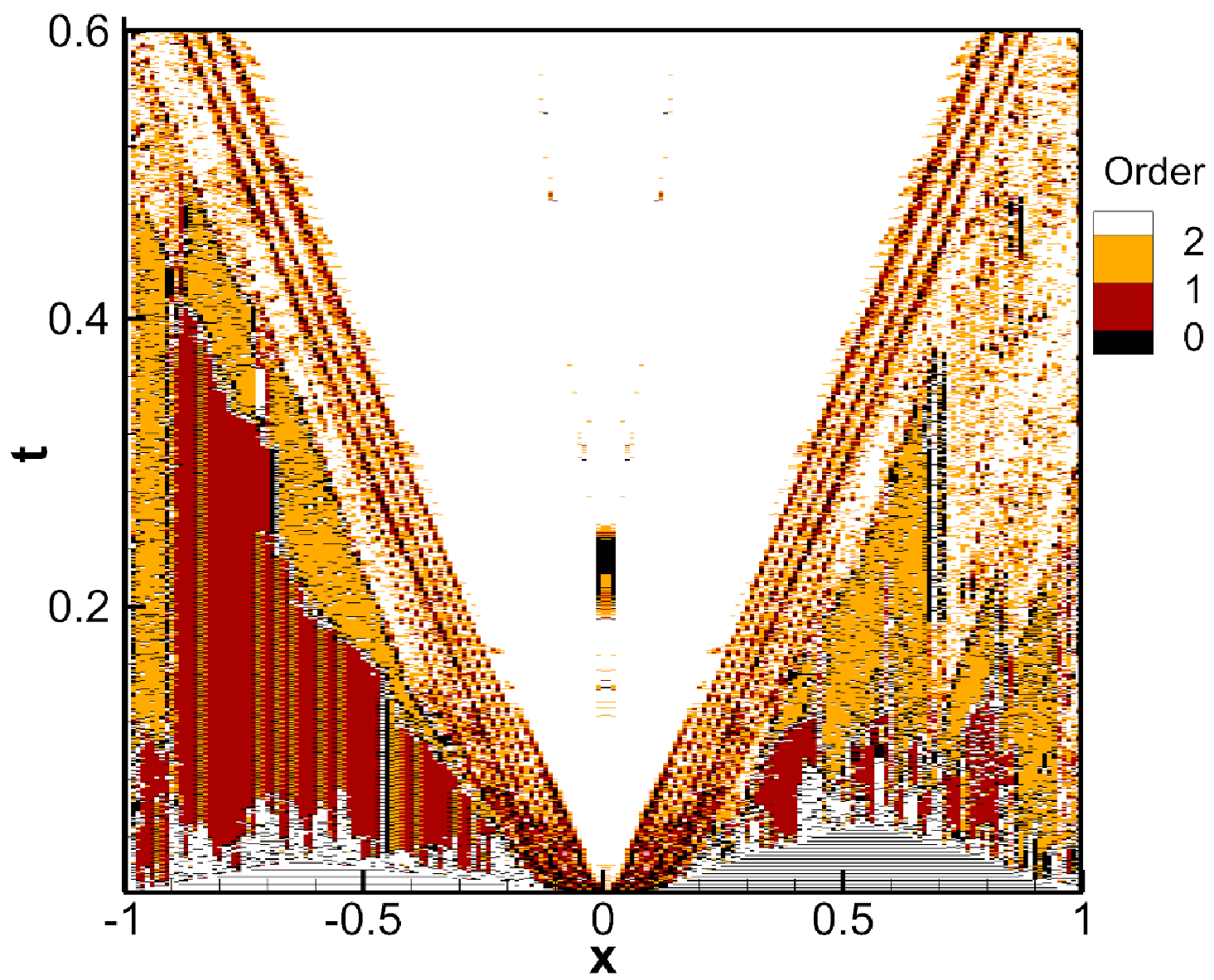}}
  \subfigure[$k=4$]{
  \includegraphics[width=5.5 cm]{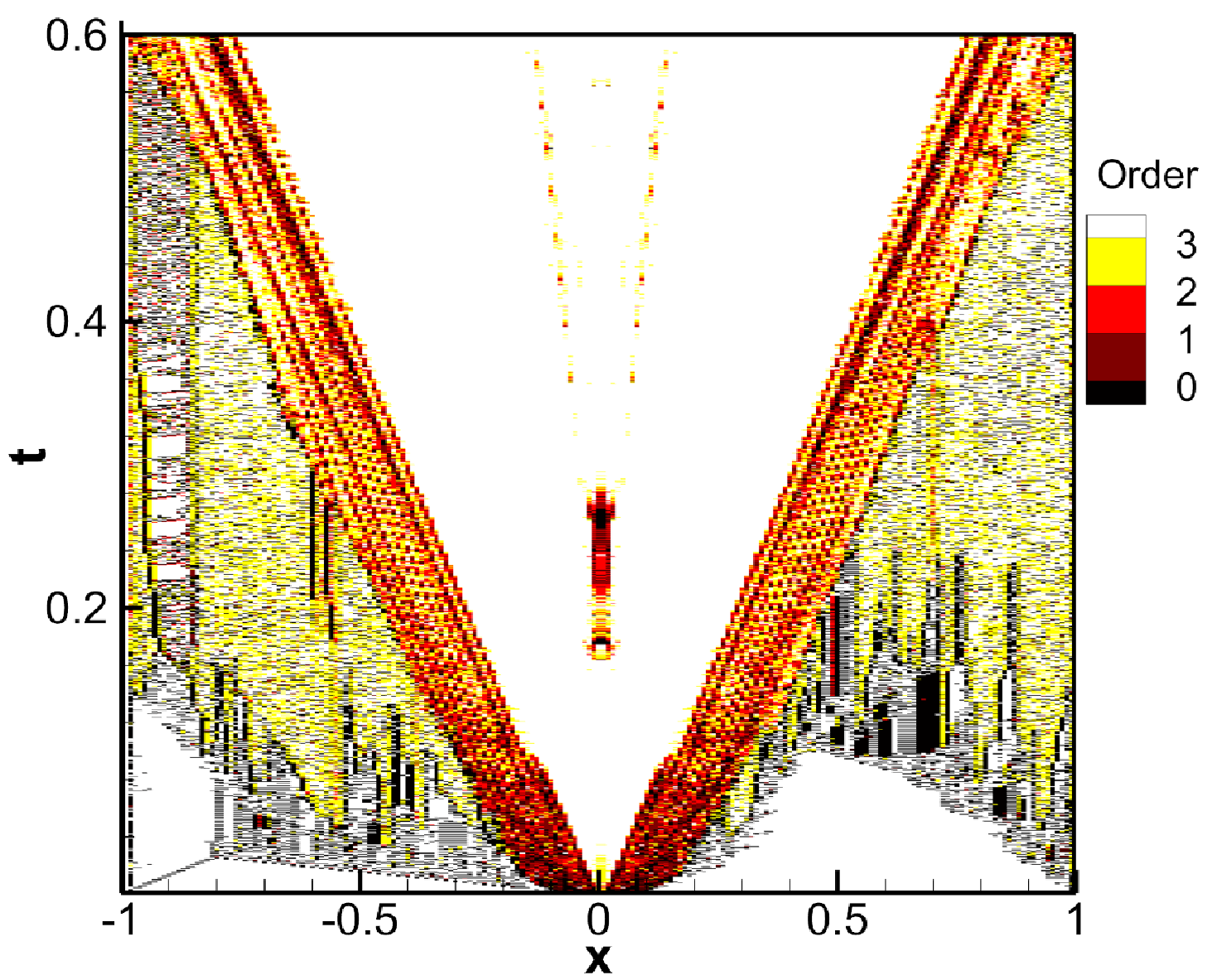}}
  \subfigure[$k=5$]{
  \includegraphics[width=5.5 cm]{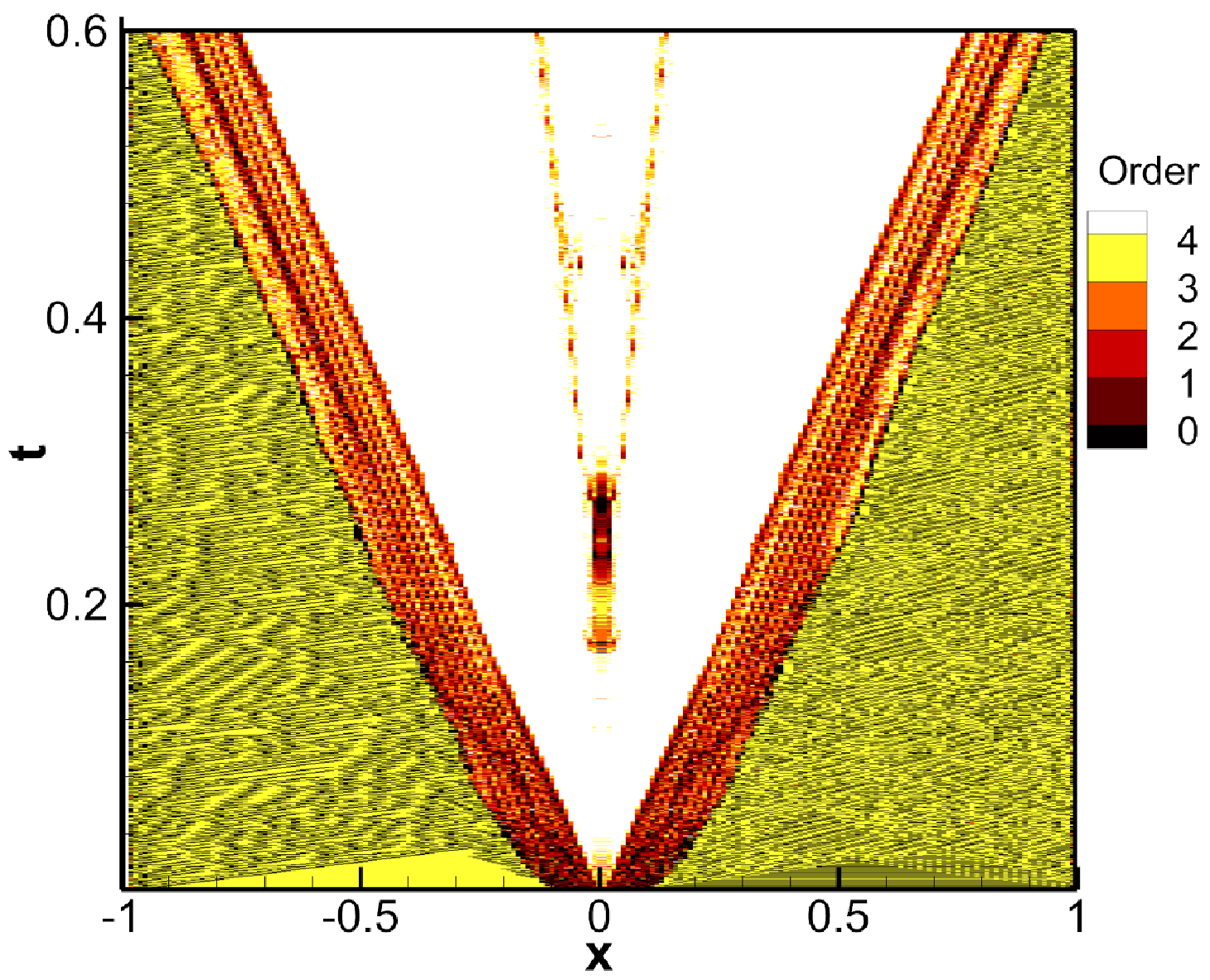}}
  \subfigure[$k=6$]{
  \includegraphics[width=5.5 cm]{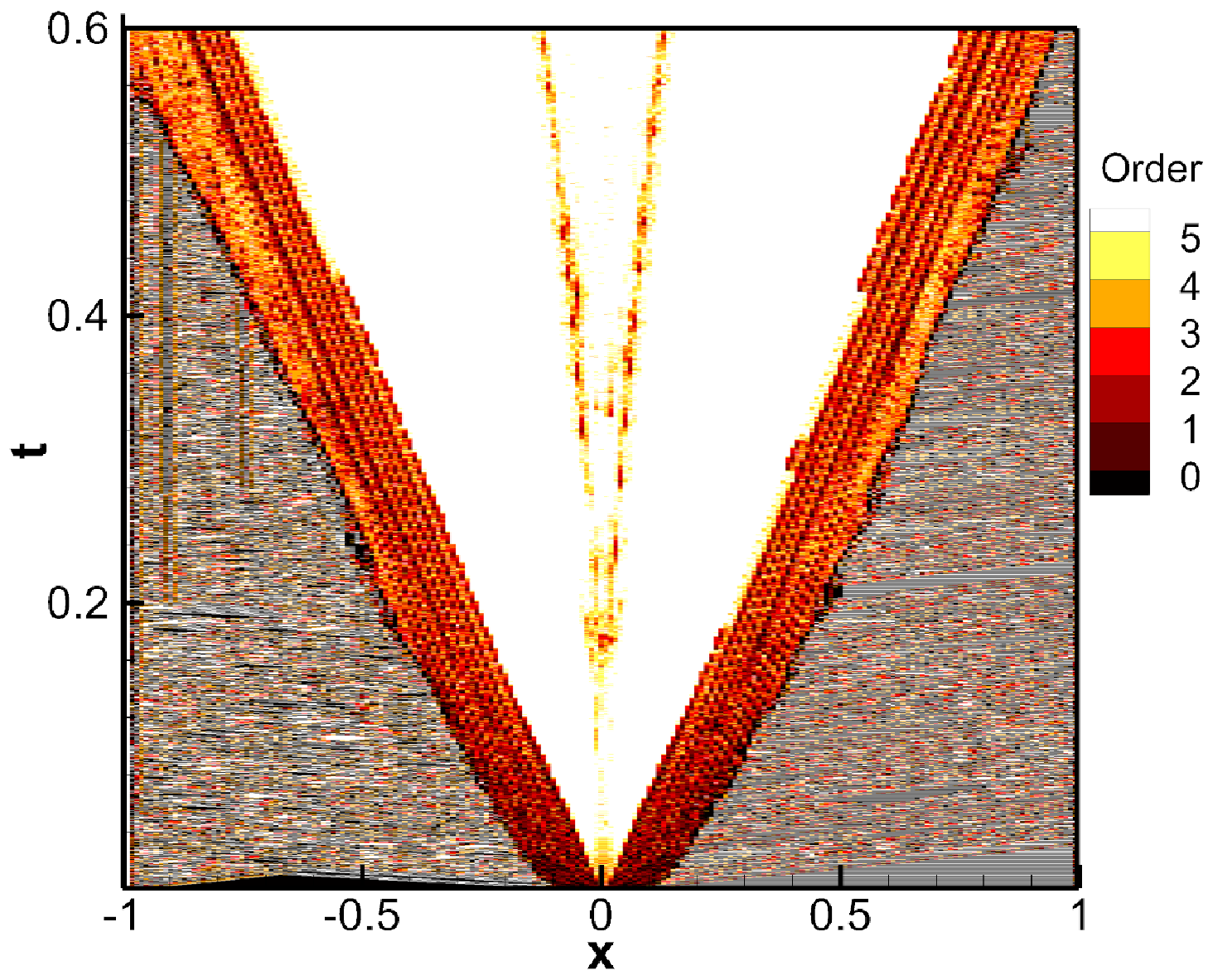}}
  \caption{The double rarefaction problem computed by RKDG schemes with the MR limiter using 200 cells. 
  The first and second rows: distributions of the density and the polynomial order at $t=0.6$;
   The third and fourth rows: the time history of the polynomial order, 
  the white parts represent the original $k$th-order DG polynomial. }
 \label{FIG:Rarefaction}
 \end{figure}

\begin{figure}[htbp]
  \centering
  \subfigure[$k=1$]{
  \includegraphics[width=5.5 cm]{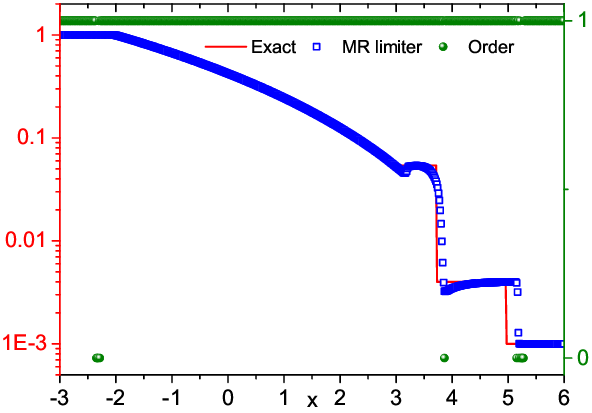}}
  \subfigure[$k=2$]{
  \includegraphics[width=5.5 cm]{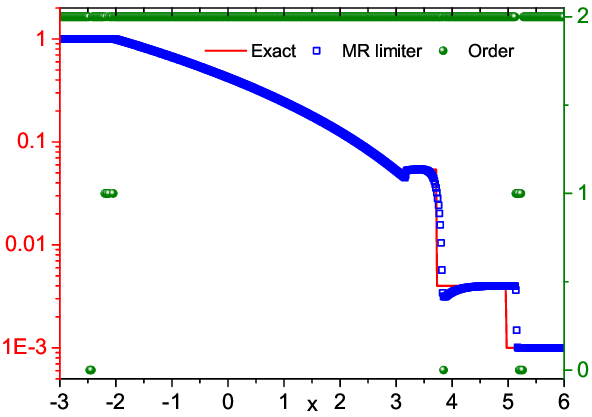}}
  \subfigure[$k=3$]{
  \includegraphics[width=5.5 cm]{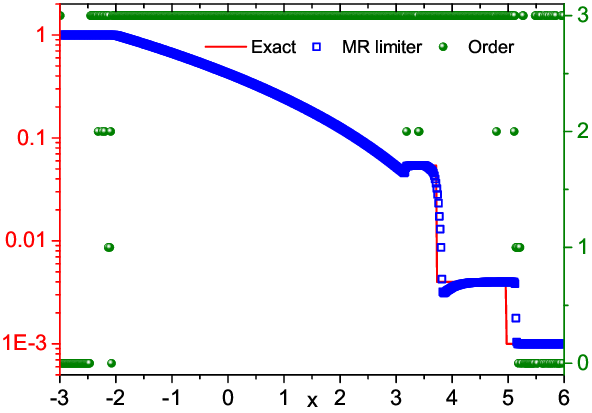}}
  \subfigure[$k=4$]{
  \includegraphics[width=5.5 cm]{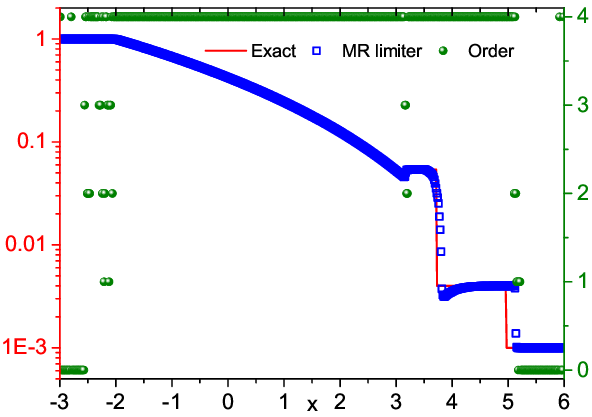}}
  \subfigure[$k=5$]{
  \includegraphics[width=5.5 cm]{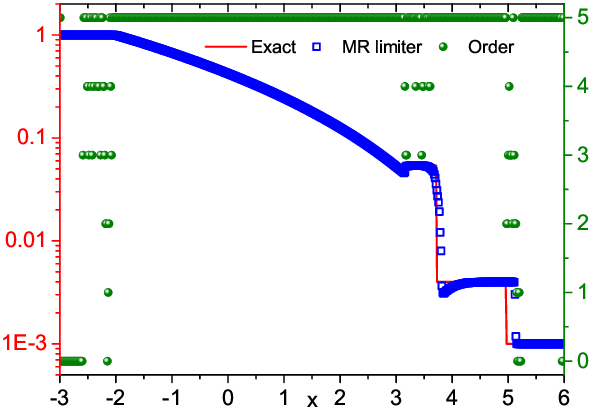}}
  \subfigure[$k=6$]{
  \includegraphics[width=5.5 cm]{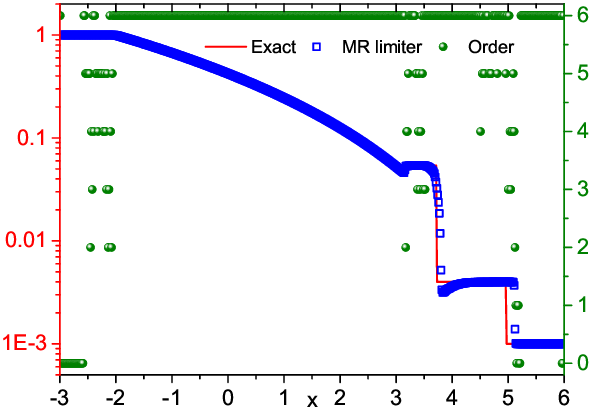}}
  \subfigure[$k=1$]{
  \includegraphics[width=5.5 cm]{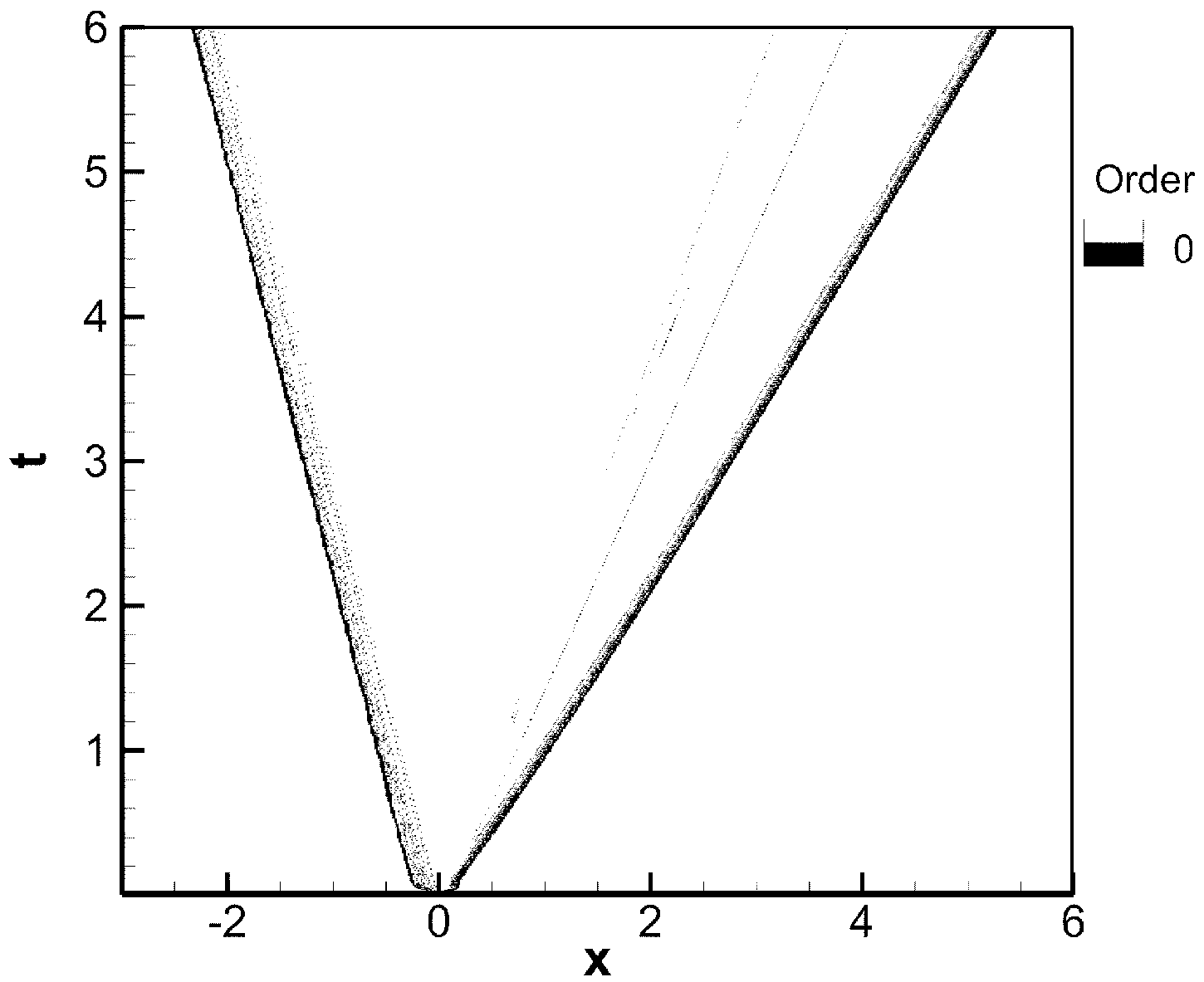}}
  \subfigure[$k=2$]{
  \includegraphics[width=5.5 cm]{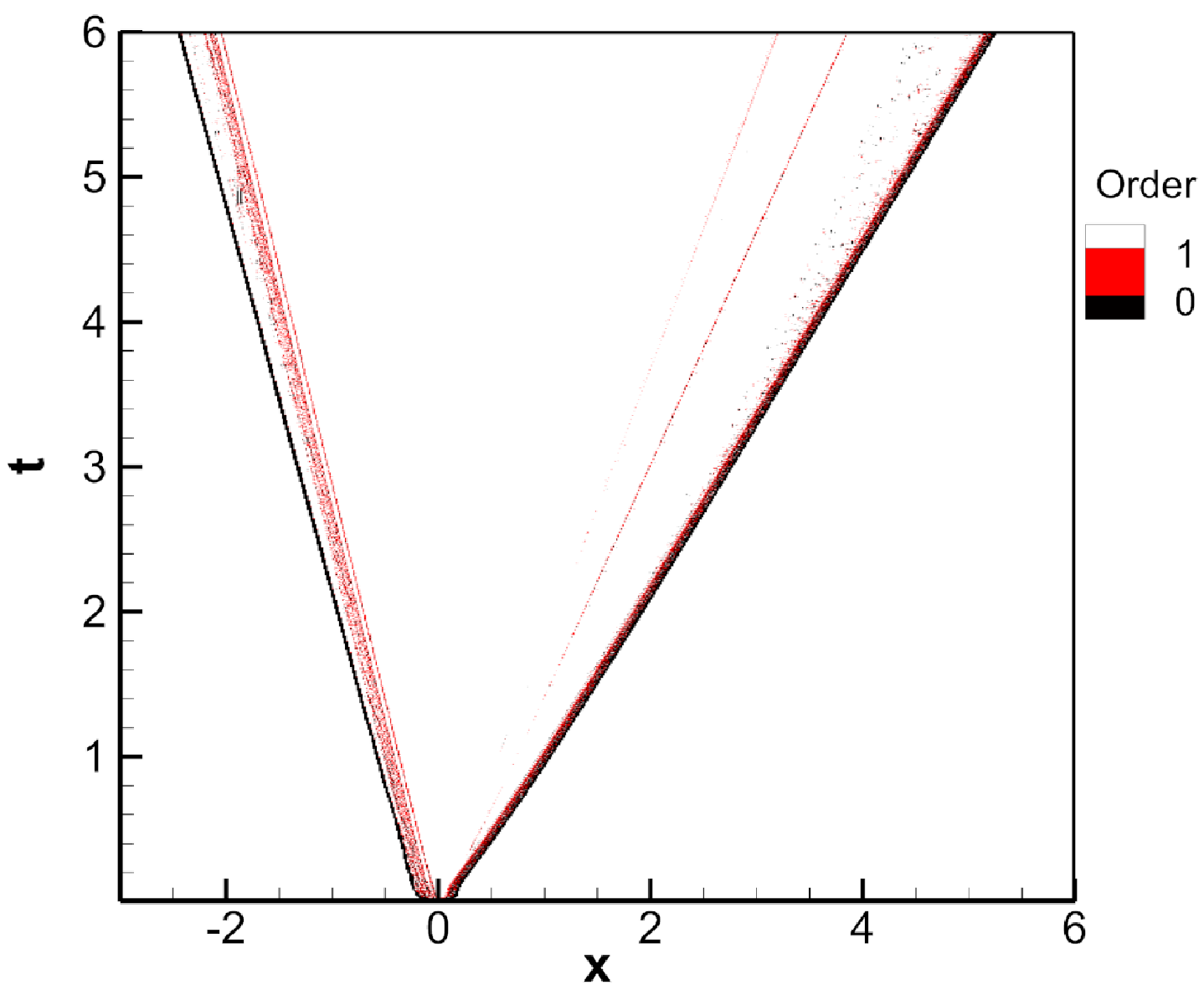}}
  \subfigure[$k=3$]{
  \includegraphics[width=5.5 cm]{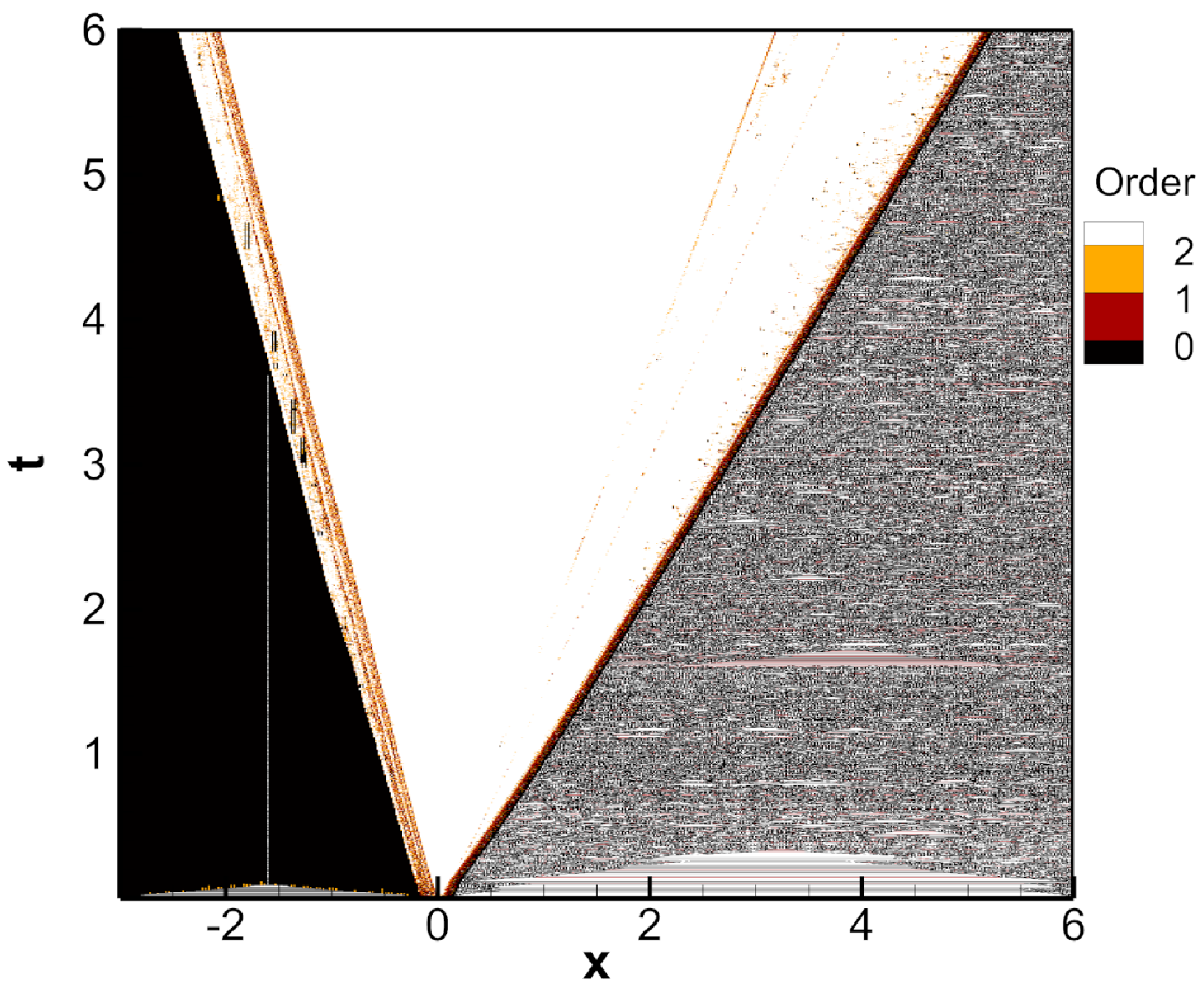}}
  \subfigure[$k=4$]{
  \includegraphics[width=5.5 cm]{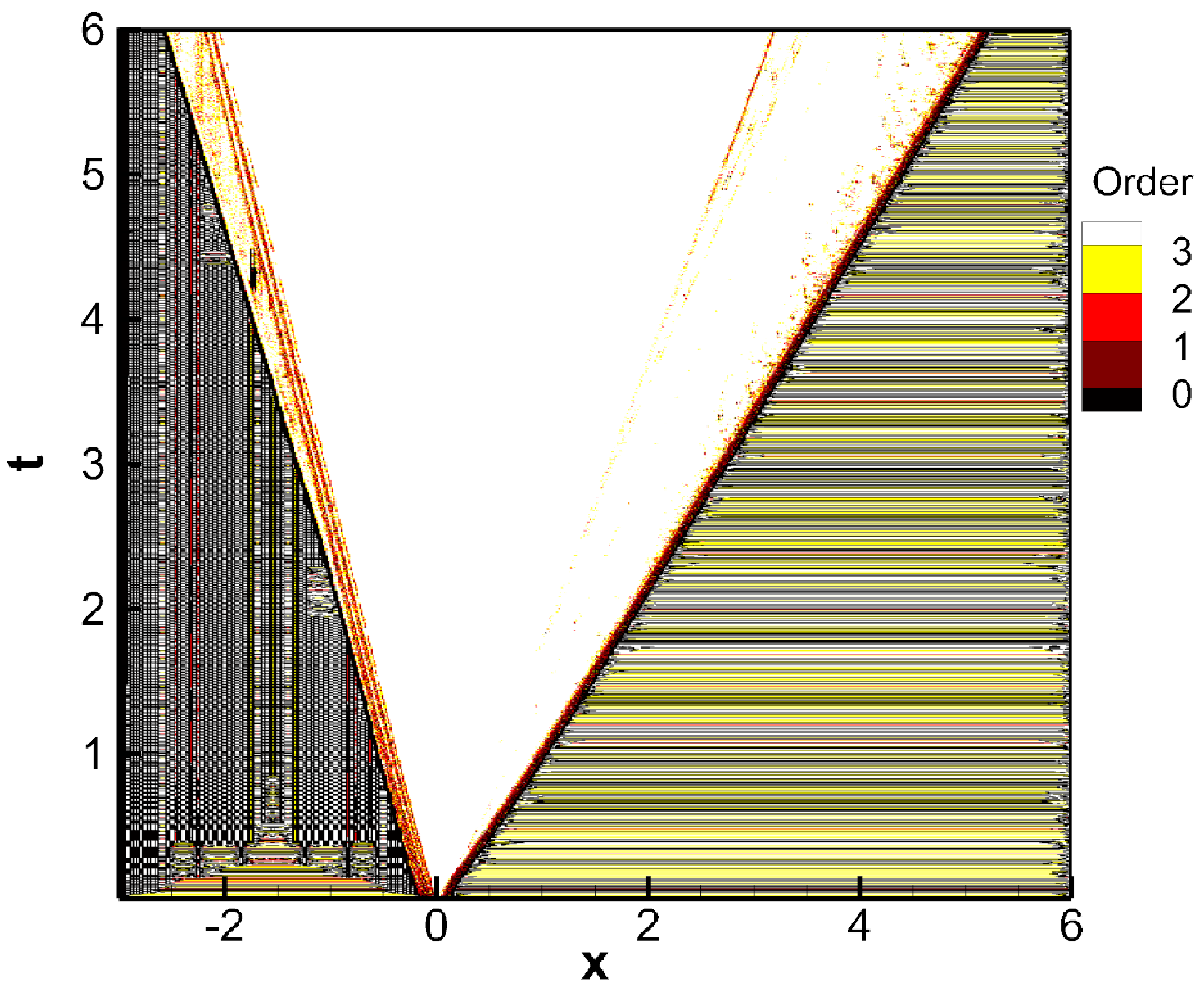}}
  \subfigure[$k=5$]{
  \includegraphics[width=5.5 cm]{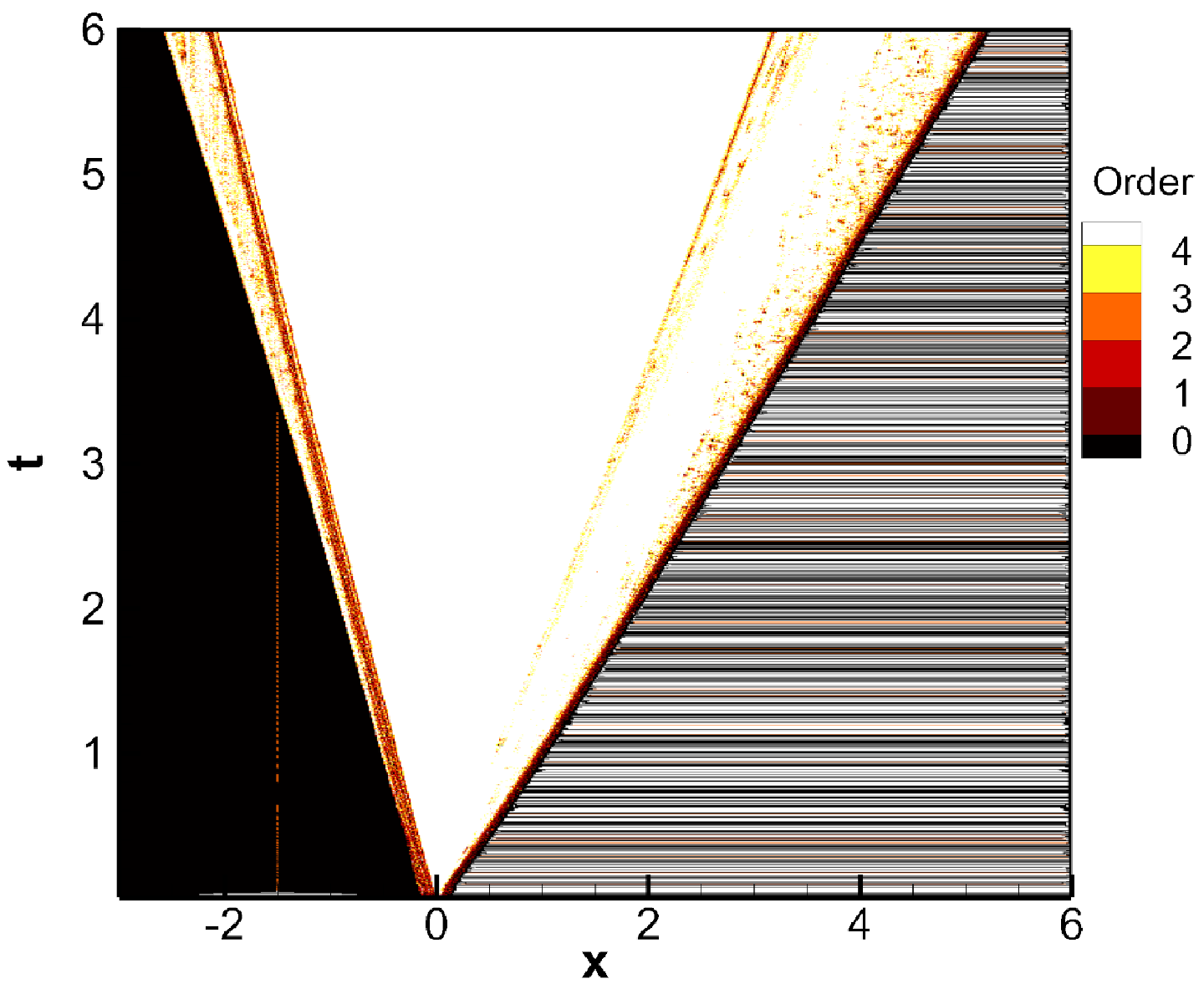}}
  \subfigure[$k=6$]{
  \includegraphics[width=5.5 cm]{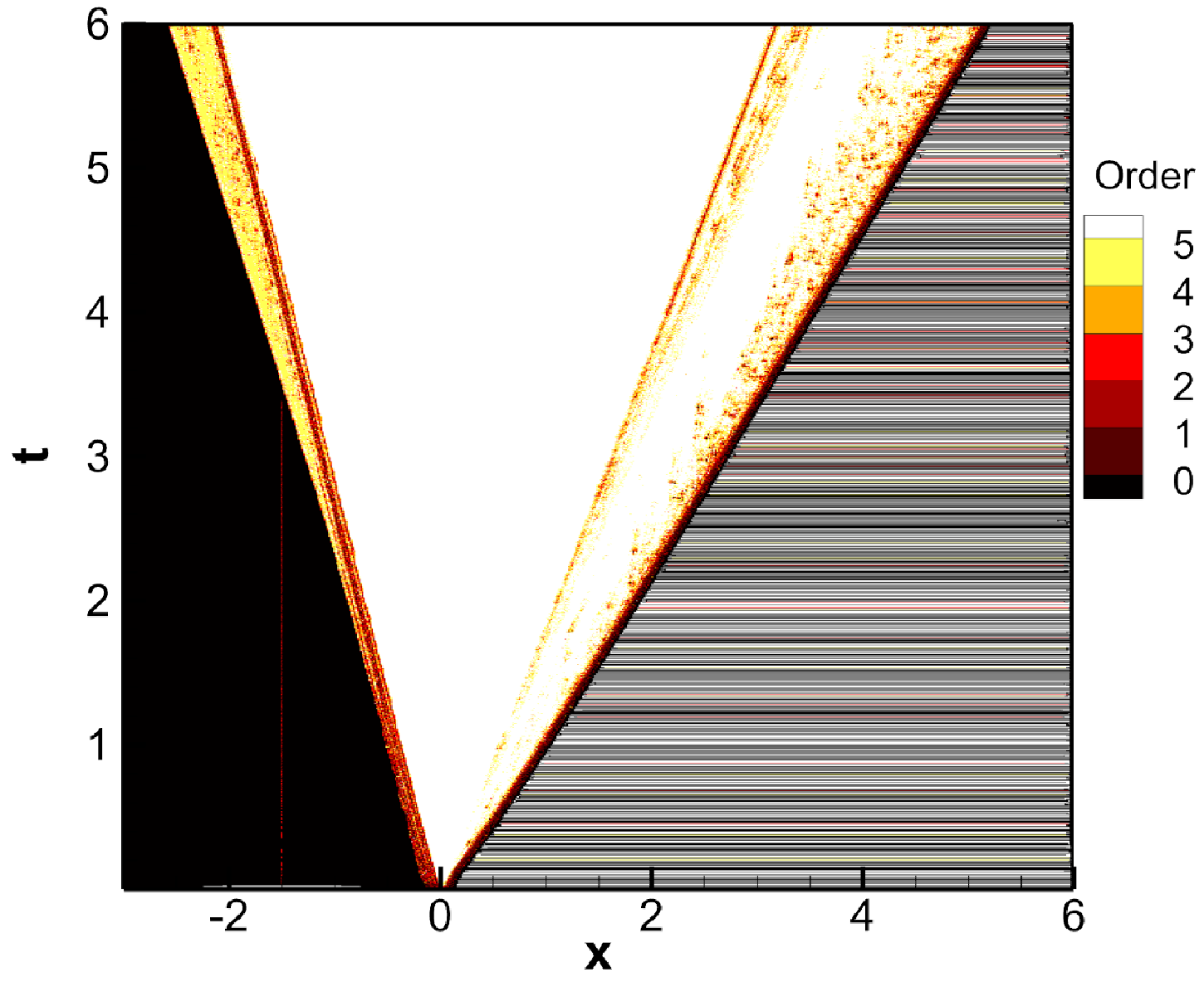}}
  \caption{The LeBlanc problem computed by RKDG schemes with the MR limiter using 600 cells. 
  The first and second rows: distributions of the density and the polynomial order at $t=6$;
   The third and fourth rows: the time history of the polynomial order, 
  the white parts represent the original $k$th-order DG polynomial. }
 \label{FIG:LeBlanc}
 \end{figure}

\paragraph{Example 4.2.4} The LeBlanc problem.
The initial condition is given by
\begin{equation*}
 (\rho,u,p)=\begin{cases}
              (1,0,\frac{0.2}{3}), & \mbox{if } x<0, \\
              (10^{-3},1,\frac{2}{3}\times 10^{-10}), & \mbox{otherwise},
            \end{cases}
\end{equation*}
The computational domain [-3,6], and $\gamma$ is set to $5/3$.
We compute this problem by the RKDG schemes with the MR limiter using 600 cells.
Fig. \ref{FIG:LeBlanc} shows the distributions of the density and the polynomial order at $t=6$
and the time history of the polynomial order.
The MR limiter preserves high-order for the expansion wave 
and flags the cells near the shock as troubled cells.
We note that the shock speed is not very accurate for this extreme problem,
which is similar to the results of the FS indicator \cite{Fu2017NewLimiter}.

\begin{figure}[htbp]
  \centering
  \subfigure[$k=1$]{
  \includegraphics[width=5.5 cm]{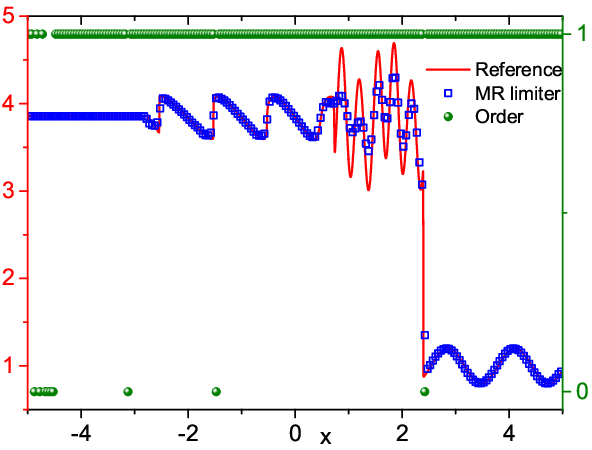}}
  \subfigure[$k=2$]{
  \includegraphics[width=5.5 cm]{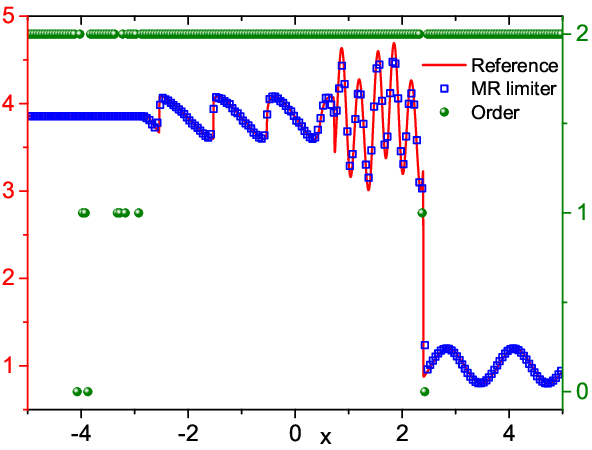}}
  \subfigure[$k=3$]{
  \includegraphics[width=5.5 cm]{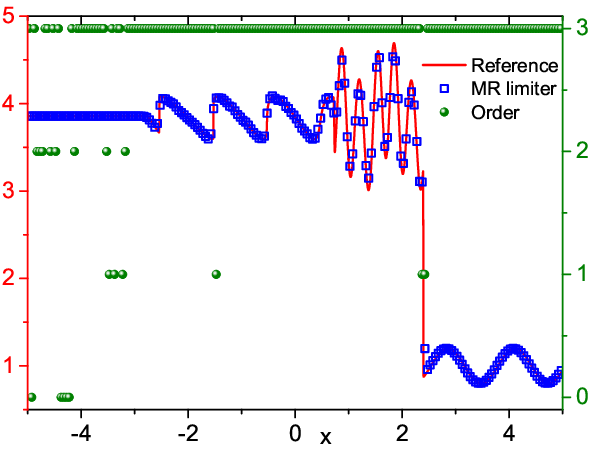}}
  \subfigure[$k=4$]{
  \includegraphics[width=5.5 cm]{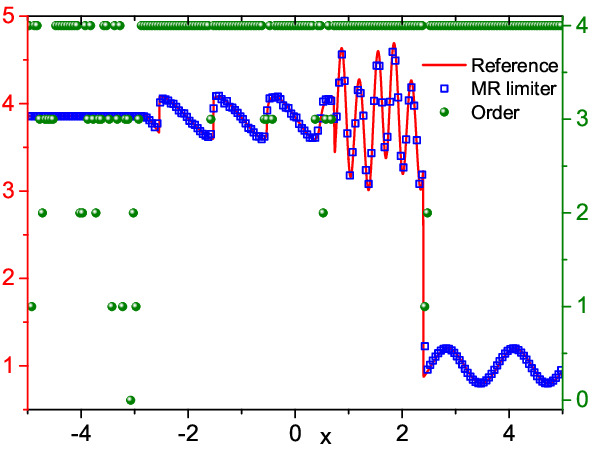}}
  \subfigure[$k=5$]{
  \includegraphics[width=5.5 cm]{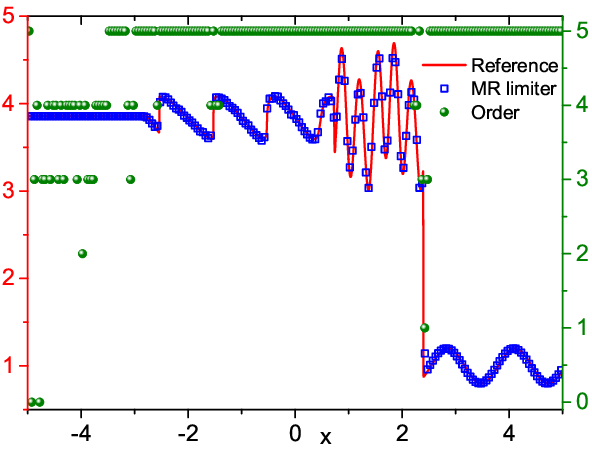}}
  \subfigure[$k=6$]{
  \includegraphics[width=5.5 cm]{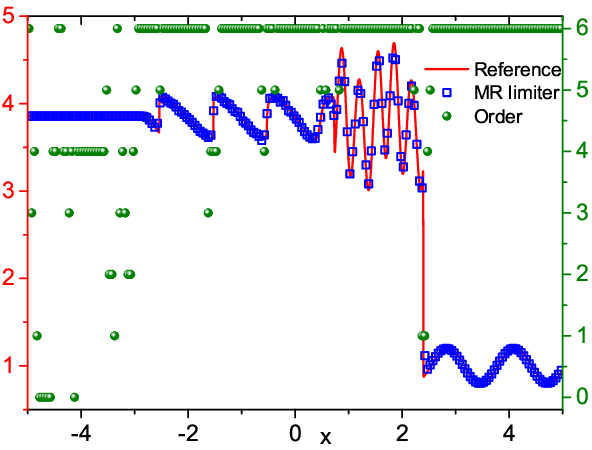}}
  \subfigure[$k=1$]{
  \includegraphics[width=5.5 cm]{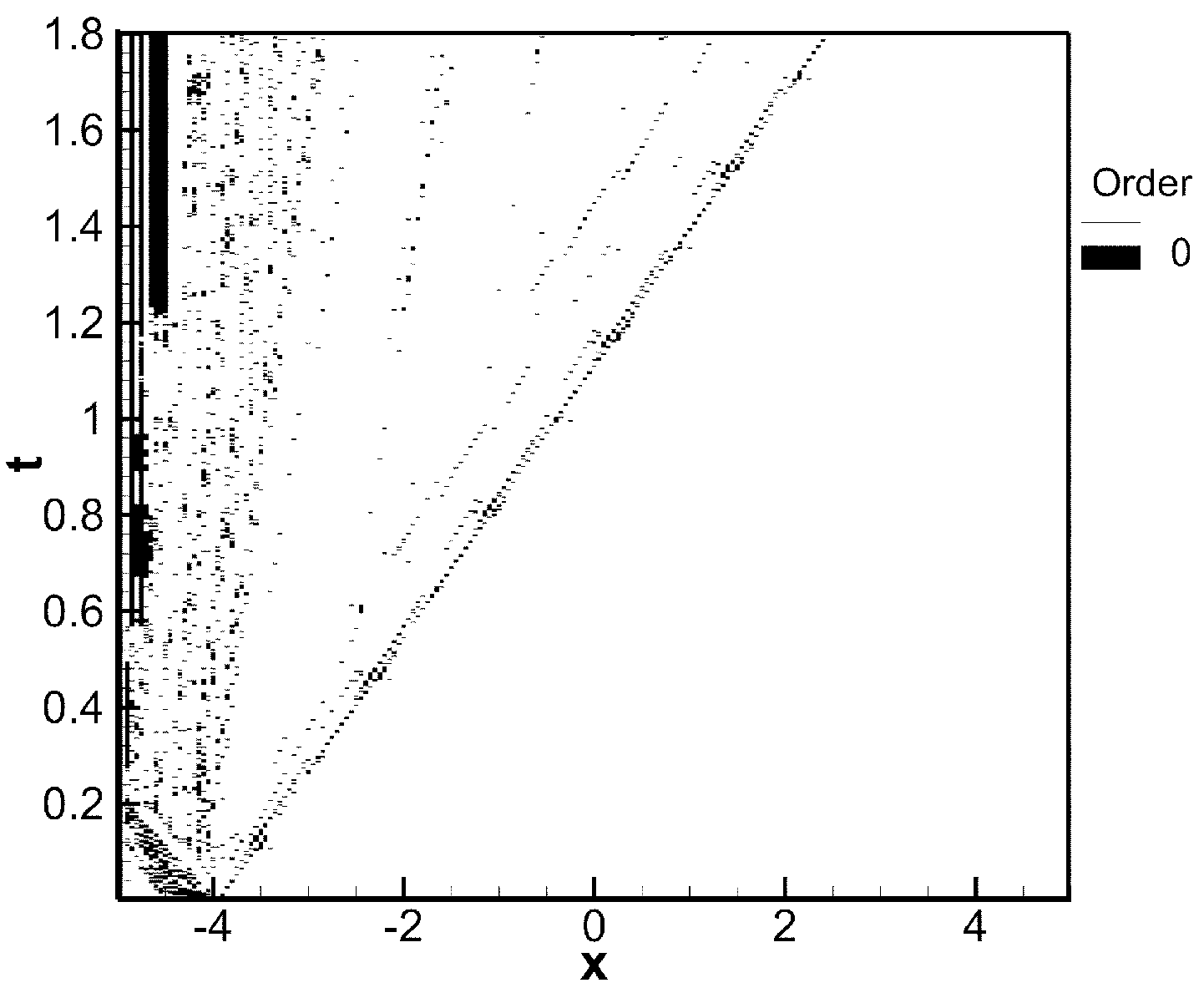}}
  \subfigure[$k=2$]{
  \includegraphics[width=5.5 cm]{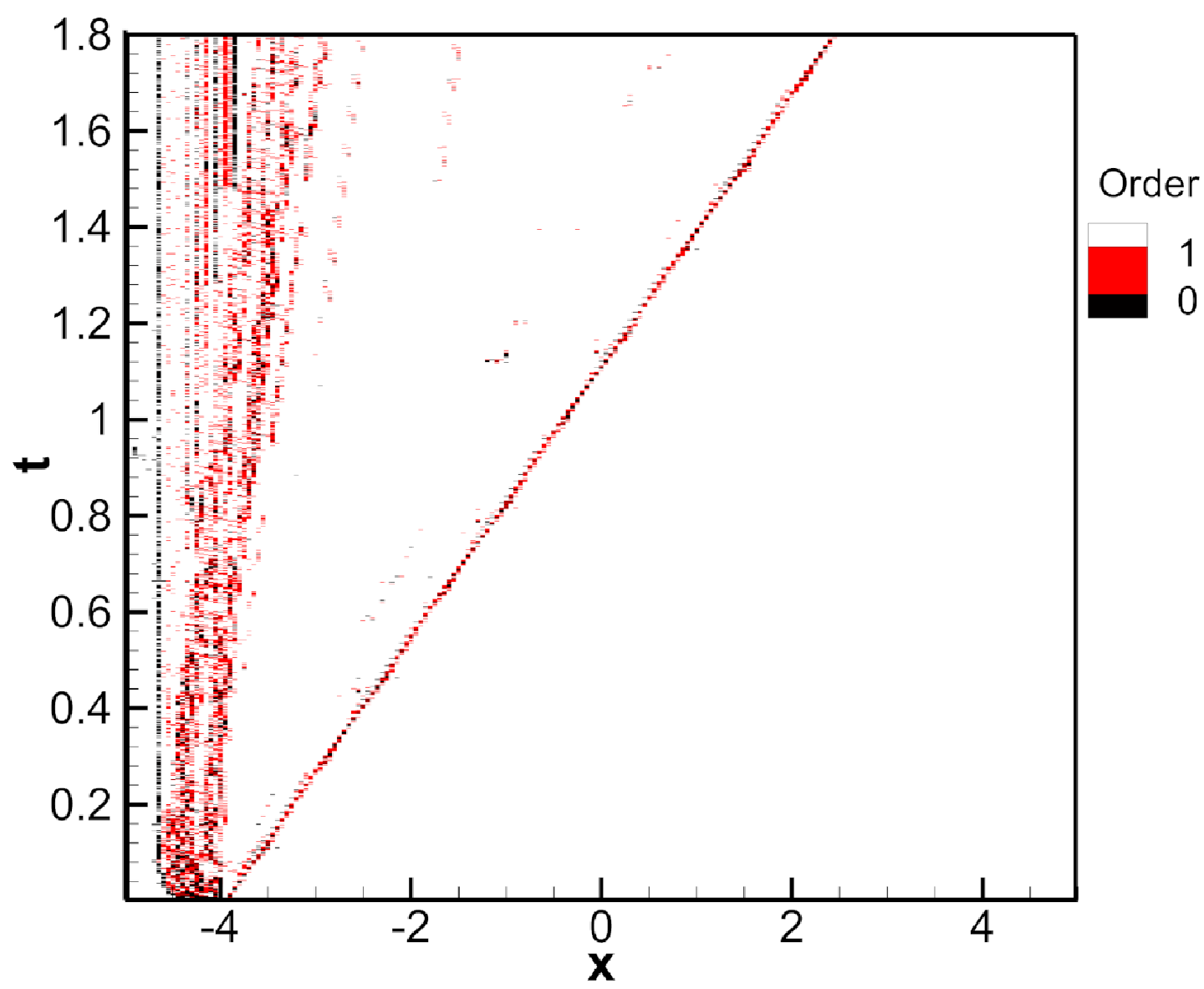}}
  \subfigure[$k=3$]{
  \includegraphics[width=5.5 cm]{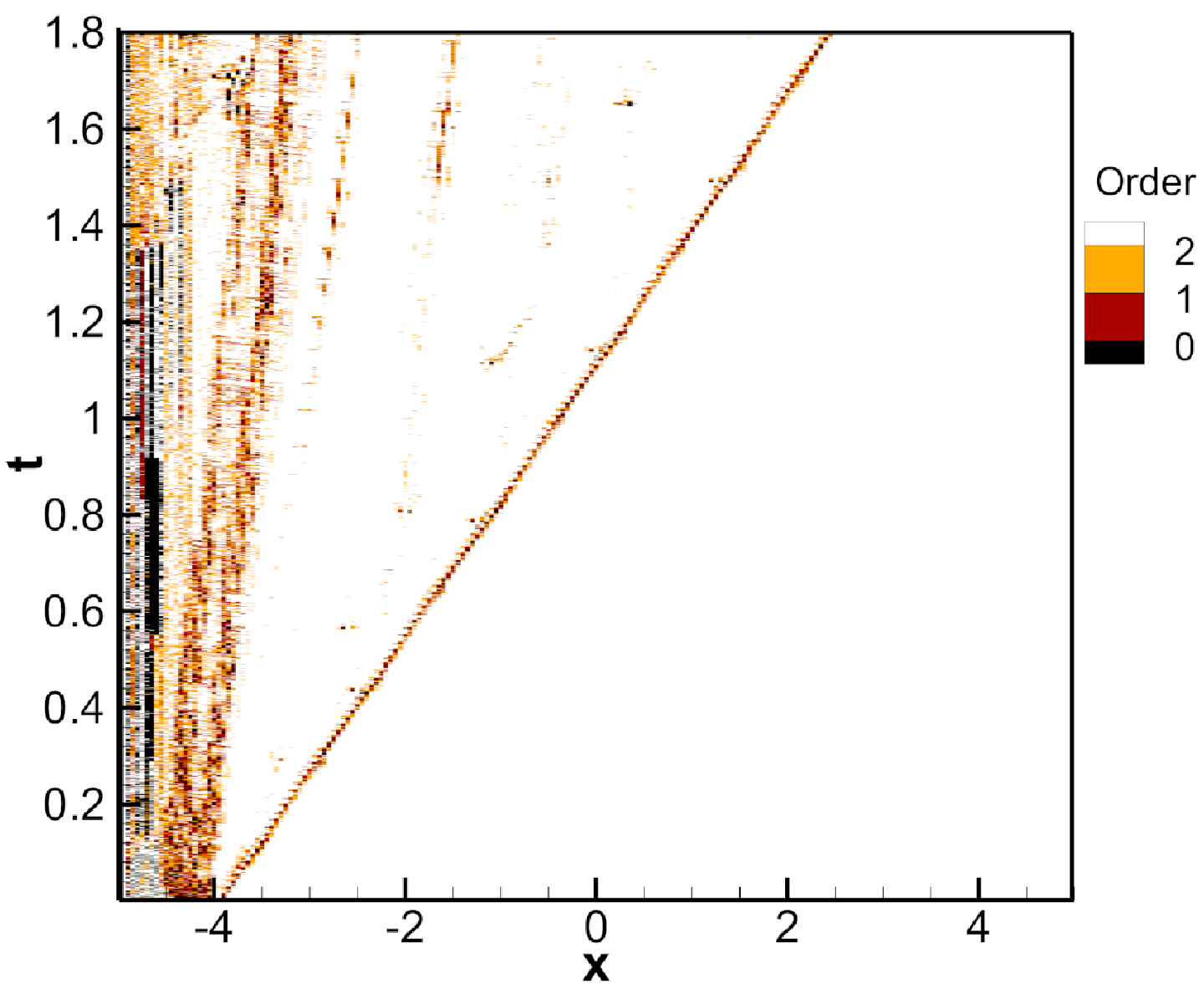}}
  \subfigure[$k=4$]{
  \includegraphics[width=5.5 cm]{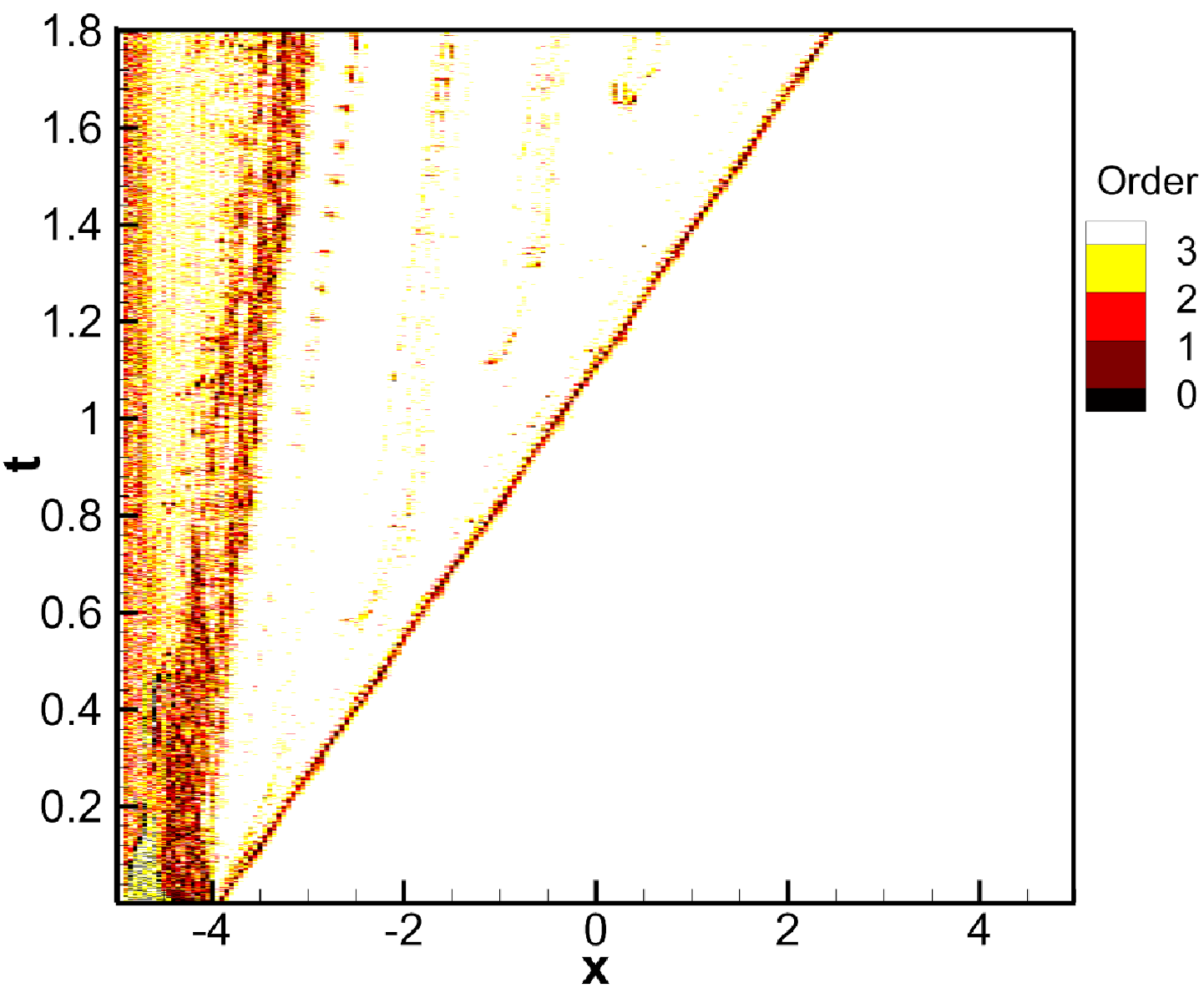}}
  \subfigure[$k=5$]{
  \includegraphics[width=5.5 cm]{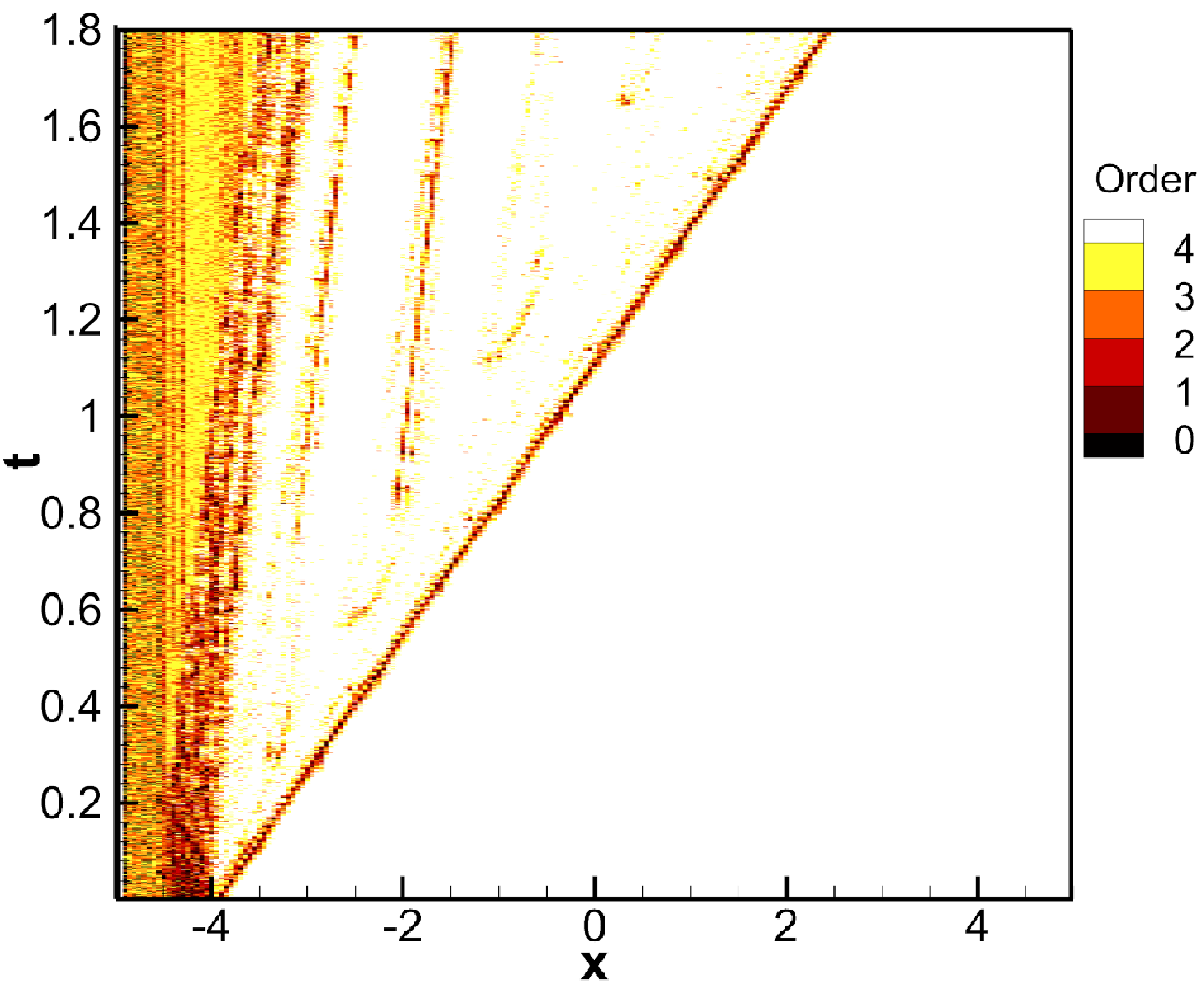}}
  \subfigure[$k=6$]{
  \includegraphics[width=5.5 cm]{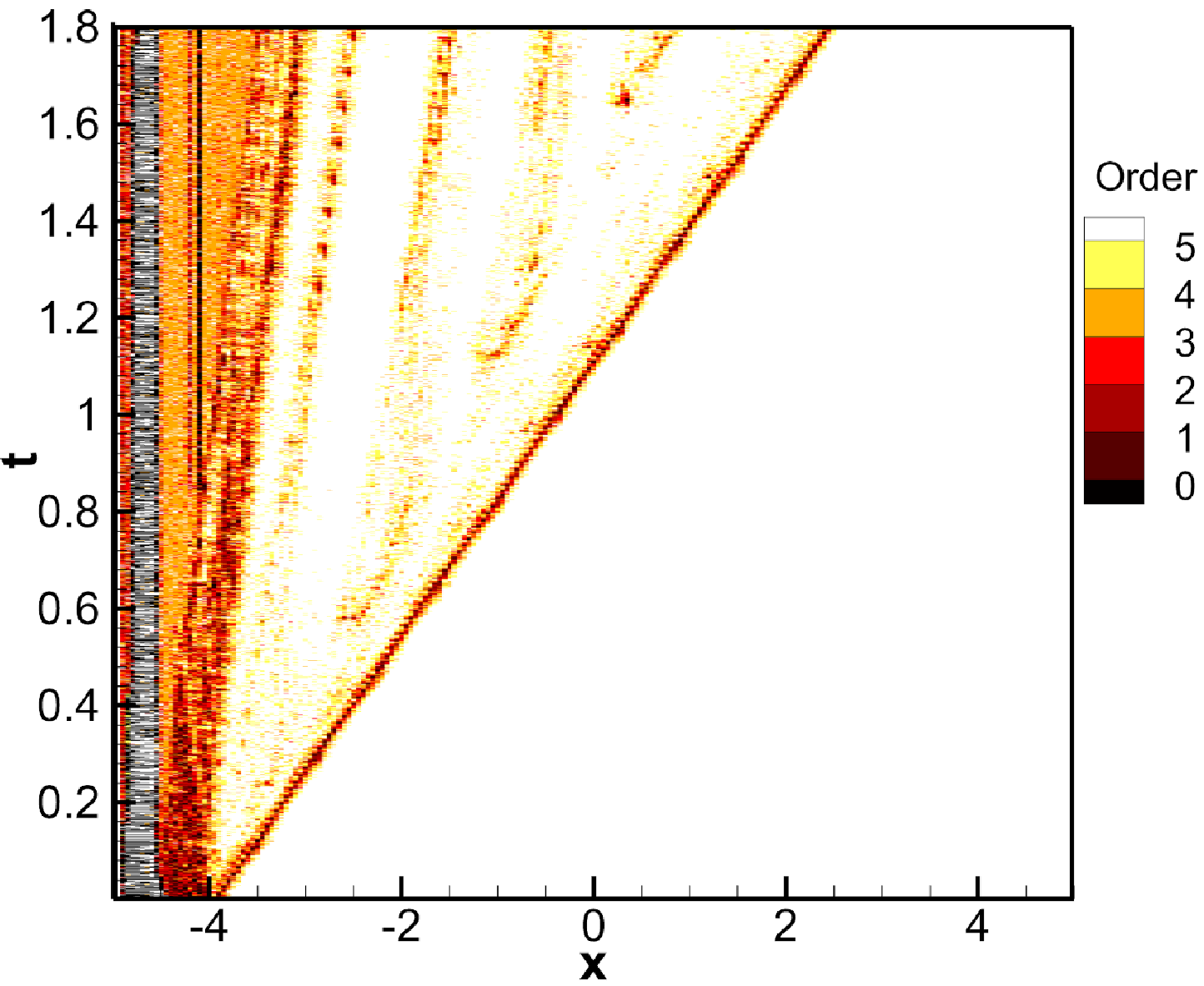}}
  \caption{The Shu-Osher problem computed by RKDG schemes with the MR limiter using 200 cells. 
  The first and second rows: distributions of the density and the polynomial order at $t=1.8$;
   The third and fourth rows: the time history of the polynomial order, 
  the white parts represent the original $k$th-order DG polynomial. }
 \label{FIG:Shu_Osher}
 \end{figure}

 \paragraph{Example 4.2.5} The Shu-Osher problem.
The computational domain is [-5,5], and the initial condition is given by
\begin{equation*}
 (\rho,u,p)=\begin{cases}
              (\frac{27}{7},\frac{4\sqrt{35}}{9},\frac{31}{3}), & \mbox{if } x<-4 \\
              (1+0.2\mbox{sin}(5\pi x),0,1), & \mbox{otherwise}.
            \end{cases}
\end{equation*}
It depicts the interaction of a shock wave with a high-frequency entropy sine wave 
which induces very complicated flow structures.
Fig. \ref{FIG:Shu_Osher} shows the distributions of the density and the polynomial order at $t=1.8$
and the time history of the polynomial order
computed by RKDG schemes with the MR limiter using 200 cells.
The reference solution is calculated by the fifth-order WENO scheme with $h=1/1000$.
We can see that there are a few troubled cells in the high-frequency region for $k\le2$.
When $k\ge3$, although the limiter detects more troubled cells, 
the order of the polynomial still retains a high level in the high frequency region,
which is beneficial to capture the complex wave structures.
The density profiles show that the MR limiter well resolves the complex flow structures.
When $k\ge2$, the results are significantly better than those computed by the KXRCF indicator and the FS indicator
under the same grid resolution as shown by Fig. 3.10 in \cite{Fu2017NewLimiter}.
The results are also much better than the recently published results computed by 
the hybrid limiter, as shown by Fig. 7 in \cite{Wei2024HybridLimiter}.

\begin{figure}[htbp]
  \centering
  \subfigure[$k=1$]{
  \includegraphics[width=5.5 cm]{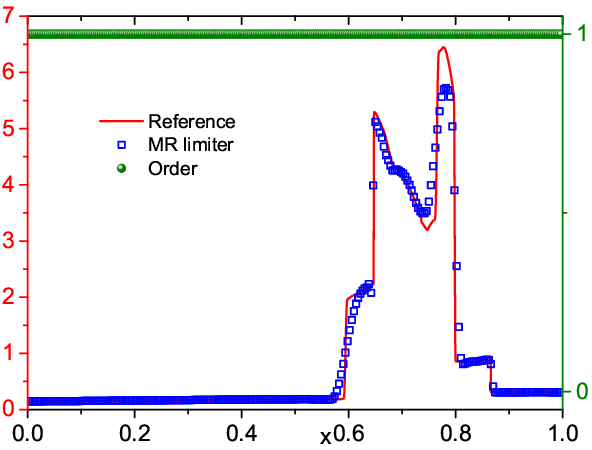}}
  \subfigure[$k=2$]{
  \includegraphics[width=5.5 cm]{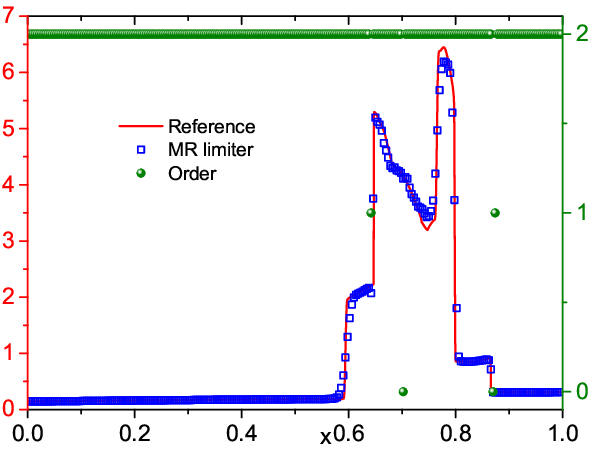}}
  \subfigure[$k=3$]{
  \includegraphics[width=5.5 cm]{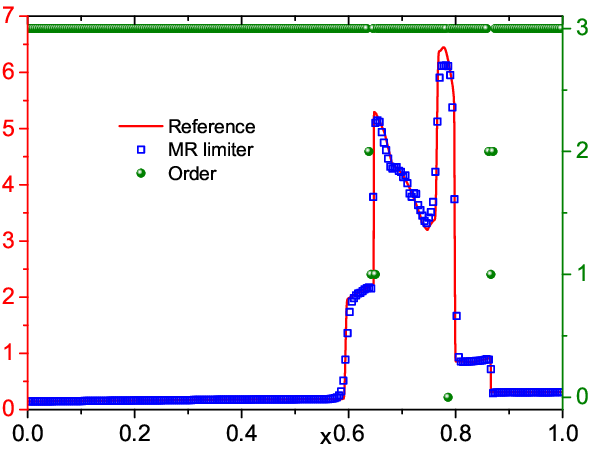}}
  \subfigure[$k=4$]{
  \includegraphics[width=5.5 cm]{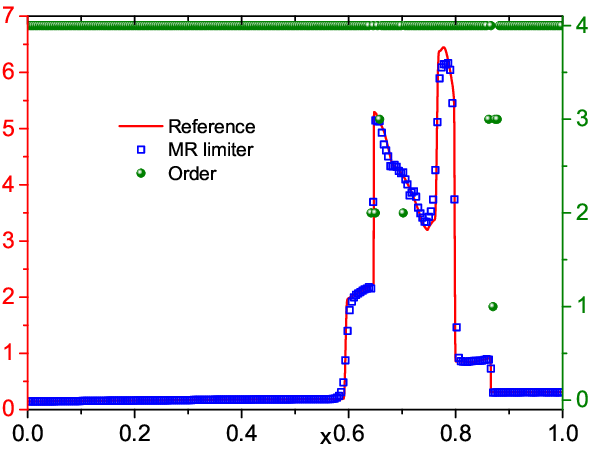}}
  \subfigure[$k=5$]{
  \includegraphics[width=5.5 cm]{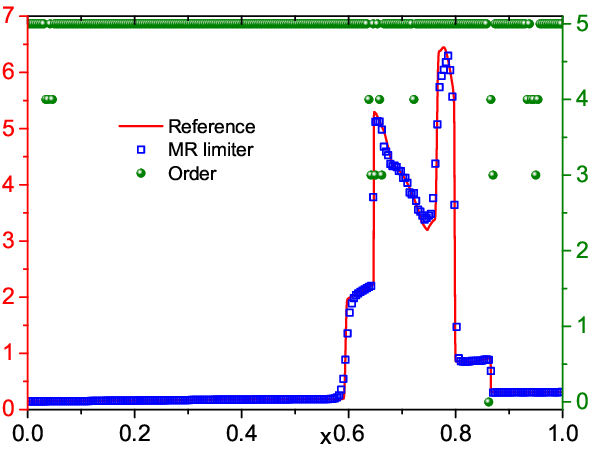}}
  \subfigure[$k=6$]{
  \includegraphics[width=5.5 cm]{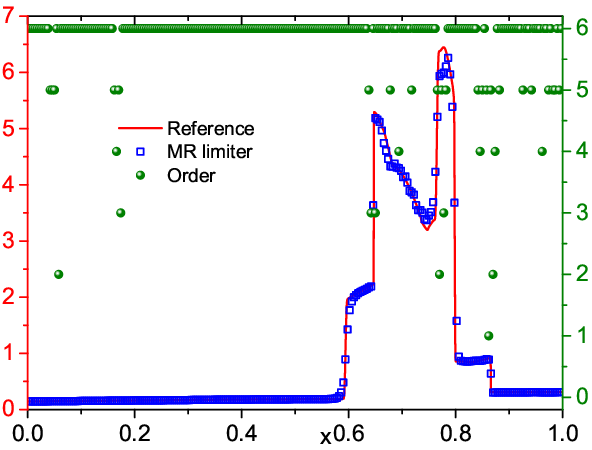}}
  \subfigure[$k=1$]{
  \includegraphics[width=5.5 cm]{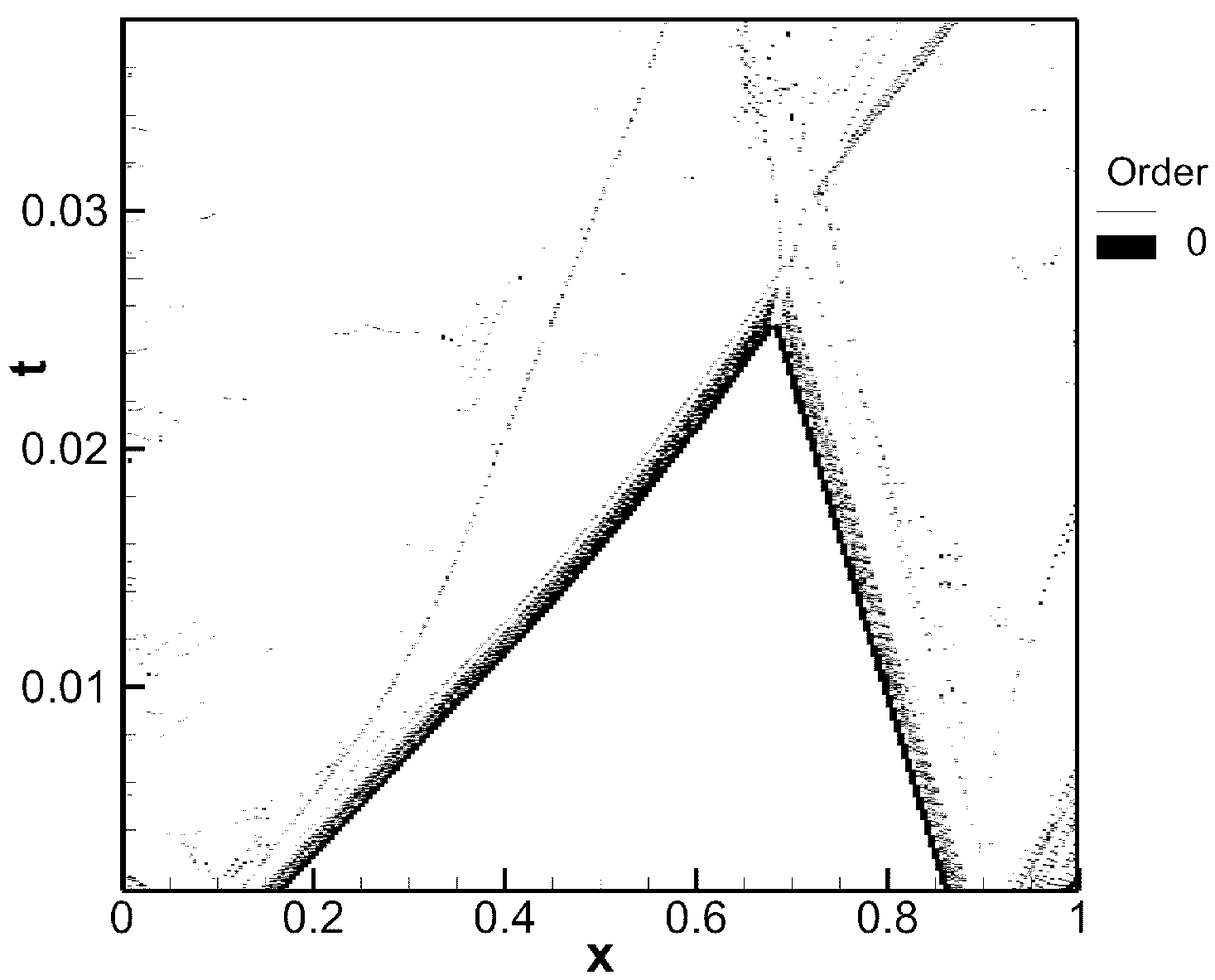}}
  \subfigure[$k=2$]{
  \includegraphics[width=5.5 cm]{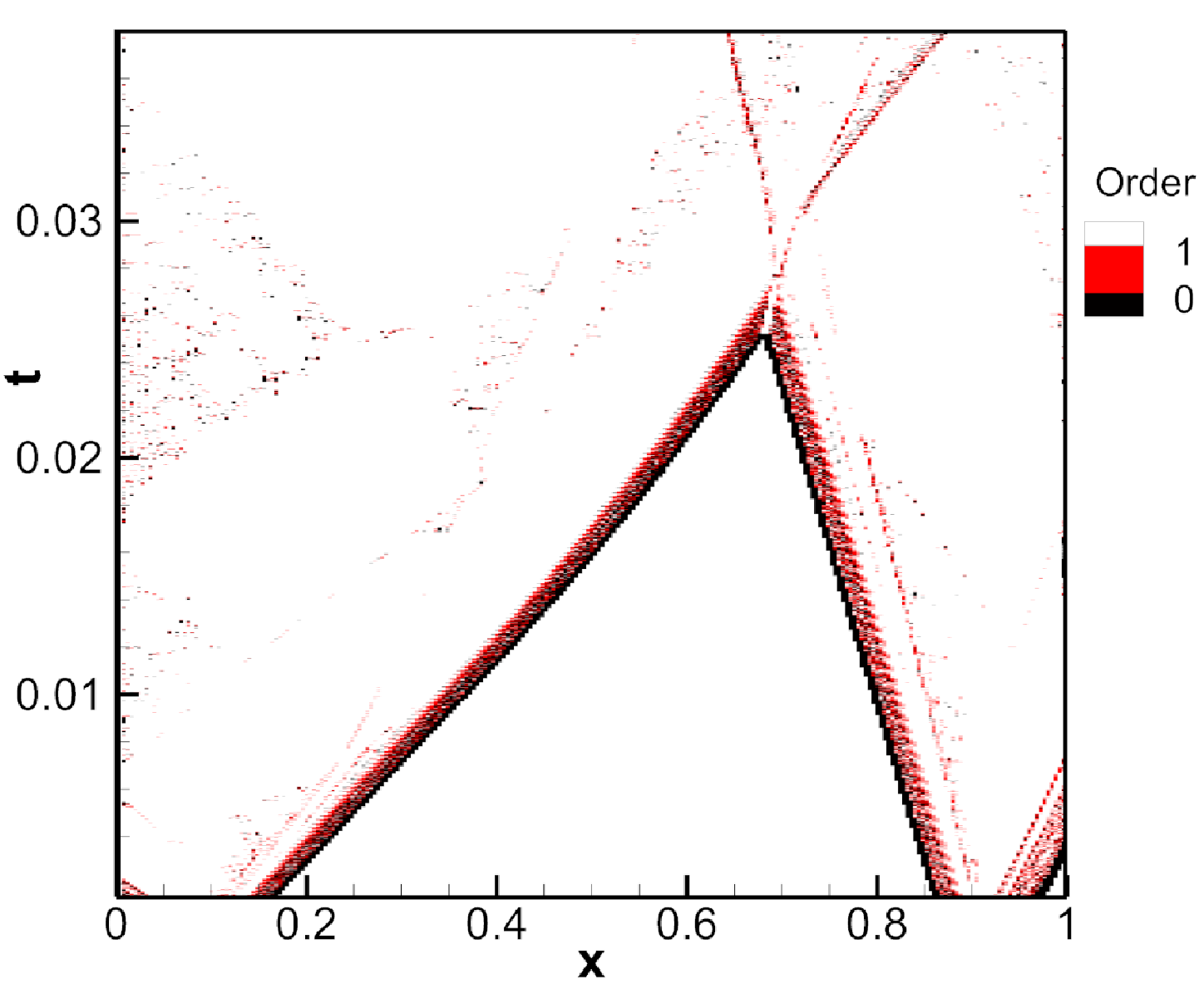}}
  \subfigure[$k=3$]{
  \includegraphics[width=5.5 cm]{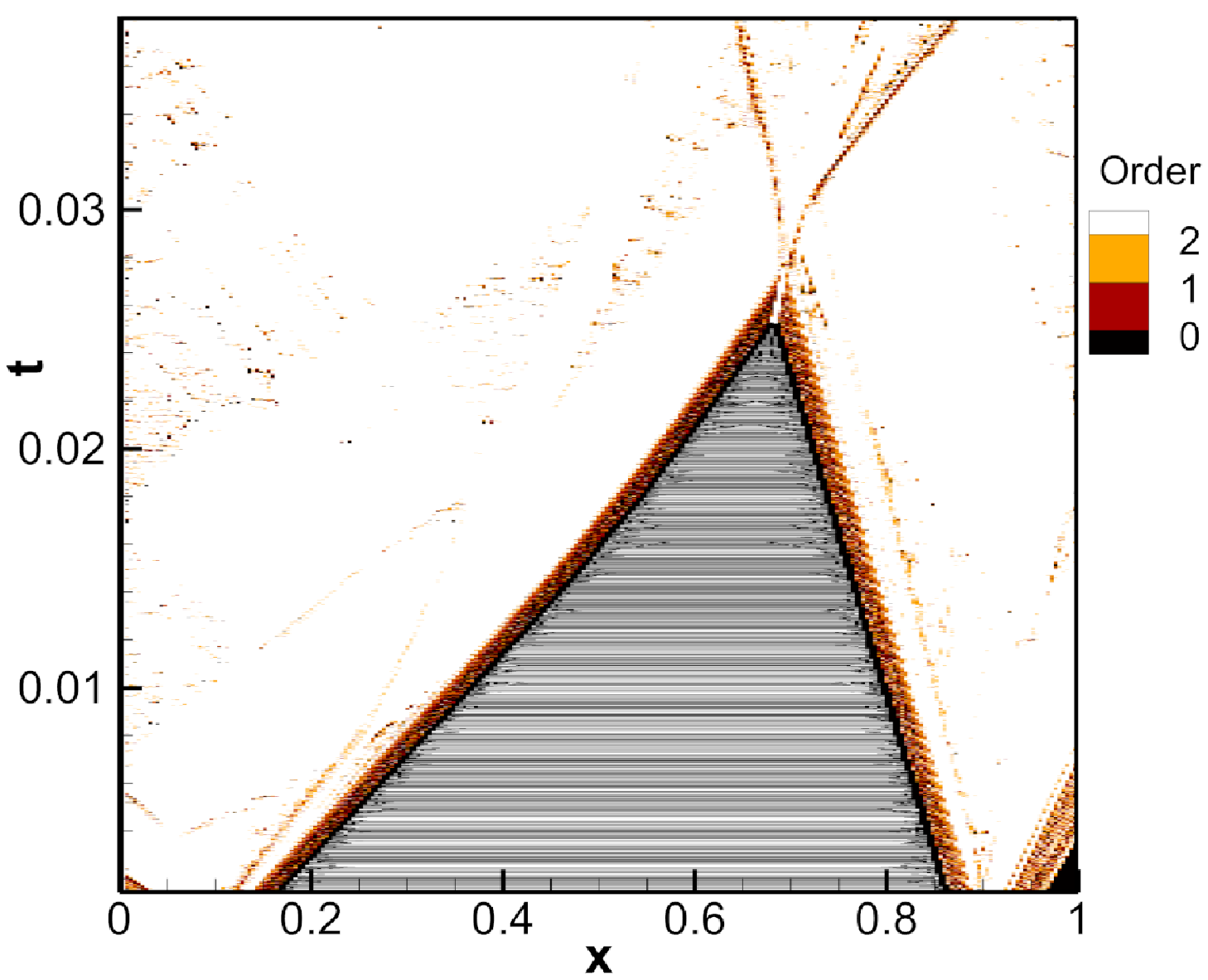}}
  \subfigure[$k=4$]{
  \includegraphics[width=5.5 cm]{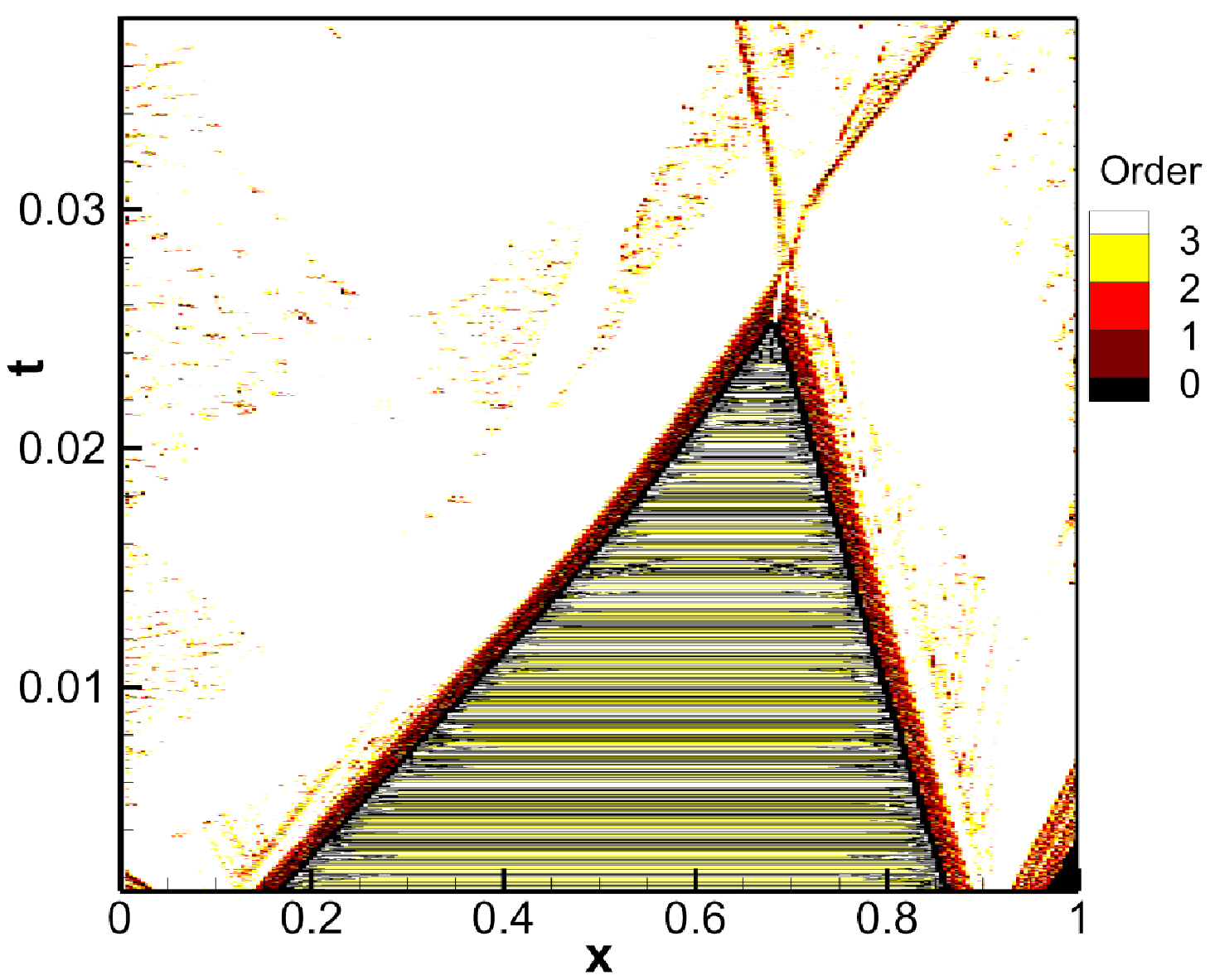}}
  \subfigure[$k=5$]{
  \includegraphics[width=5.5 cm]{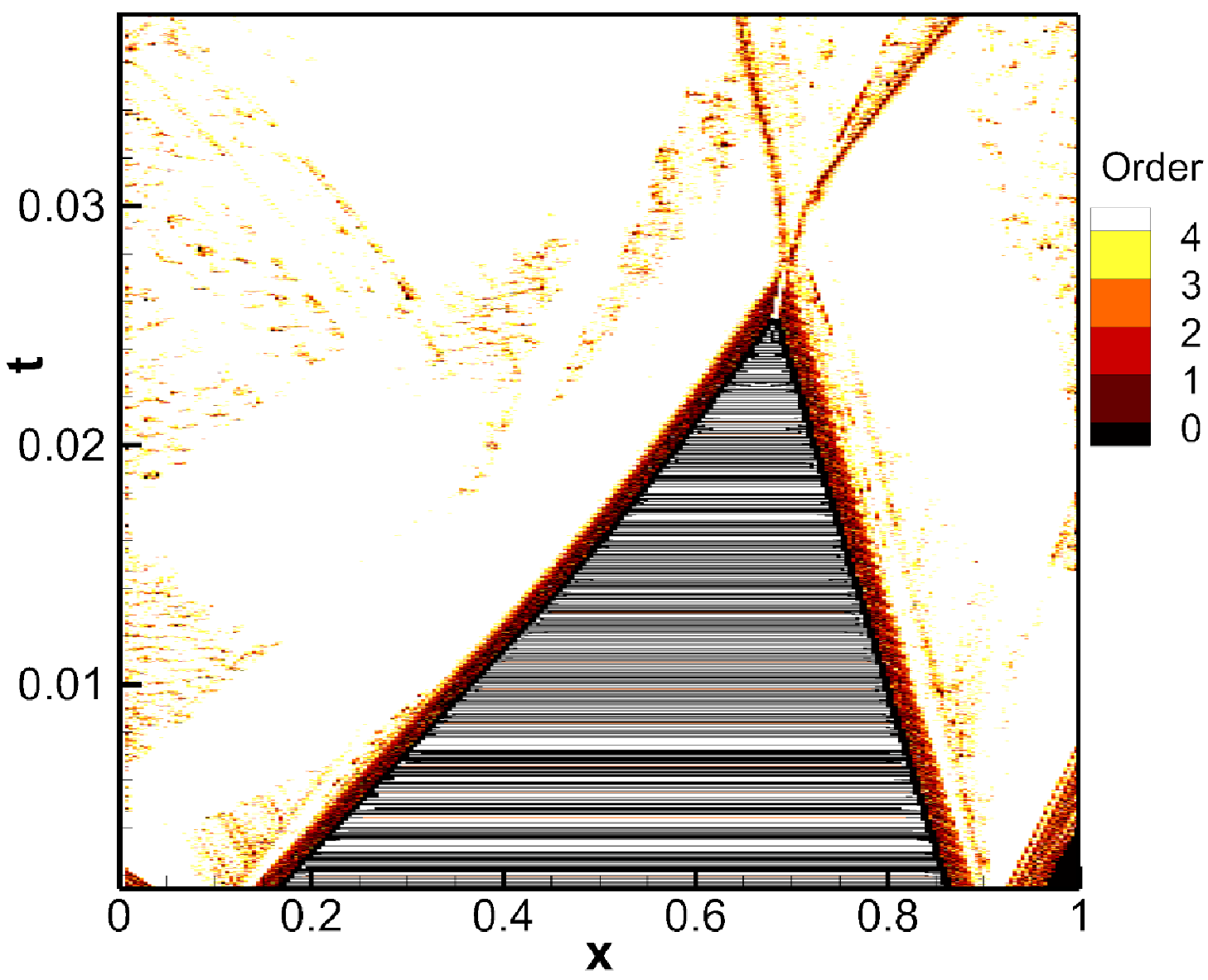}}
  \subfigure[$k=6$]{
  \includegraphics[width=5.5 cm]{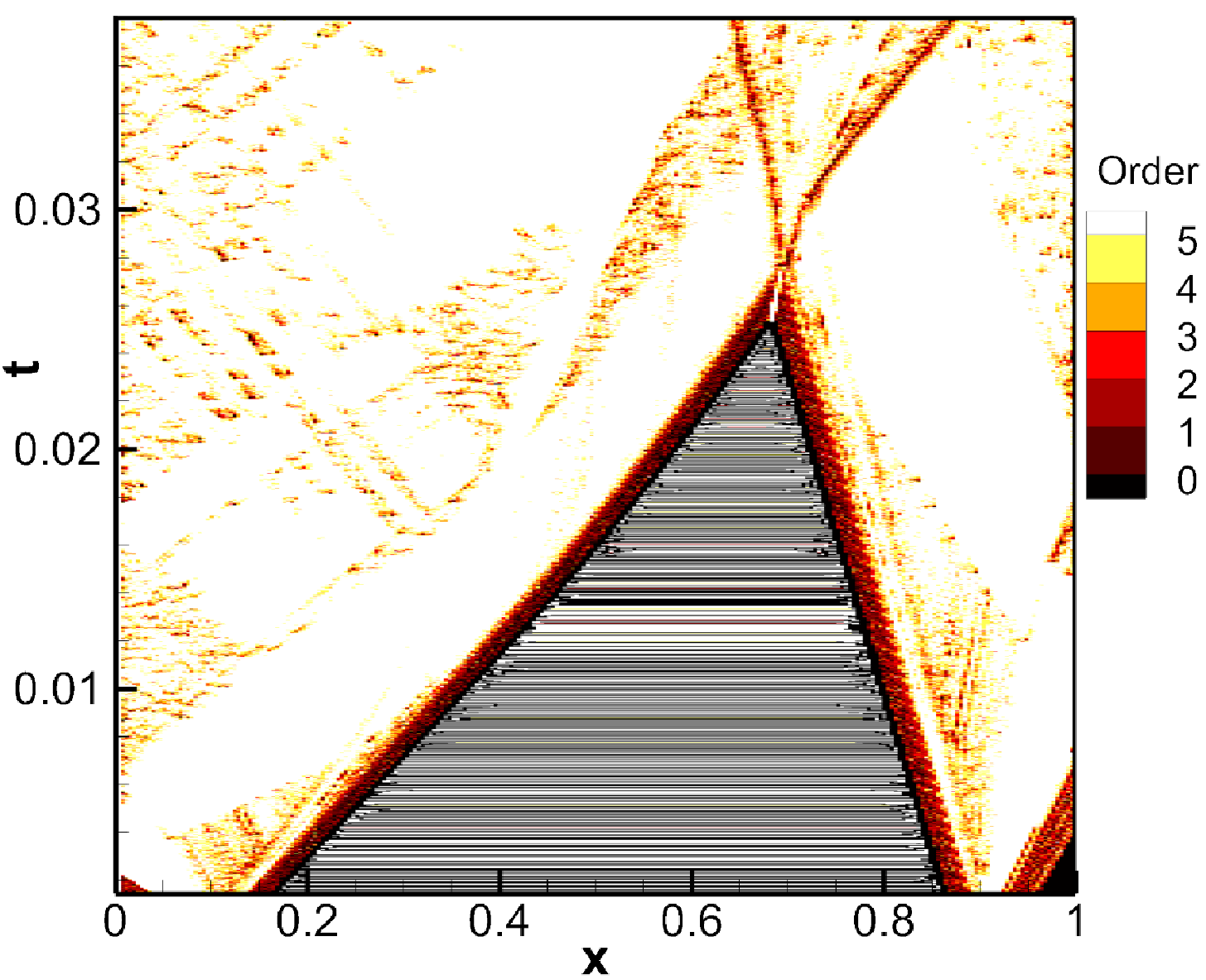}}
  \caption{The Blast wave problem computed by RKDG schemes with the MR limiter using 250 cells. 
  The first and second rows: distributions of the density and the polynomial order at $t=0.038$;
   The third and fourth rows: the time history of the polynomial order, 
  the white parts represent the original $k$th-order DG polynomial. }
 \label{FIG:Blast}
 \end{figure}

 \paragraph{Example 4.2.6} The blast wave problem.
 Two blast waves are initialized as
 \begin{equation*}
   (\rho,u,p)=\begin{cases}
                (1,0,1000),  0\le x<0.1, \\
                (1,0,0.01), 0.1\le x< 0.9,\\
                (1,0,100), 0.9\le x\le 1,\\
              \end{cases}
  \end{equation*}
in the computational domain $[-1,1]$ with
reflective boundary conditions are implemented on both sides.
 When the simulation starts, the two blast waves move towards each other and interact after the collision.
 It is difficult for numerical methods to resolve the peak value appearing 
 at the center of the interaction zone.
We compute this problem with the RKDG schemes by using 250 cells, 
which is a very coarse mesh compared to most simulations in the literature.
Fig. \ref{FIG:Blast} shows the distributions of the density and the polynomial order at $t=1.8$
and the time history of the polynomial order.
The reference solution is calculated by the fifth-order WENO scheme with $h=1/5000$.
We observe that the result of $k=1$ with the MR limiter can 
even compete with most results computed by higher-order schemes with 400 cells
in the literature, see for example Fig. 3.12 in \cite{Fu2017NewLimiter},
Fig. 8 in \cite{Wei2024HybridLimiter}, and Fig. 3.2 in \cite{Wei2025JumpFilter}.
When $k\ge2$, the MR limiter can almost recover the peak value.

\subsection{Numerical examples for the 2D compressible Euler equations}
The 2D compressible Euler equations read
\begin{equation}
 \frac{\partial \mathbf{U}}{\partial t}+\frac{\partial \mathbf{F}}{\partial x}+\frac{\partial \mathbf{G}}{\partial y}=0,
\end{equation}
where $\mathbf{U}=[\rho,\rho u,\rho v,\rho e]^T$, $\mathbf{F}=[{\rho u,\rho u^2+p,\rho uv,(\rho e+p)u}]^T$
and $\mathbf{G}=[{\rho v,\rho uv,\rho v^2+p,(\rho e+p)v}]^T$
with $\rho$, $u$, $v$, $p$, and $e$ respectively representing the density, $x$-velocity, $y$-velocity, pressure, and specific total energy.
The specific total energy is calculated as $e=\frac{p}{\rho(\gamma-1)}+\frac{1}{2} (u^2+v^2)$,
where $\gamma$ is the specific heat ratio and is set to 1.4 for all the following 2D examples.

\paragraph{Example 4.3.1} The double Mach reflection problem.

This example was originally proposed by Woodward and Colella \cite{Woodward1984JCP}
and became a benchmark to test the shock-capturing capability and high fidelity of high-order schemes.
This problem describes the double Mach reflection phenomenon induced by a Mach 10 shock passing through a $30^\circ$ wedge.
First, we consider an equivalent problem used in many other papers.
The initial condition is given by
\begin{equation*}
 \left(\rho, u, v, p\right)=
   \begin{cases}
     \left(1.4, 0, 0, 1\right) & \text{if } x>\frac{1}{6}+\frac{y}{\sqrt{3}}, \\
     \left(8, 8.25\sin(60^\circ), -8.25\cos(60^\circ), 116.5\right), & \text{otherwise},
   \end{cases}
\end{equation*}
in the computational domain $[0, 4]\times[0, 1]$,
which describes a Mach 10 oblique shock 
inclined at an angle of $60^\circ$ to the horizontal direction.
The post-shocked states are enforced on the left boundary.
Nonreflective boundary conditions are applied to the right boundary.
The exact motion of the Mach 10 oblique shock is imposed on the top boundary.
On the bottom, nonreflective boundary conditions are imposed for $x\le\frac{1}{6}$ and 
reflective wall boundary conditions are imposed for $x>\frac{1}{6}$.

\begin{figure}[htbp]
  \centering
  \subfigure[$k=1$]{
  \includegraphics[height=3 cm]{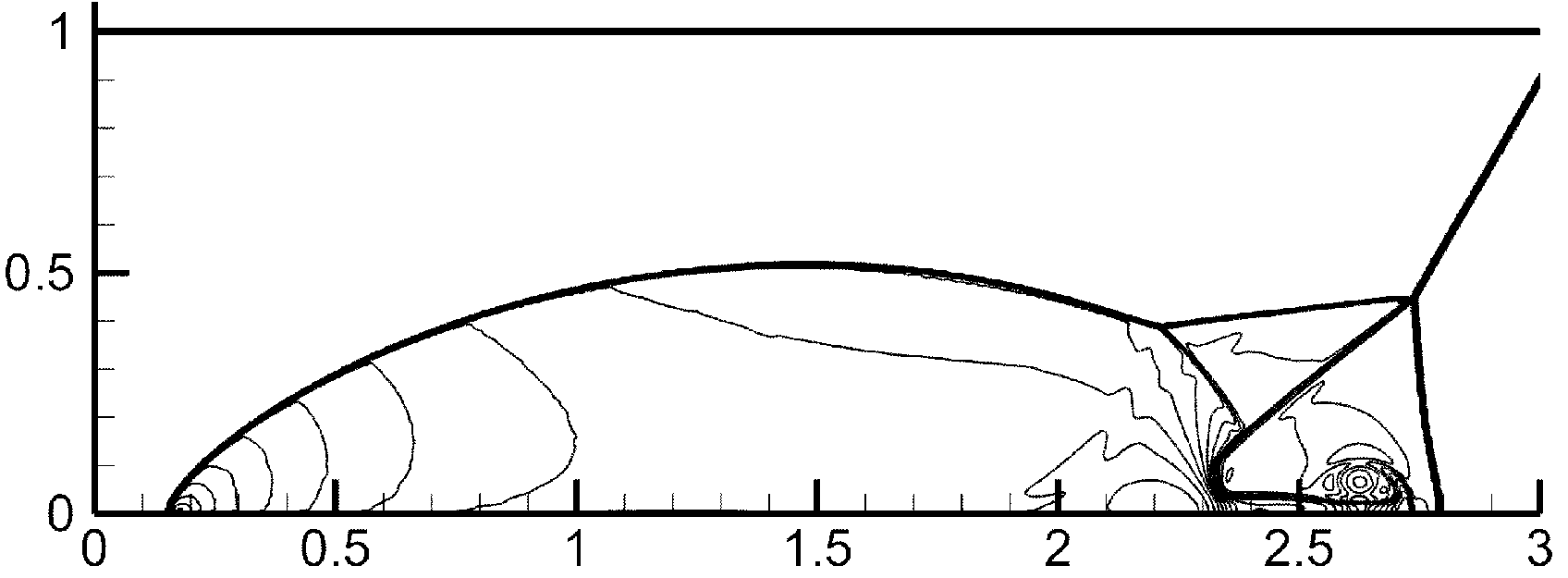}}
  \subfigure[$k=1$]{
  \includegraphics[height=3 cm]{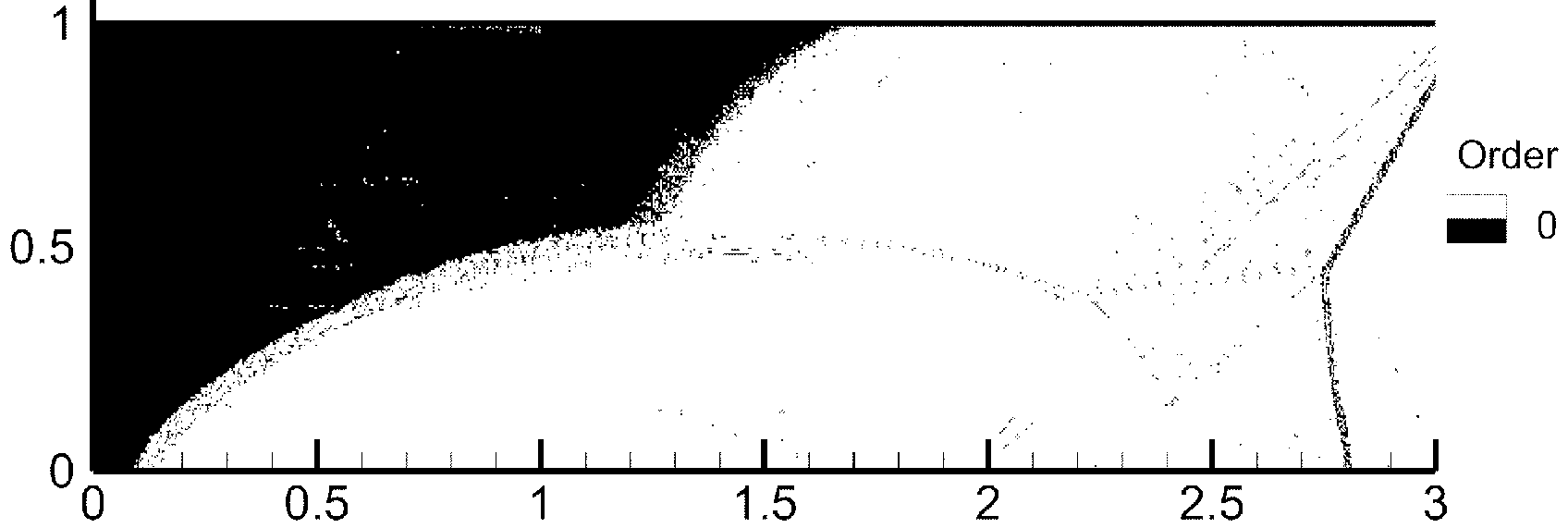}}
  \subfigure[$k=2$]{
  \includegraphics[height=3 cm]{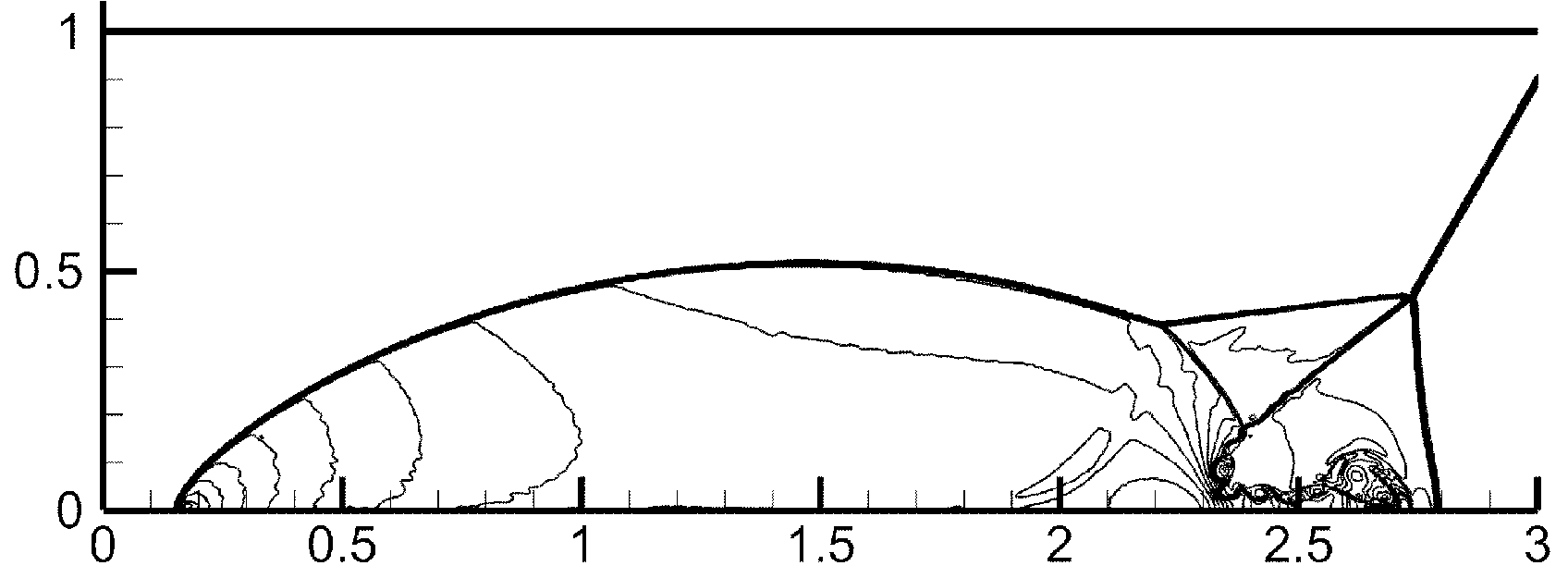}}
  \subfigure[$k=2$]{
  \includegraphics[height=3 cm]{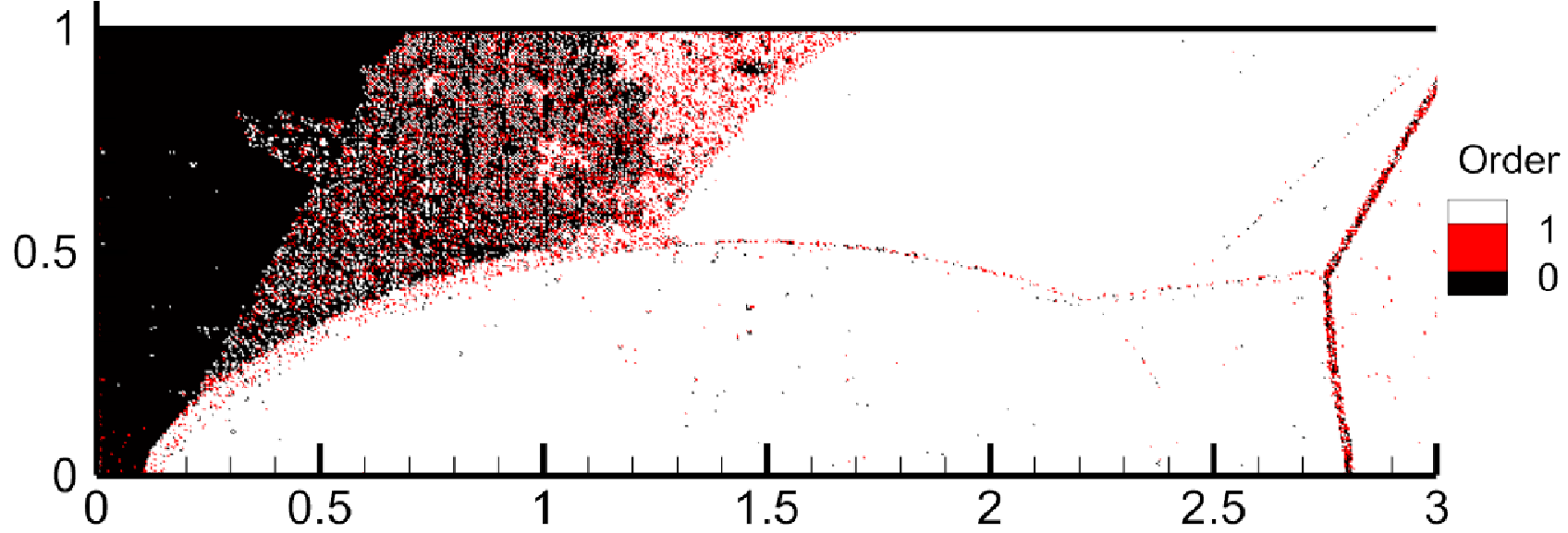}}
  \subfigure[$k=3$]{
  \includegraphics[height=3 cm]{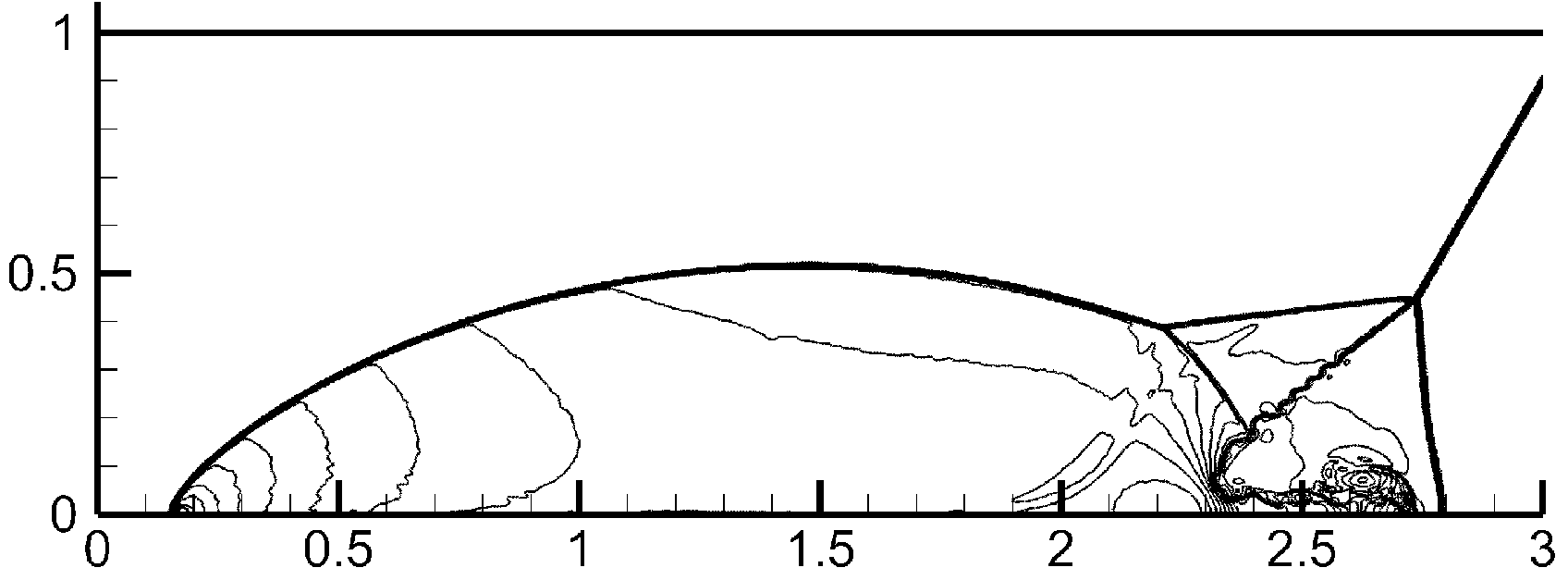}}
  \subfigure[$k=3$]{
  \includegraphics[height=3 cm]{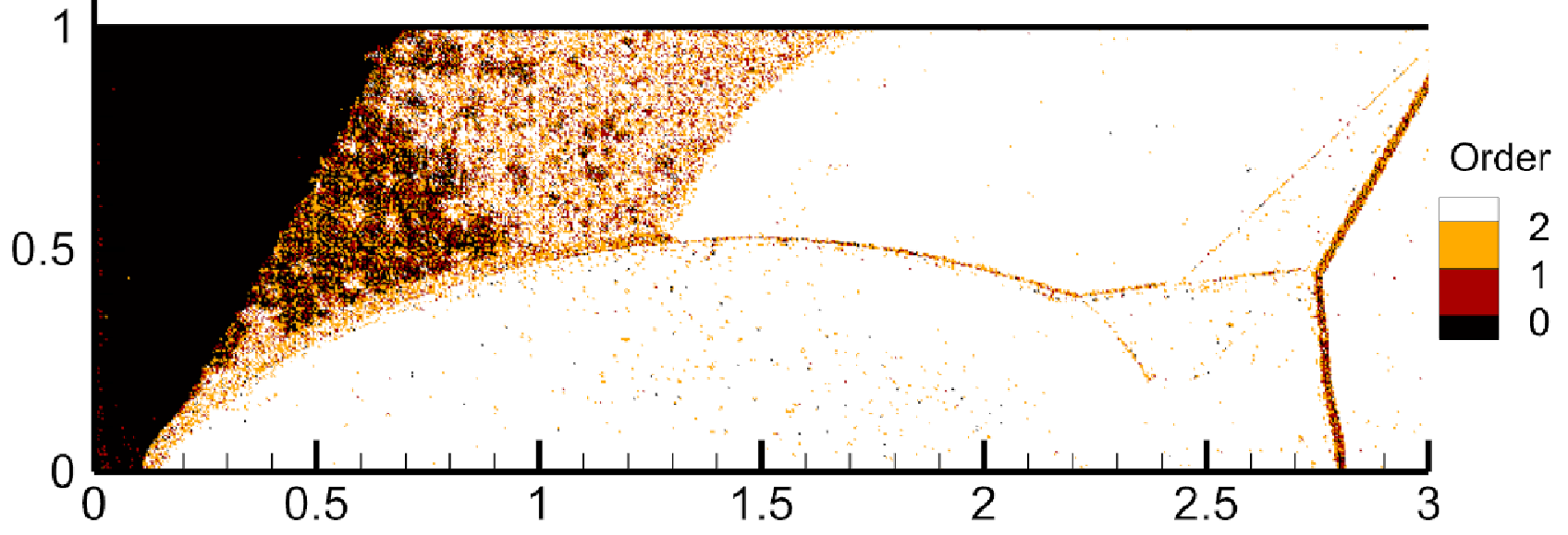}}
  \caption{Density contours (left column) and the polynomial order (right column) of the double Mach reflection problem at $t=0.2$ calculated 
  by RKDG schemes with the MR limiter using $960\times240$ uniform quadrilateral cells.
  The density contours contain 30 equidistant contours from 1.5 to 21.5. 
  The white parts of the polynomial order plots represent the original $k$th-order DG polynomial.}
 \label{FIG:DMR_Qua}
 \end{figure}

\begin{figure}[htbp]
  \centering
  \subfigure[$k=1$]{
  \includegraphics[height=3 cm]{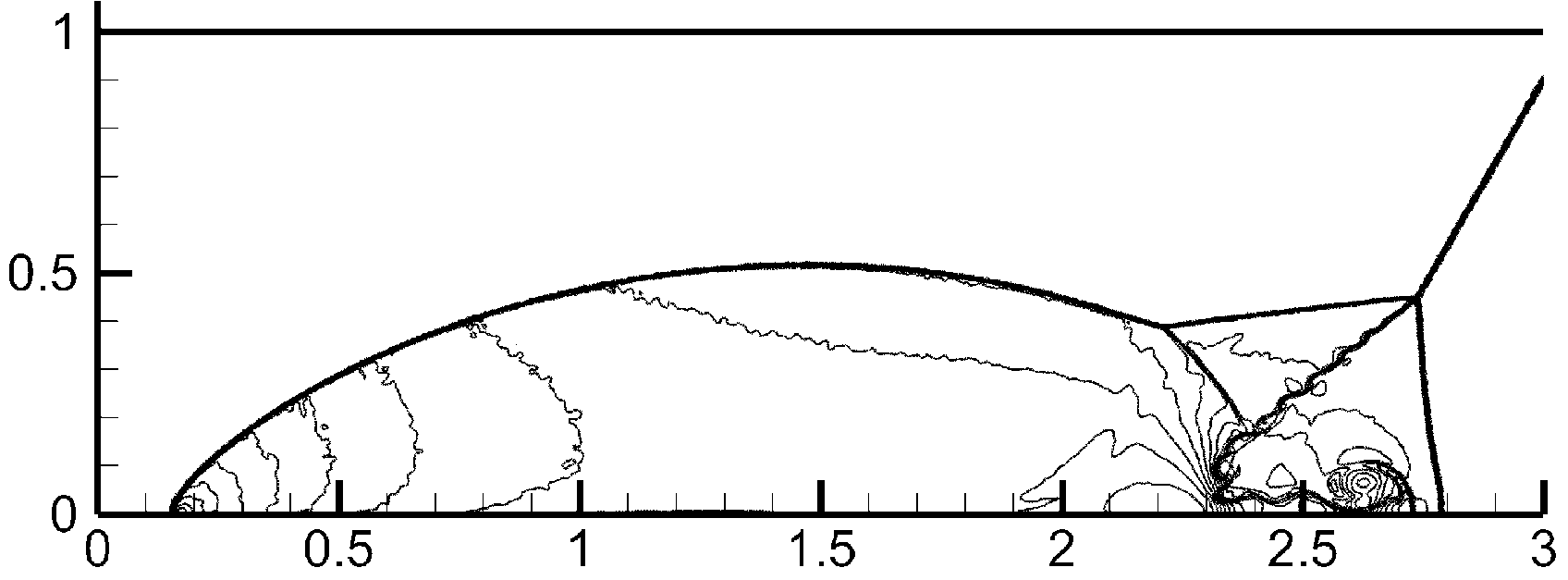}}
  \subfigure[$k=1$]{
  \includegraphics[height=3 cm]{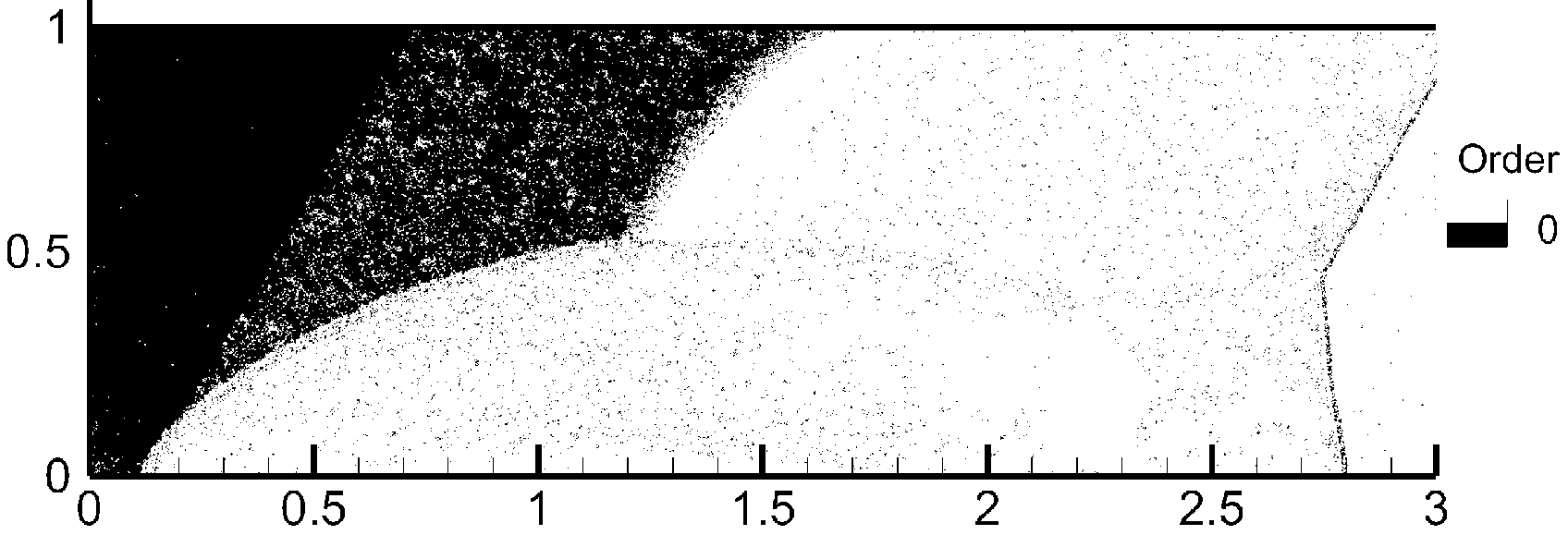}}
  \subfigure[$k=2$]{
  \includegraphics[height=3 cm]{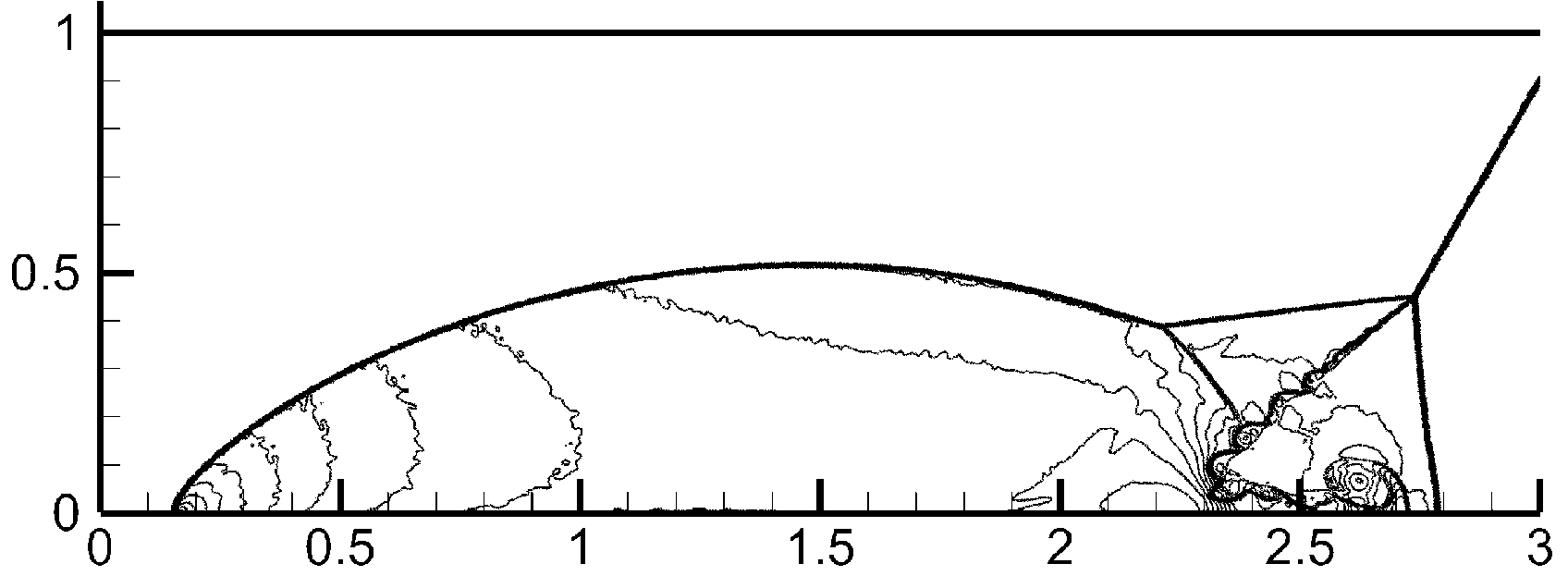}}
  \subfigure[$k=2$]{
  \includegraphics[height=3 cm]{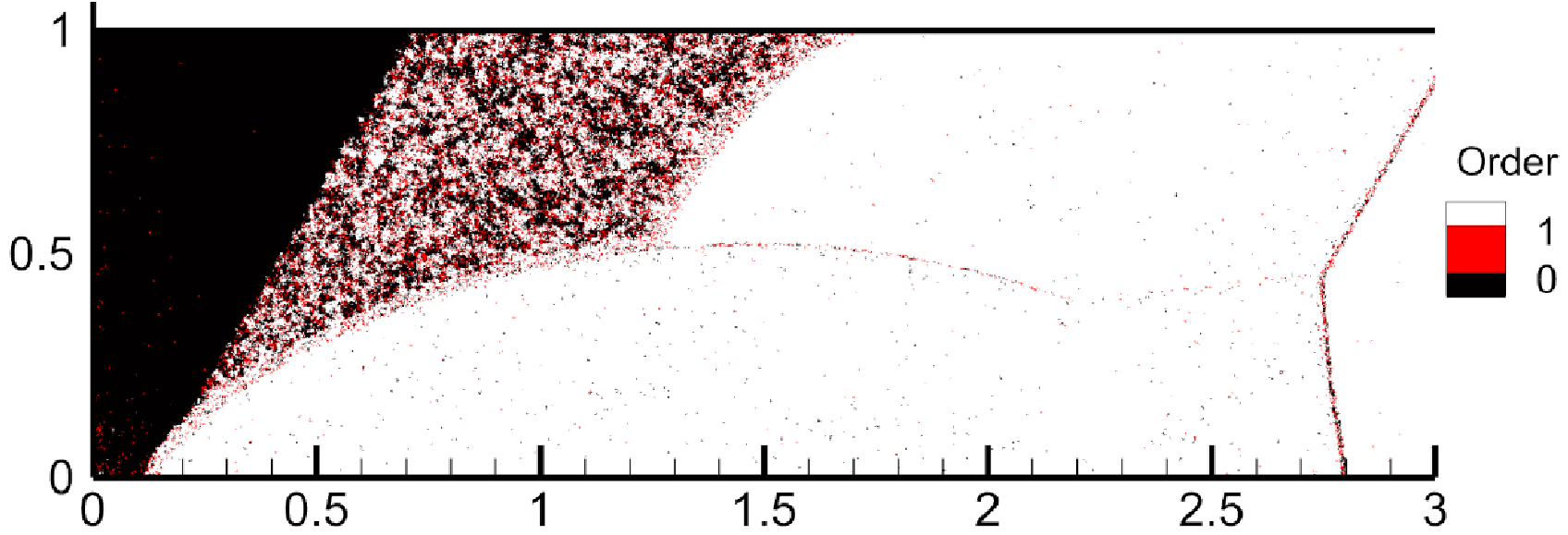}}
  \subfigure[$k=3$]{
  \includegraphics[height=3 cm]{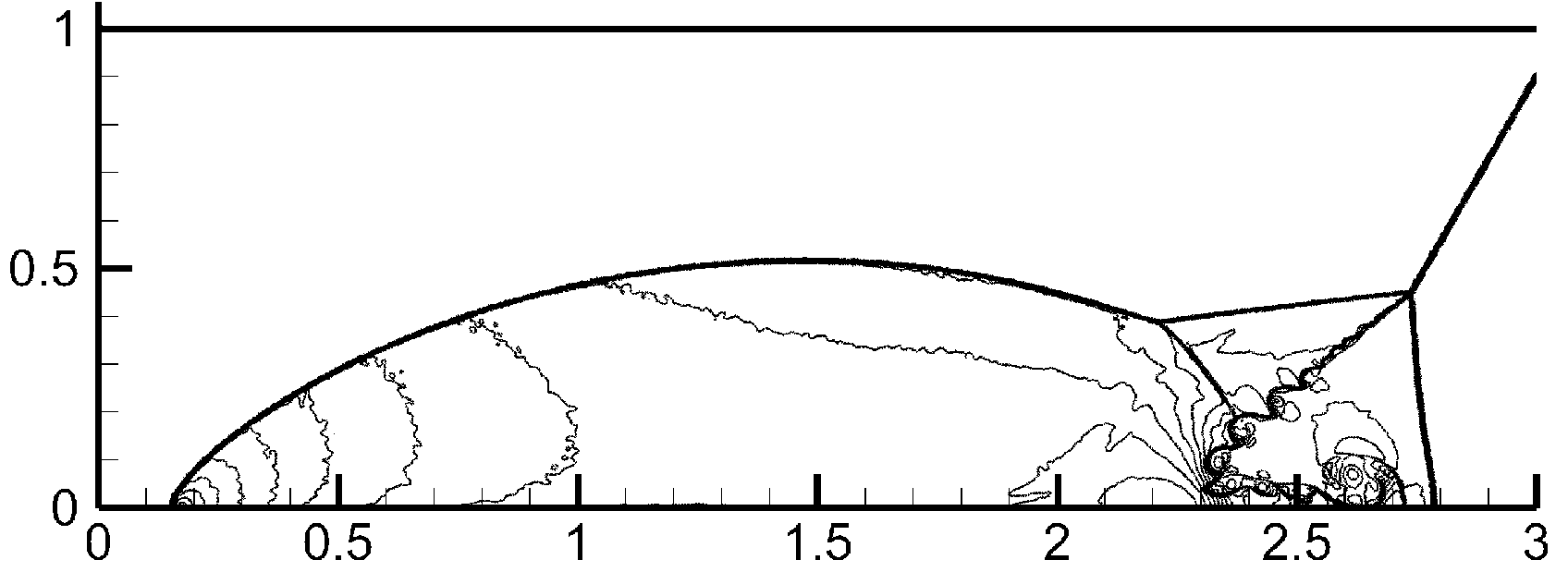}}
  \subfigure[$k=3$]{
  \includegraphics[height=3 cm]{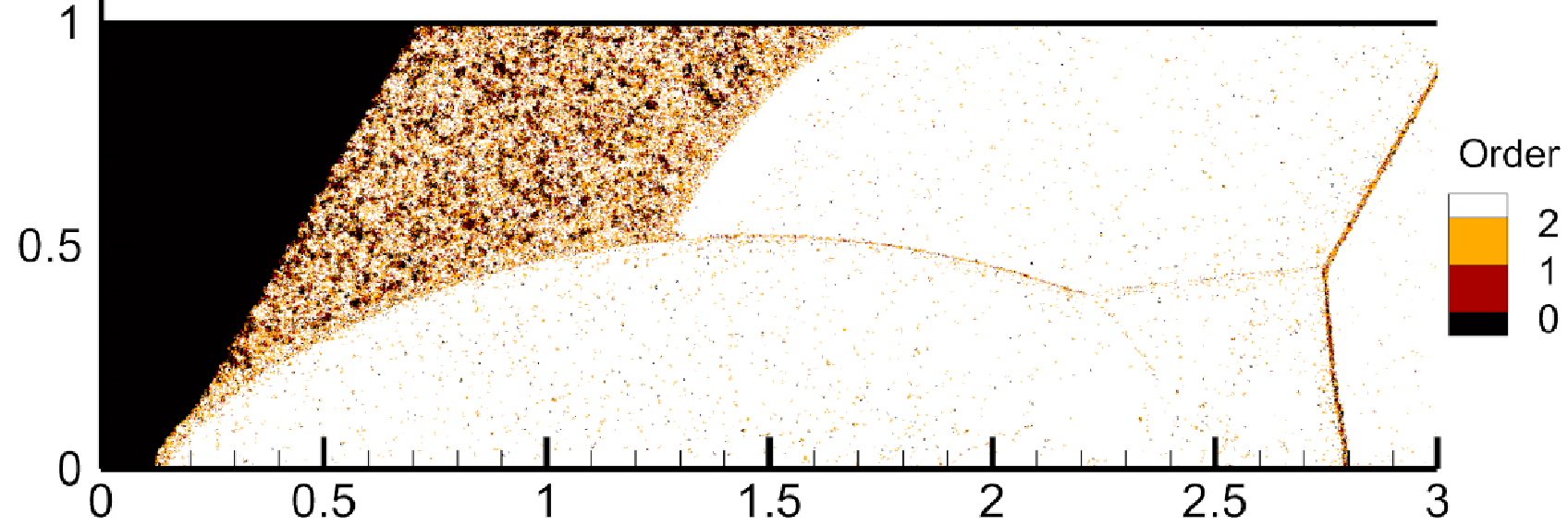}}
  \caption{Density contours (left column) and the polynomial order (right column) of the double Mach reflection problem at $t=0.2$ calculated 
  by RKDG schemes with the MR limiter using quasi-uniform triangular cells 
  with the edge length equal to $1/240$ on the boundaries.
  The density contours contain 30 equidistant contours from 1.5 to 21.5.
  The white parts of the polynomial order plots represent the original $k$th-order DG polynomial.}
 \label{FIG:DMR_Tri}
 \end{figure}

Fig. \ref{FIG:DMR_Qua} shows the density contours and the polynomial order at $t=0.2$ calculated 
by RKDG schemes with the MR limiter using $960\times240$ uniform quadrilateral cells.
Fig. \ref{FIG:DMR_Tri} shows the corresponding results calculated by RKDG schemes with the MR limiter using quasi-uniform triangular cells 
with the edge length equal to $1/240$ on the boundaries of the computational domain.
The MR limiter can detect the discontinuities effectively,
and the DG polynomials do not always reduce to the lowest order near discontinuities.
The MR limiter tends to detect fewer troubled cells on quadrilateral meshes than on triangular meshes,
since discrete solutions on quadrilateral meshes are smoother than those on triangular meshes.
We note that many troubled cells are detected in the uniform flow field by the MR limiter,
because truncation errors lead to `non-uniform' distributions of the density
and the MR limiter can detect discontinuities of all scales.
Nevertheless, these troubled cells do not affect the computational accuracy.
The MR limiter with $k\ge 2$ using the quadrilateral mesh captures some small vortices induced by
the Kelvin-Helmholtz instability in the reflection zone.
Since the effective mesh size of the triangular mesh is smaller than that of the quadrilateral mesh,
the MR limiter with $k=1$ can already capture small vortices on the triangular mesh.
Additionally, the flow fields computed by the MR limiter are non-oscillatory on the quadrilateral mesh and the triangular mesh, 
which demonstrates the non-oscillatory property of the MR limiter.
In fact, it is difficult to strike a balance between high-resolution and non-oscillatory properties of high-order schemes.
For example, the FS limiter (Fig. 3.16 in \cite{Fu2017NewLimiter}), the hybrid limiter (Fig. 11 in \cite{Wei2024HybridLimiter}), 
and the hybrid limiter with jump filter (Fig. 3.5 in \cite{Wei2025JumpFilter}) can resolve small vortices very well, 
but the flow fields contain some obvious noises.

 Next, we compute the standard double Mach reflection problem.
 The computational domain is bounded by segments connecting the points 
 $(0,0)$, $(0.1,0)$, $(2.5,2.4/\sqrt{3})$, $(2.5,2)$, $(0,2)$.
 The initial condition is
 \begin{equation*}
 \left(\rho, u, v, p\right)=
   \begin{cases}
     \left(1.4, 0, 0, 1\right) & \text{if } x>0.05, \\
     \left(8, 8.25, 0, 116.5\right), & \text{otherwise},
   \end{cases}
\end{equation*}
which depicts a right-moving Mach 10 normal shock.
The post-shocked states are enforced on the left boundary.
Nonreflective boundary conditions are applied to the right boundary.
The exact motion of the Mach 10 normal shock is imposed on the top boundary.
Reflective wall boundary conditions are imposed on the bottom and the wedge.

\begin{figure}[htbp]
  \centering
  \subfigure[$k=1$]{
  \includegraphics[height=6 cm]{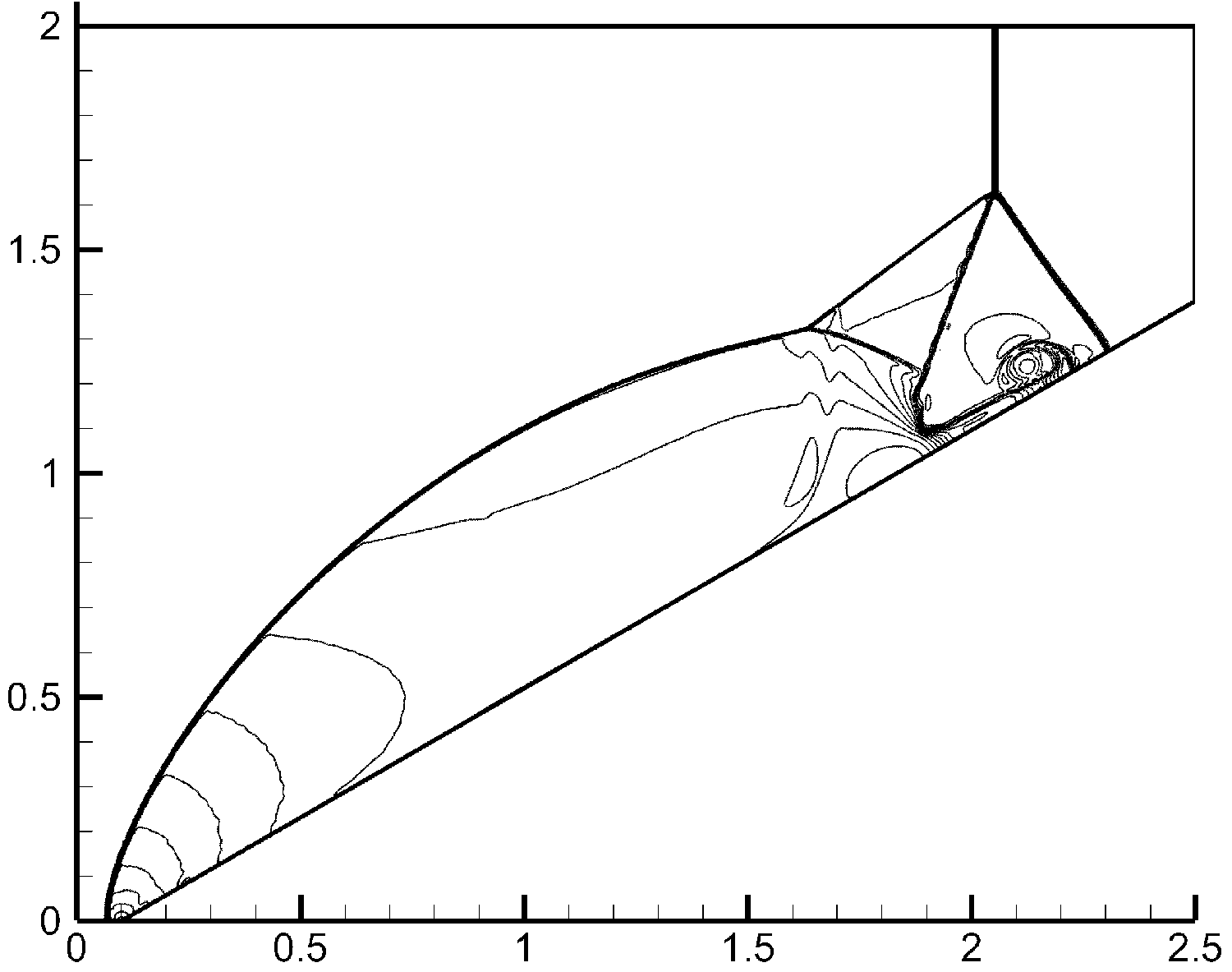}}
  \subfigure[$k=1$]{
  \includegraphics[height=6 cm]{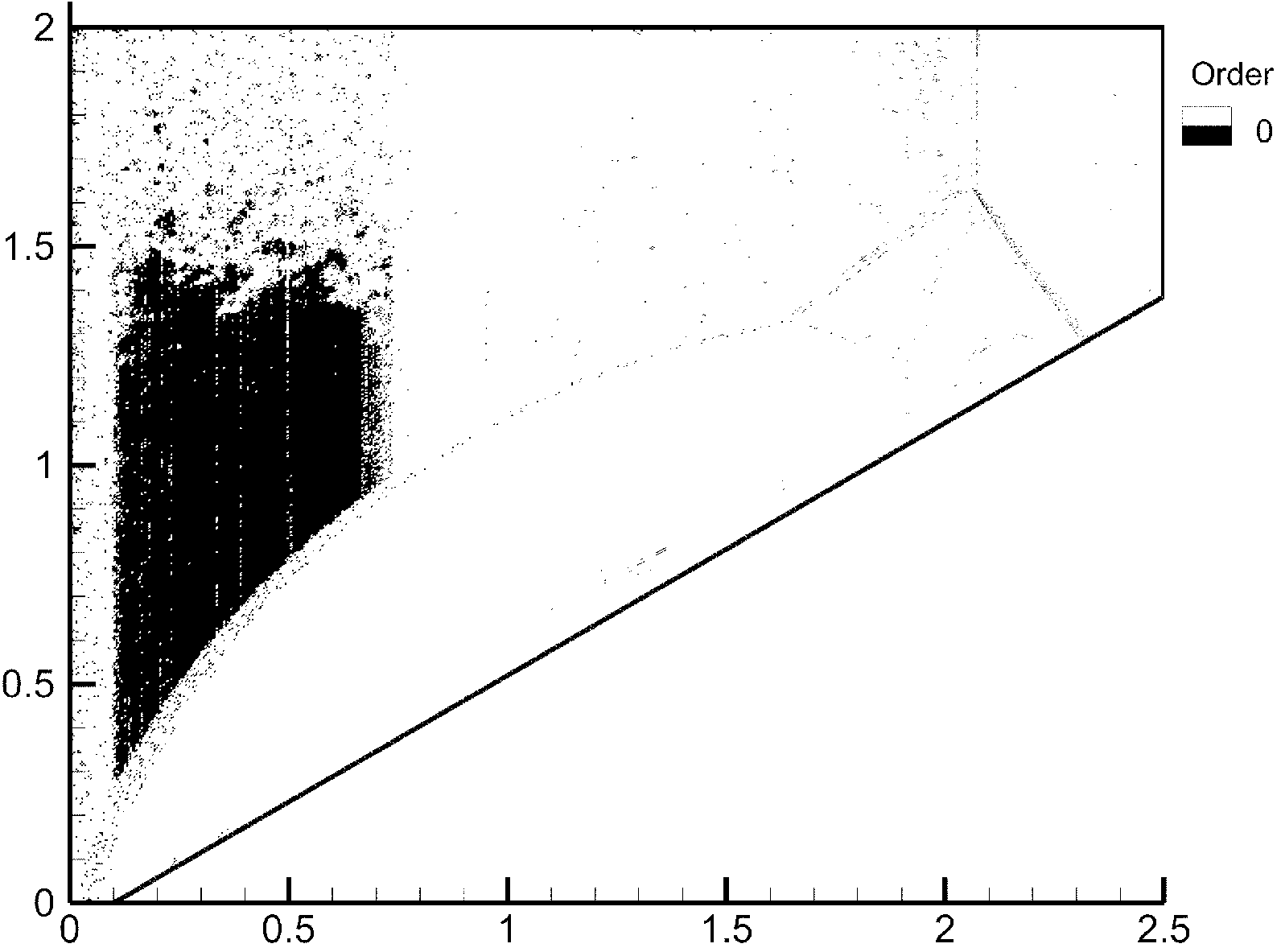}}
  \subfigure[$k=2$]{
  \includegraphics[height=6 cm]{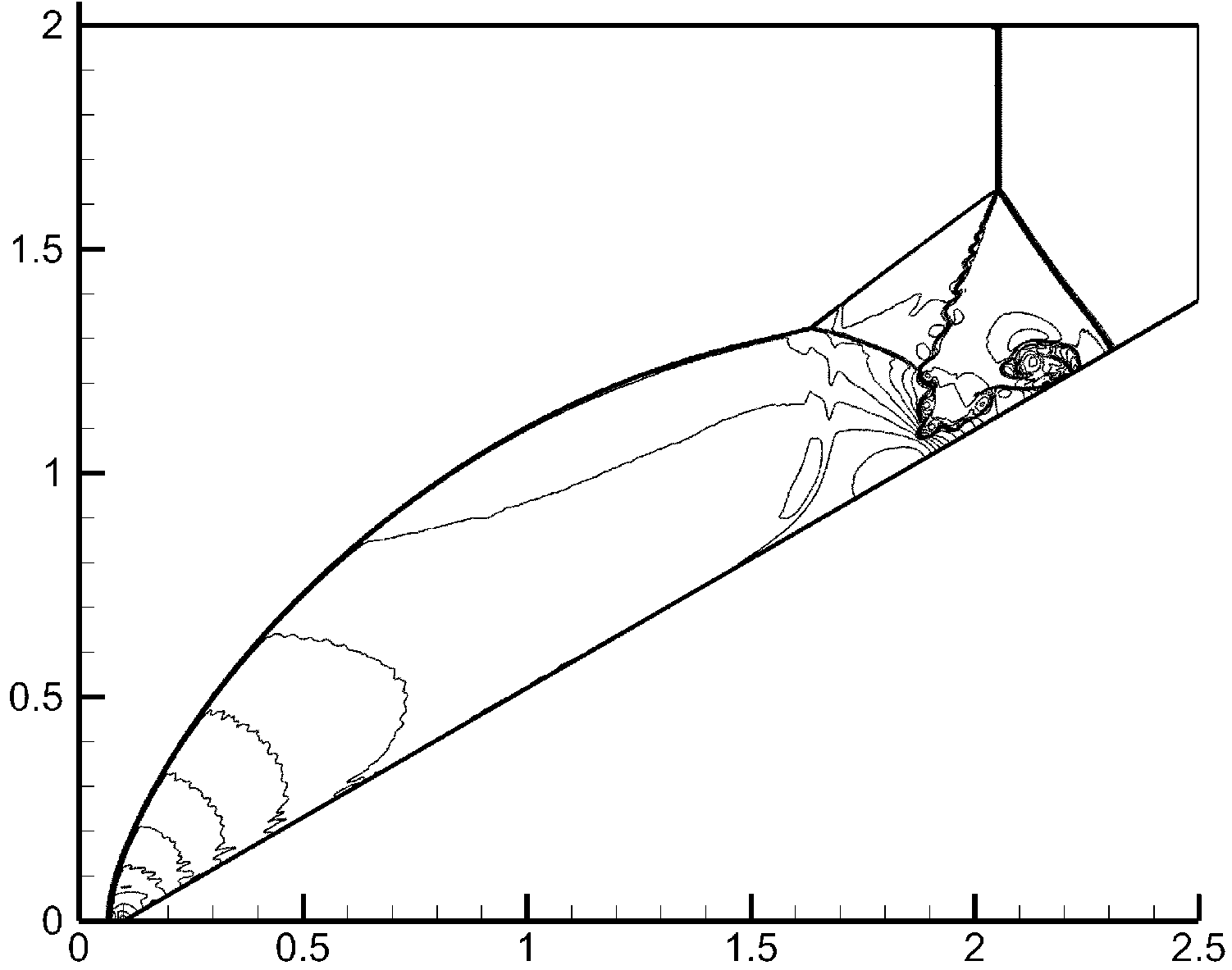}}
  \subfigure[$k=2$]{
  \includegraphics[height=6 cm]{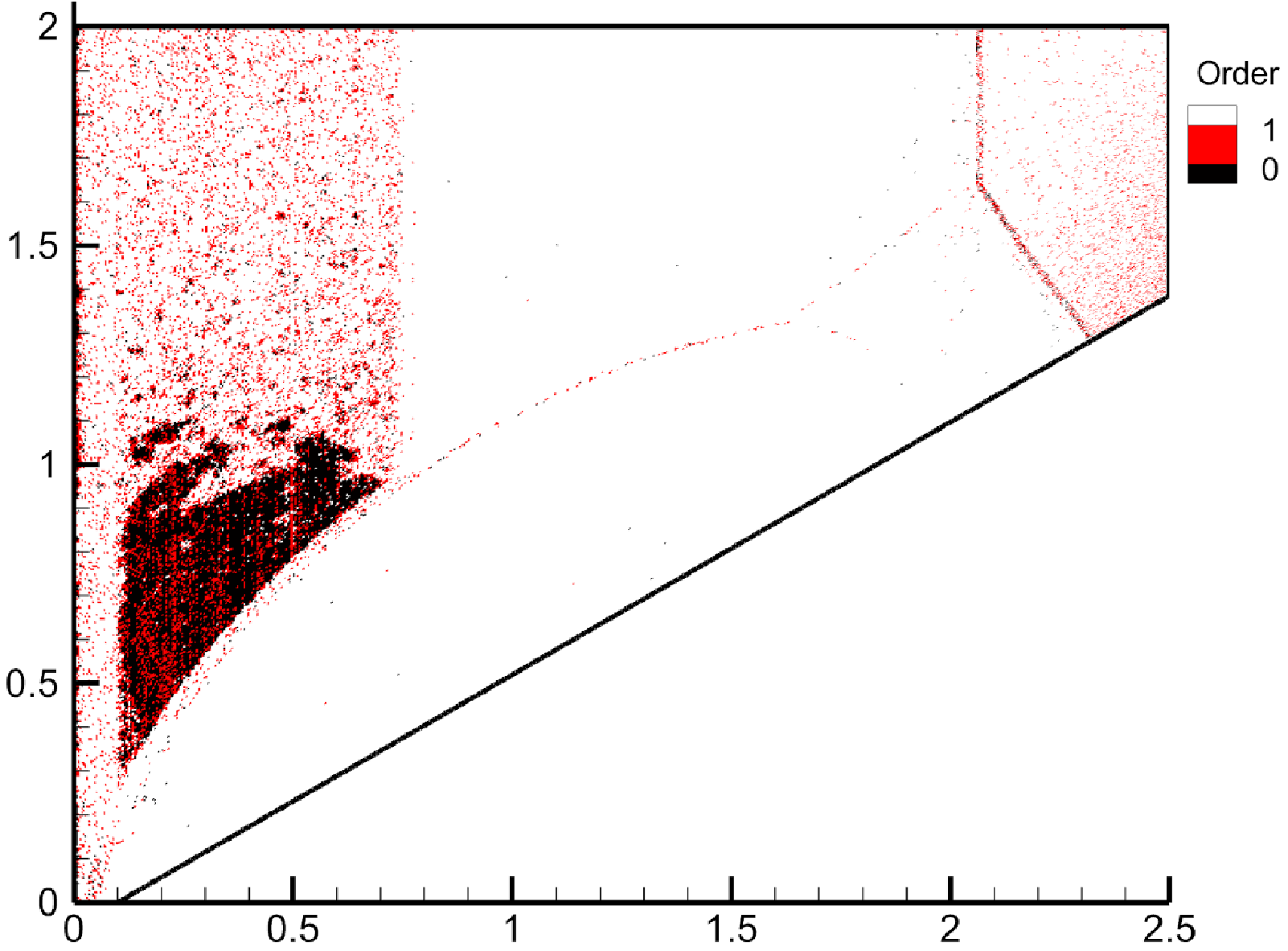}}
  \subfigure[$k=3$]{
  \includegraphics[height=6 cm]{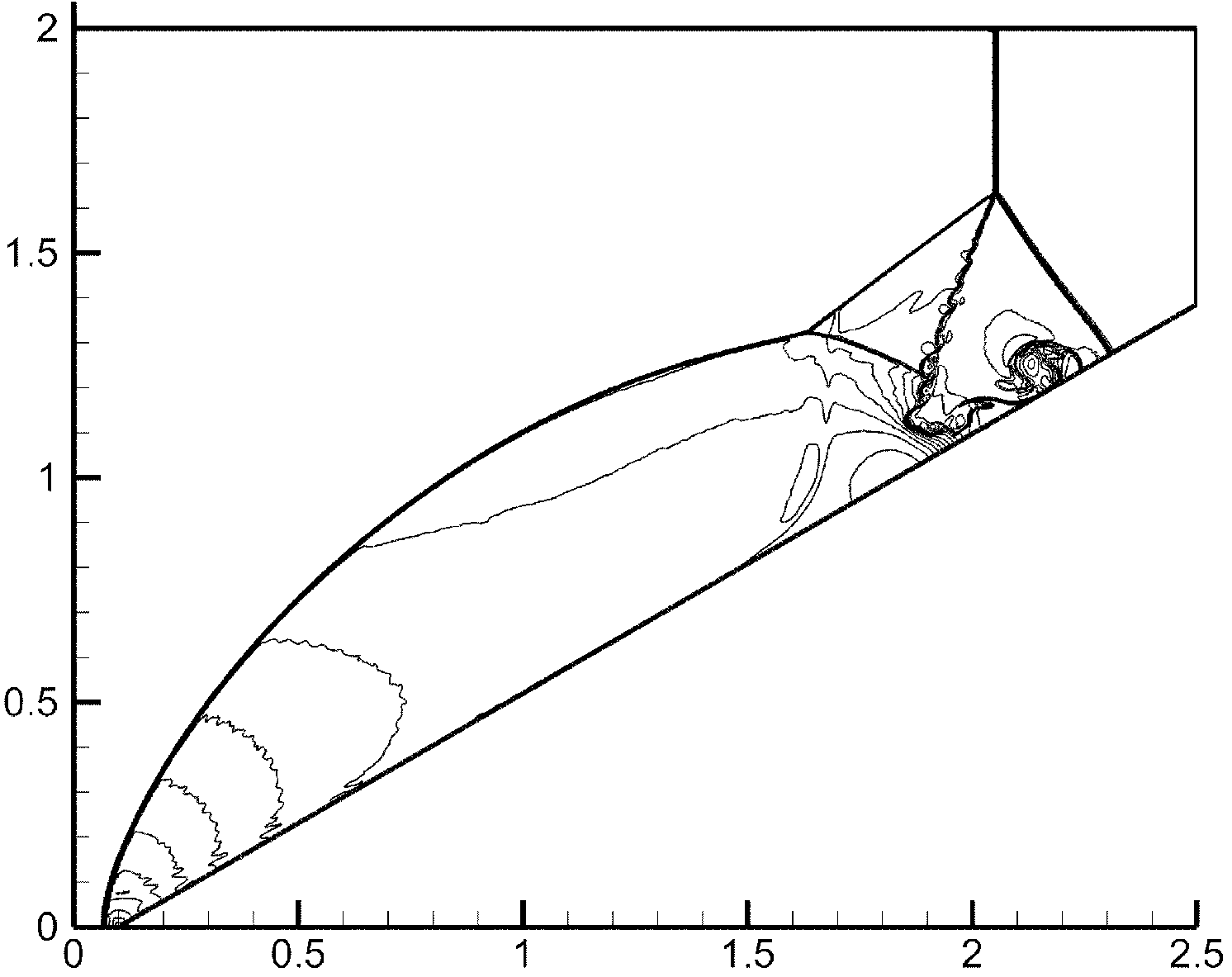}}
  \subfigure[$k=3$]{
  \includegraphics[height=6 cm]{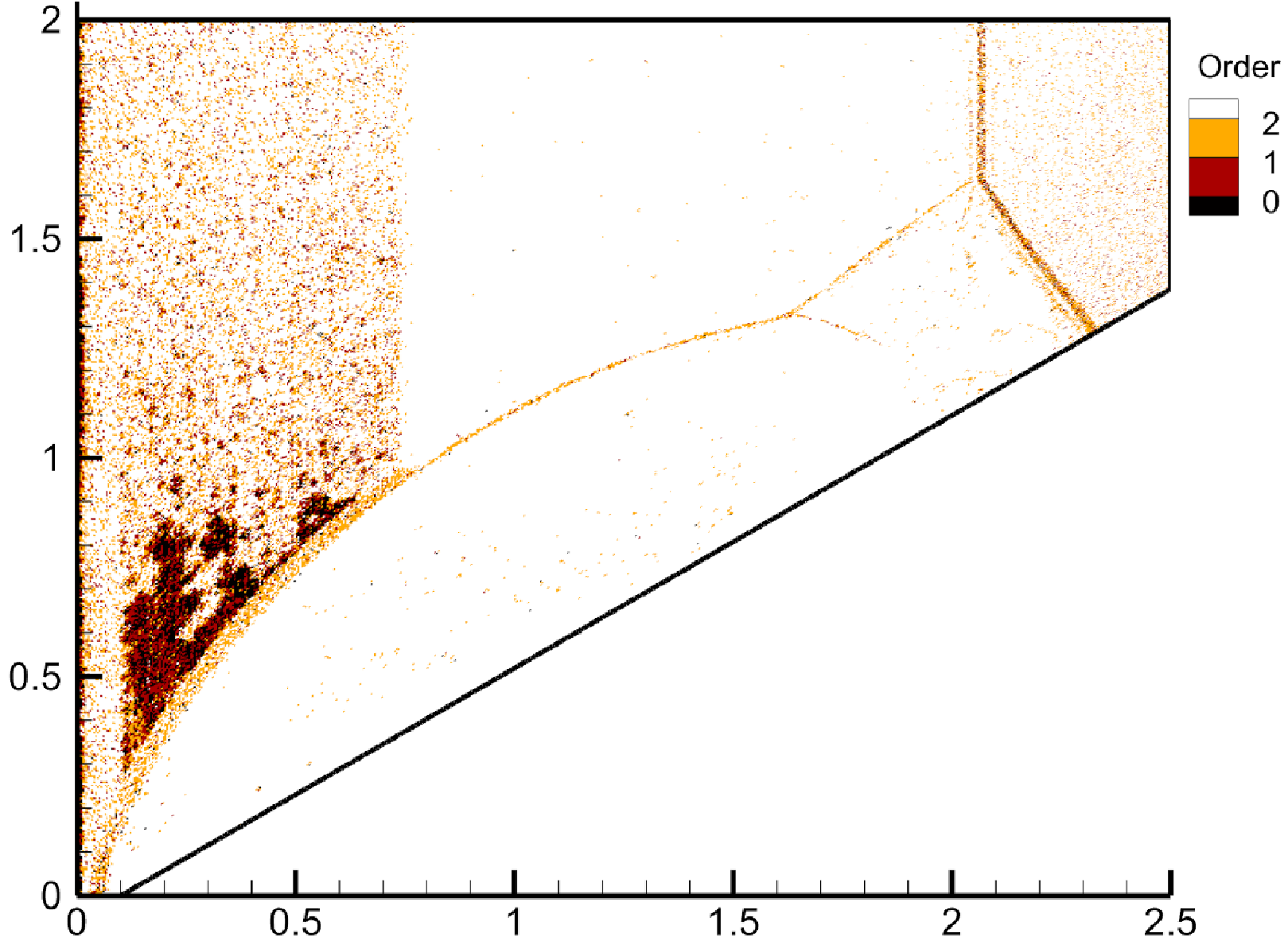}}
  \caption{Density contours (left column) and the polynomial order (right column) of the standard double Mach reflection problem at $t=0.2$ calculated 
  by RKDG schemes with the MR limiter using $600\times480$ structured quadrilateral cells.
  The density contours contain 30 equidistant contours from 1.5 to 21.5.
  The white parts of the polynomial order plots represent the original $k$th-order DG polynomial.}
 \label{FIG:Standard_DMR_Qua}
 \end{figure}

\begin{figure}[htbp]
  \centering
  \subfigure[$k=1$]{
  \includegraphics[height=6 cm]{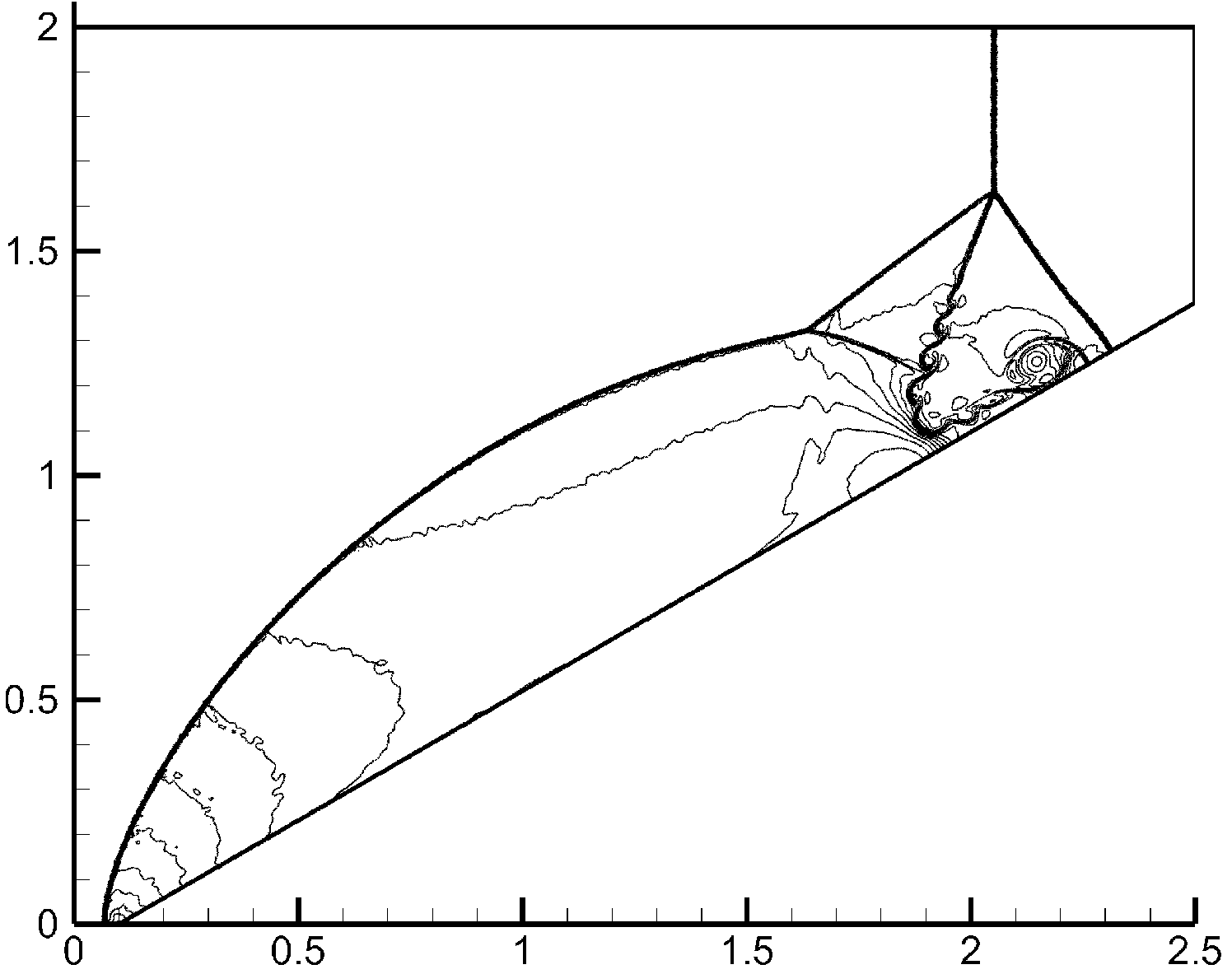}}
  \subfigure[$k=1$]{
  \includegraphics[height=6 cm]{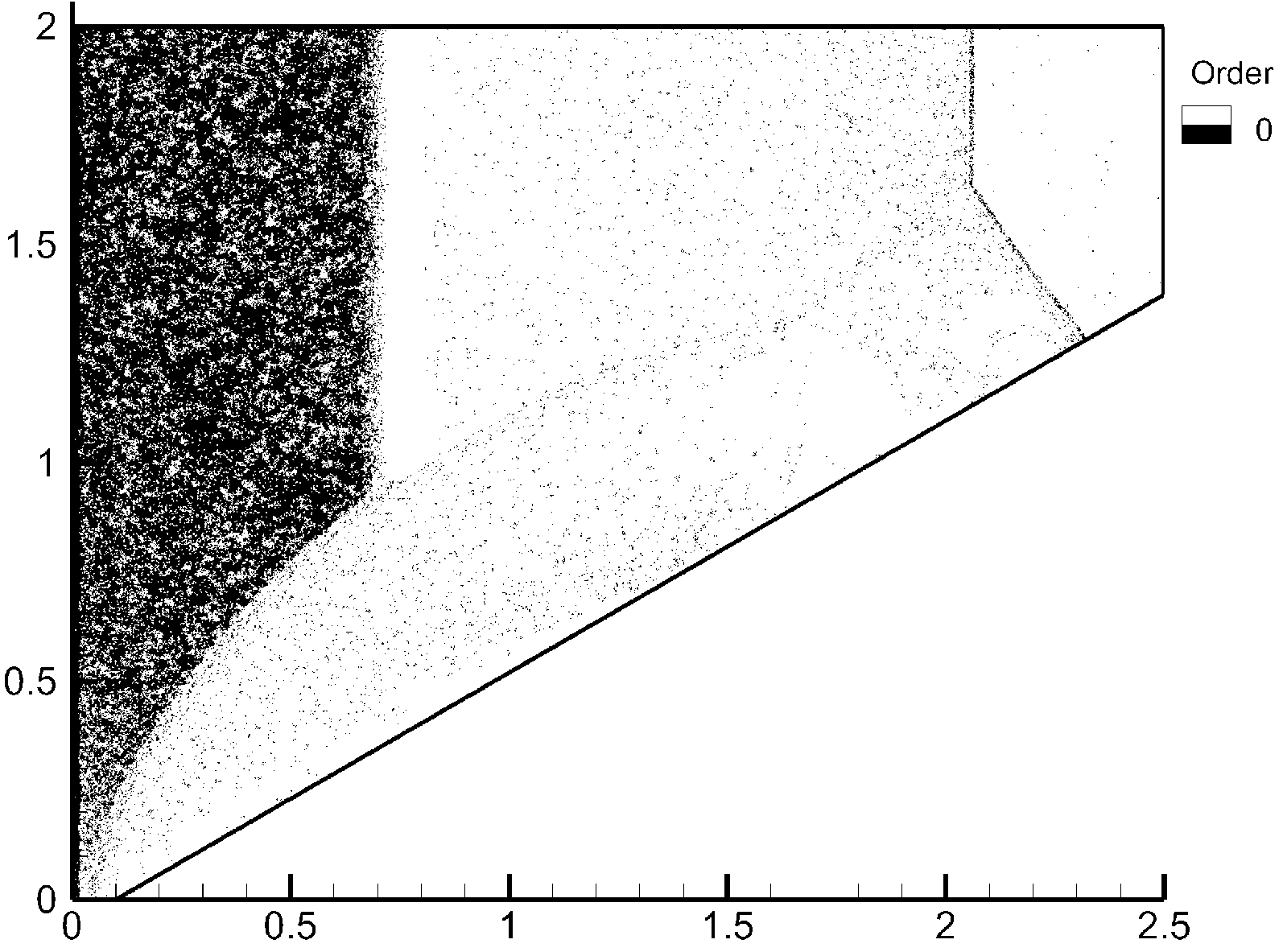}}
  \subfigure[$k=2$]{
  \includegraphics[height=6 cm]{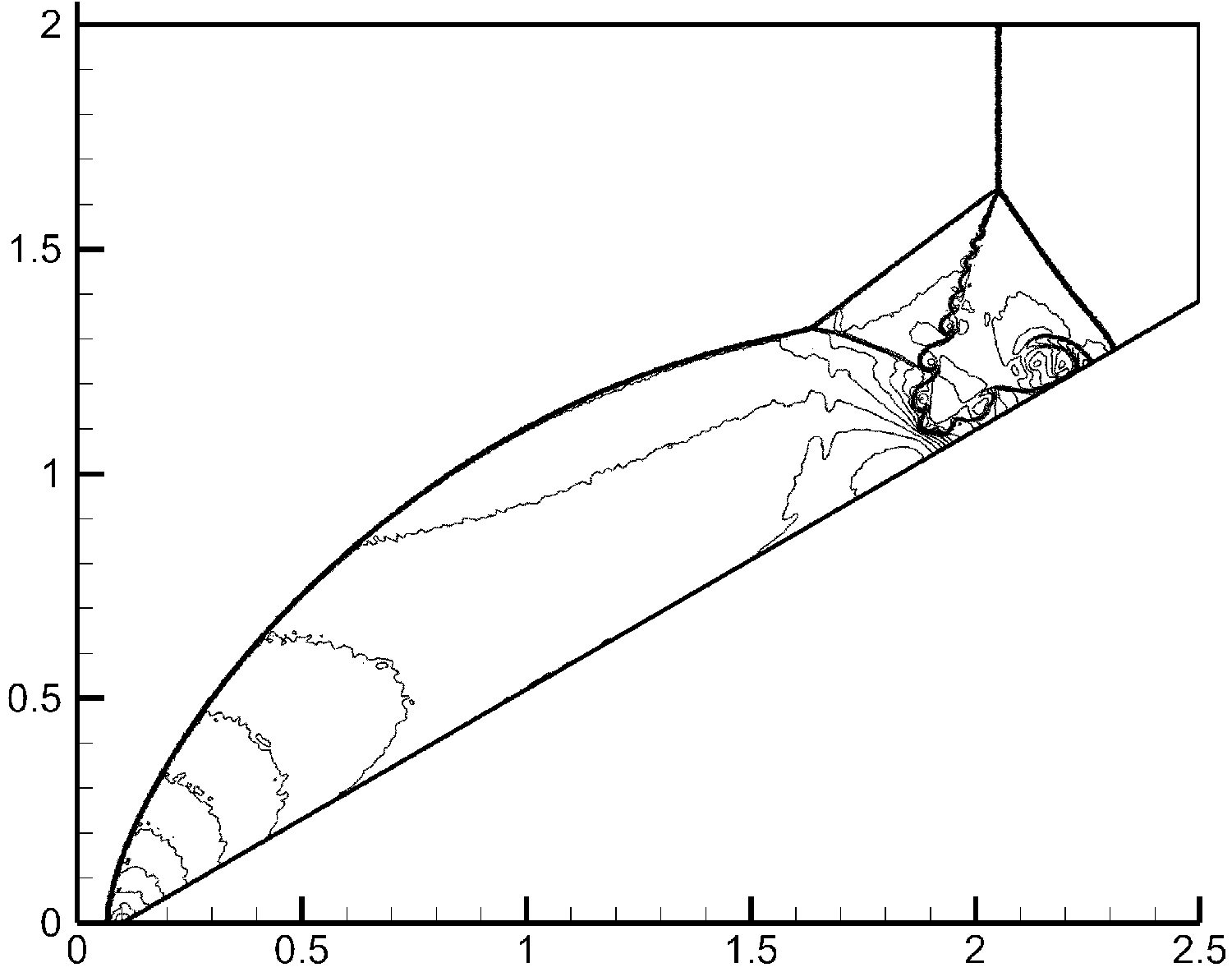}}
  \subfigure[$k=2$]{
  \includegraphics[height=6 cm]{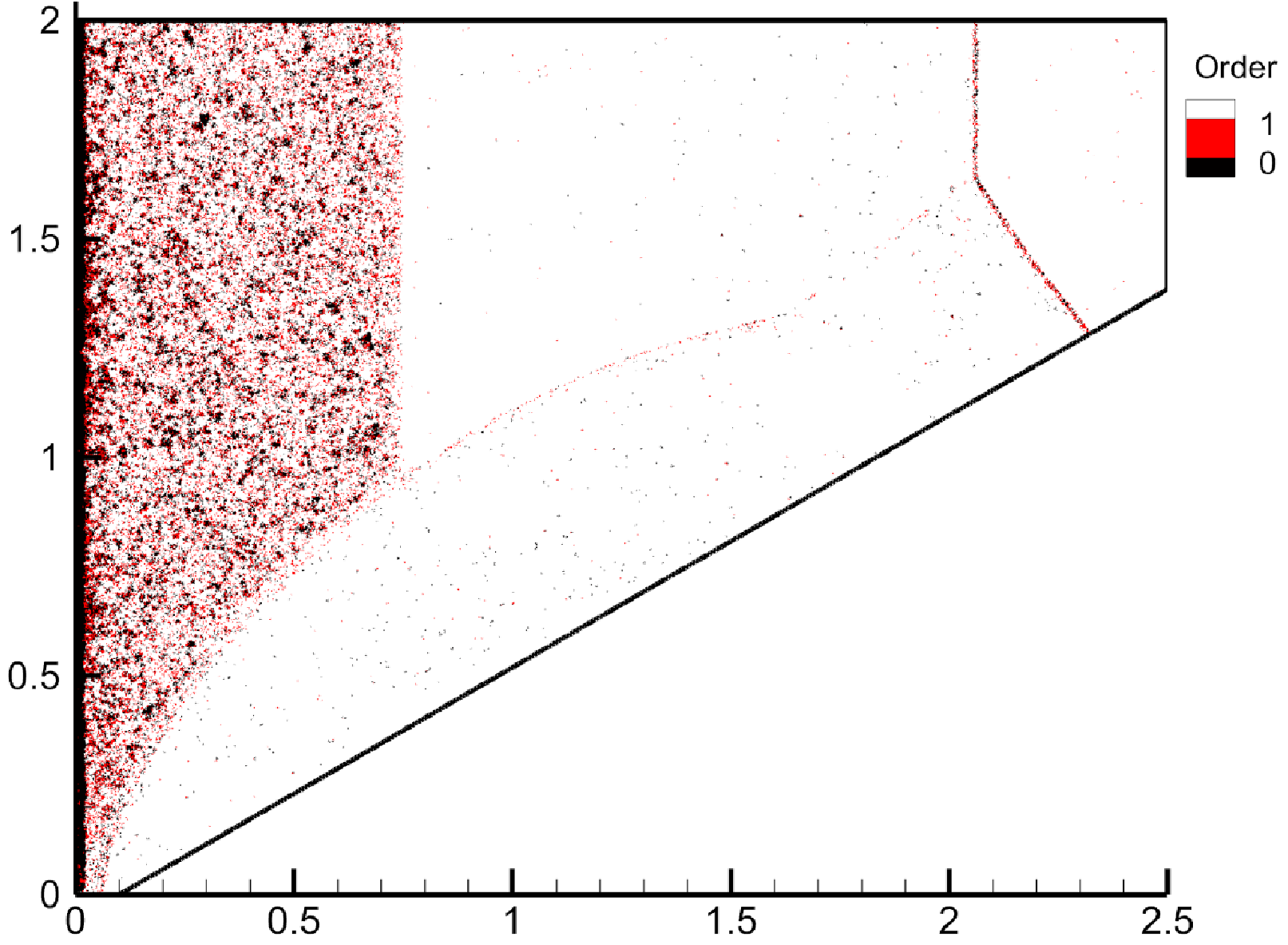}}
  \subfigure[$k=3$]{
  \includegraphics[height=6 cm]{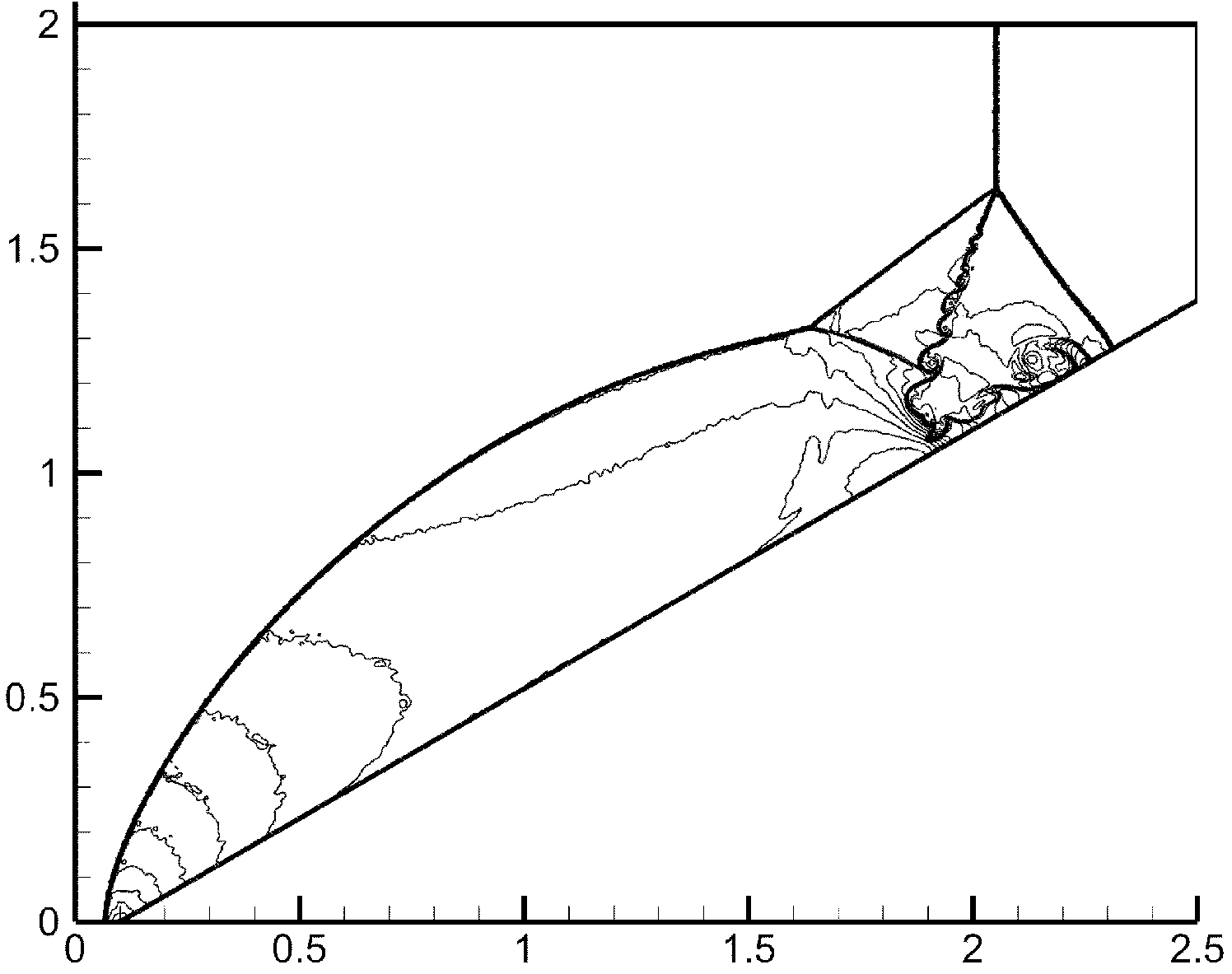}}
  \subfigure[$k=3$]{
  \includegraphics[height=6 cm]{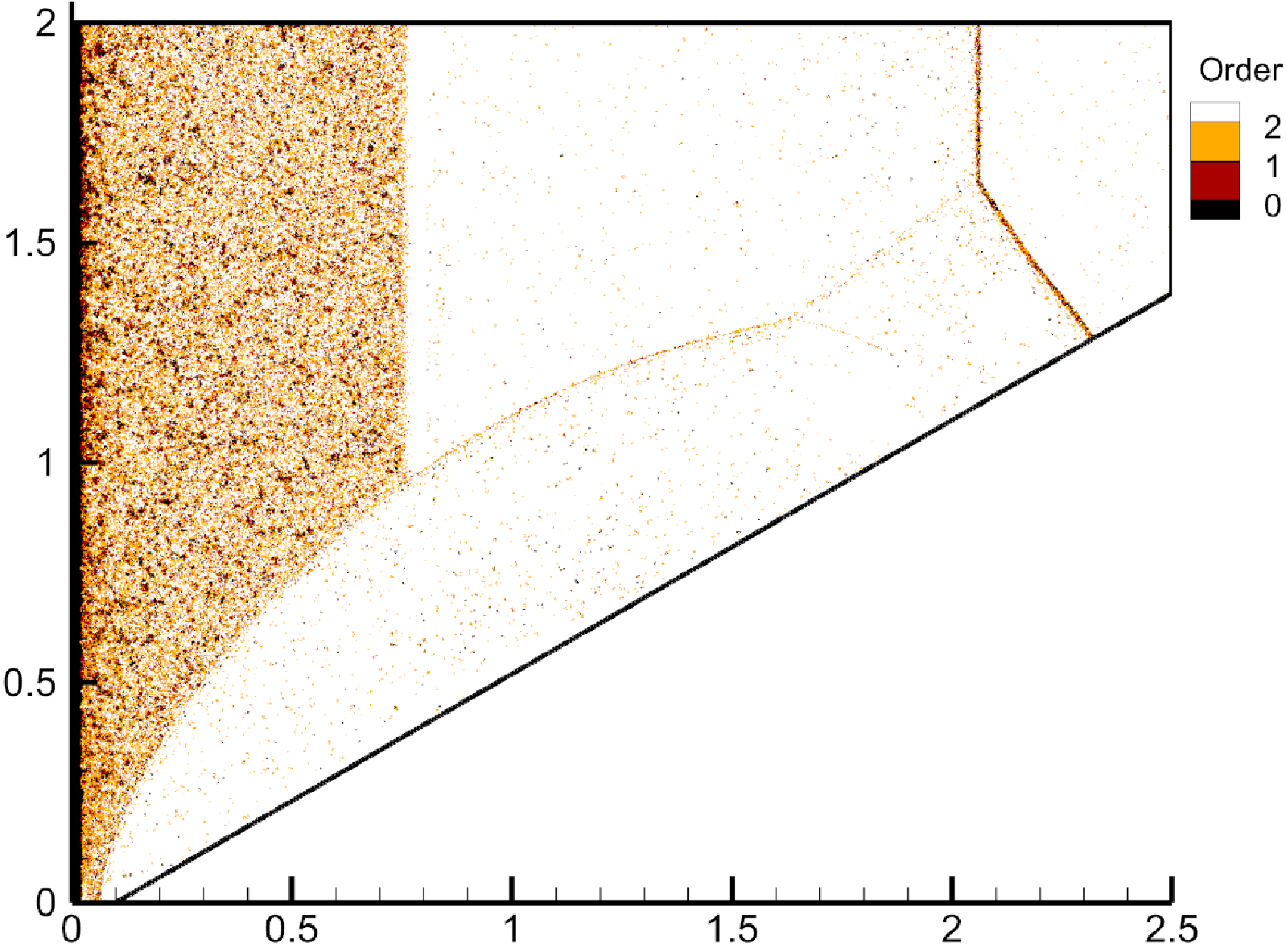}}
  \caption{Density contours (left column) and the polynomial order (right column) of the standard double Mach reflection problem at $t=0.2$ calculated 
  by RKDG schemes with the MR limiter using quasi-uniform triangular cells 
  with the edge length equal to $1/240$ on the boundaries
  The density contours contain 30 equidistant contours from 1.5 to 21.5.
  The white parts of the polynomial order plots represent the original $k$th-order DG polynomial.}
 \label{FIG:Standard_DMR_Tri}
 \end{figure}

Fig. \ref{FIG:Standard_DMR_Qua} shows the density contours and the polynomial order at $t=0.2$ calculated 
by RKDG schemes with the MR limiter using $600\times480$ structured quadrilateral cells.
Fig. \ref{FIG:Standard_DMR_Tri} shows the corresponding results calculated by RKDG schemes with the MR limiter using quasi-uniform triangular cells 
with the edge length equal to $1/240$ on the boundaries of the computational domain.
Except for the case of $k=1$ using the quadrilateral mesh, the RKDG schemes with the MR limiter 
can resolve the small vortices very well.
Meanwhile, the flow fields remain almost non-oscillatory.

 \paragraph{Example 4.3.2} Mach 10 shock passing a sharp corner.

\begin{figure}[htbp]
  \centering
  \includegraphics[width=8 cm]{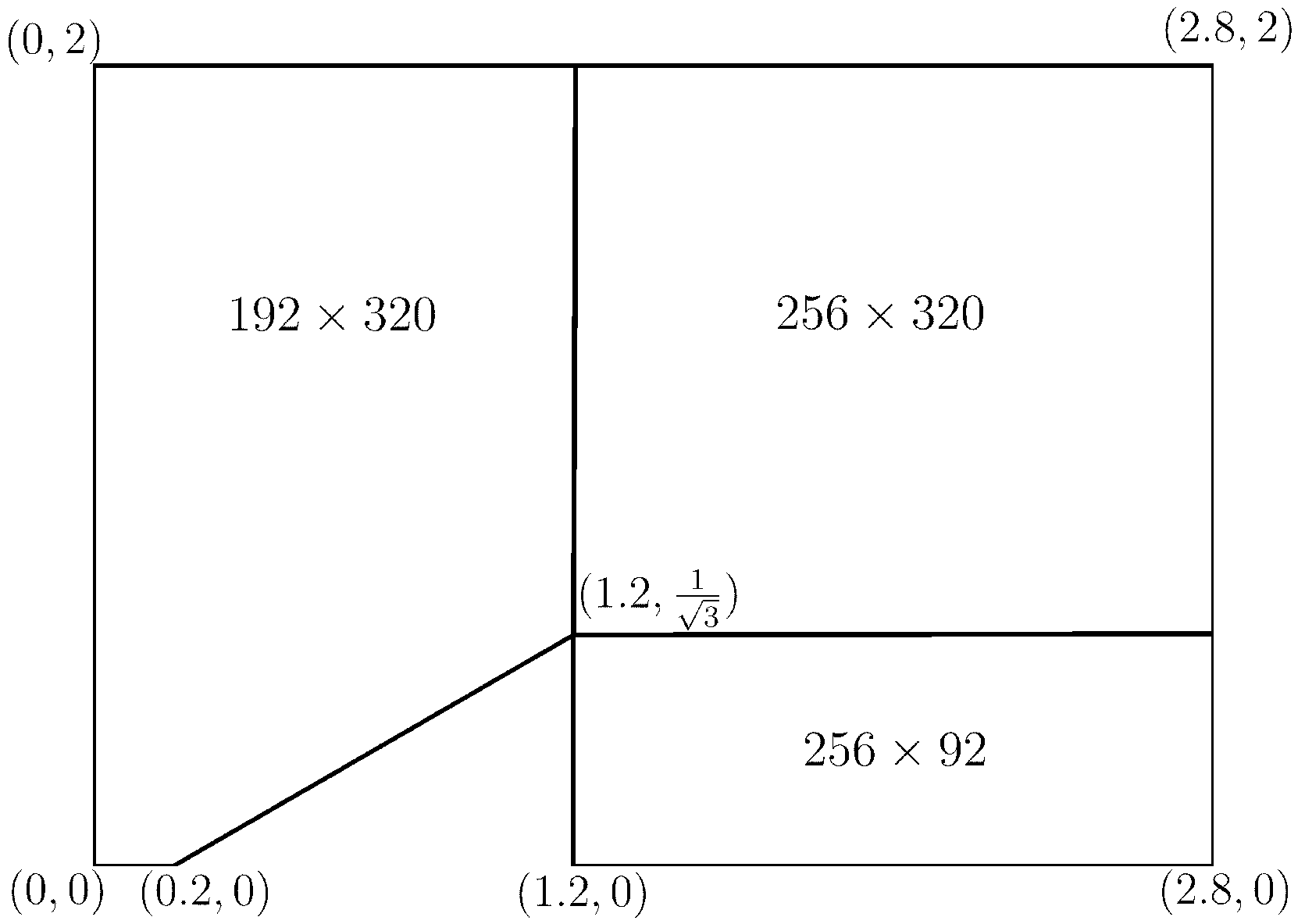}
  \caption{Mach 10 shock passing a sharp corner. 
  The computational domain and the distribution of multi-block structured quadrilateral cells.}
 \label{FIG:Shock_Diffraction_Domain}
 \end{figure}

 This case is taken from Zhang \cite{Zhang2017PositivityDGNS}.
 The computational domain is bounded by segments connecting the points 
 $(0,0)$, $(0.2,0)$, $(1.2,1/\sqrt{3})$, $(1.2,0)$, $(2.8,0)$, $(2.8,2)$, $(0,2)$ as shown in Fig. \ref{FIG:Shock_Diffraction_Domain}.
The initial condition is the same as the double Mach reflection problem.
The post-shocked states are enforced on the left boundary.
Nonreflective boundary conditions are applied to the right boundary.
Reflective wall boundary conditions are imposed on the remaining parts.

\begin{figure}[htbp]
  \centering
  \subfigure[$k=1$]{
  \includegraphics[height=5.5 cm]{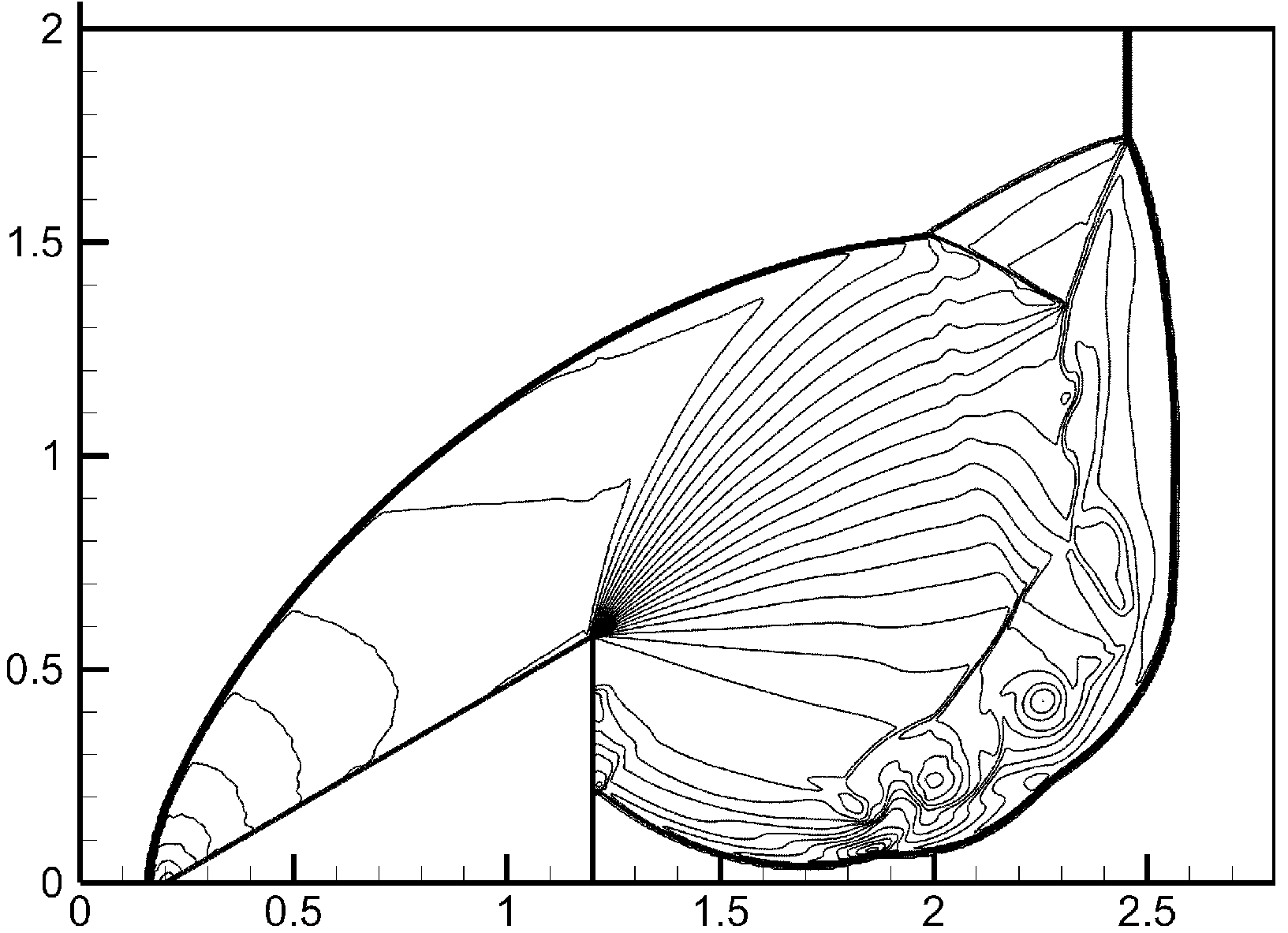}}
  \subfigure[$k=1$]{
  \includegraphics[height=5.5 cm]{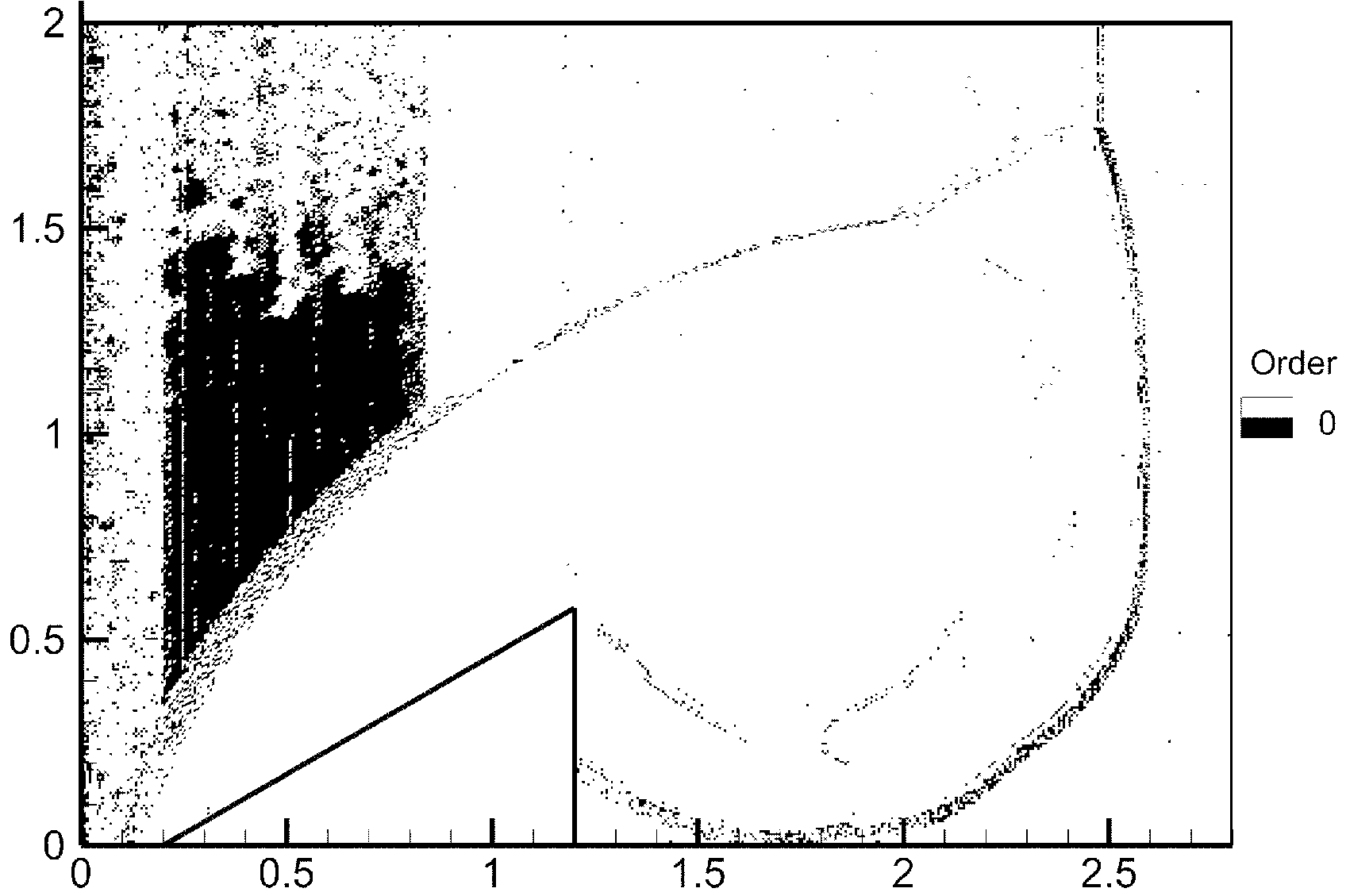}}
  \subfigure[$k=2$]{
  \includegraphics[height=5.5 cm]{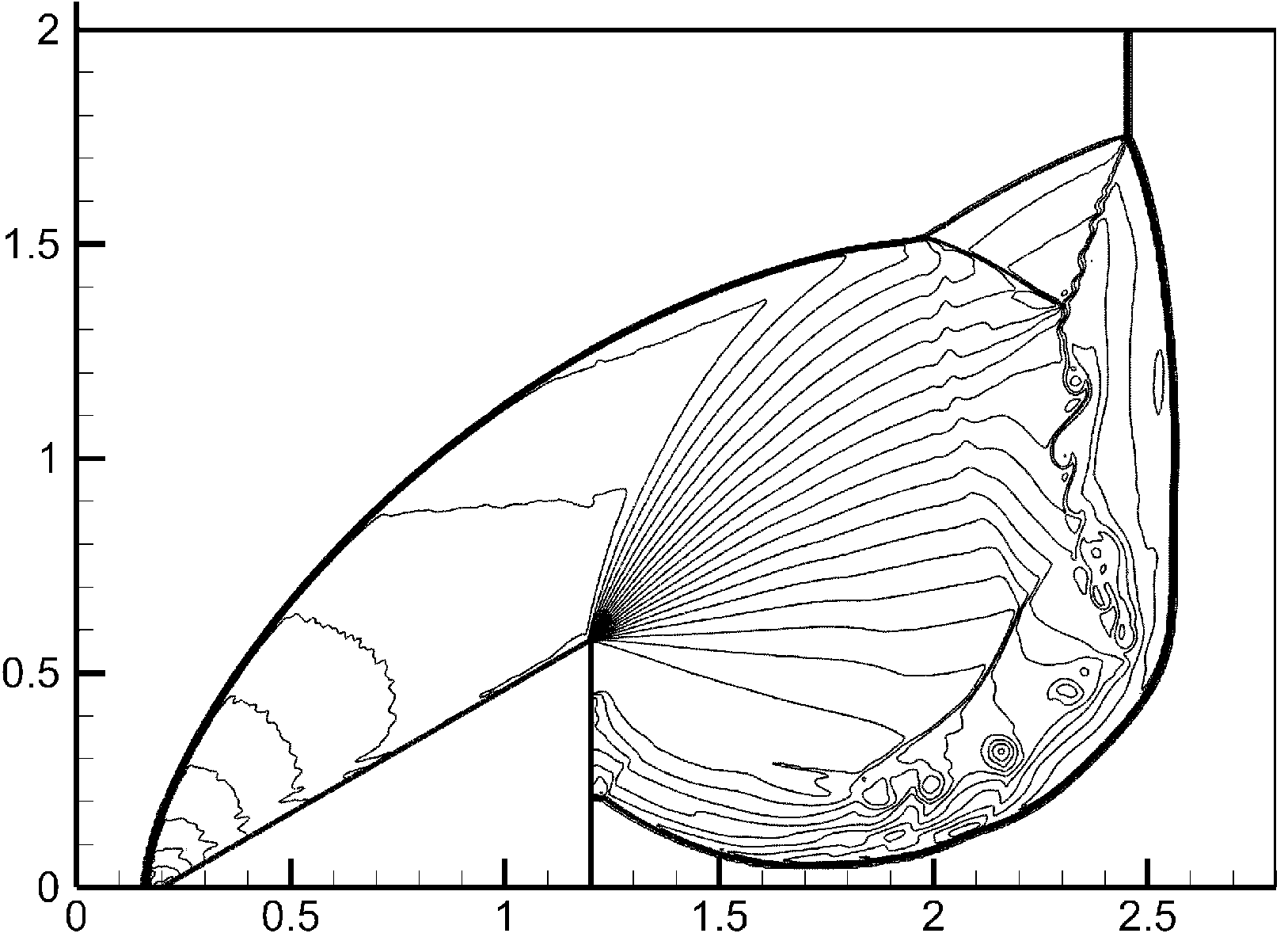}}
  \subfigure[$k=2$]{
  \includegraphics[height=5.5 cm]{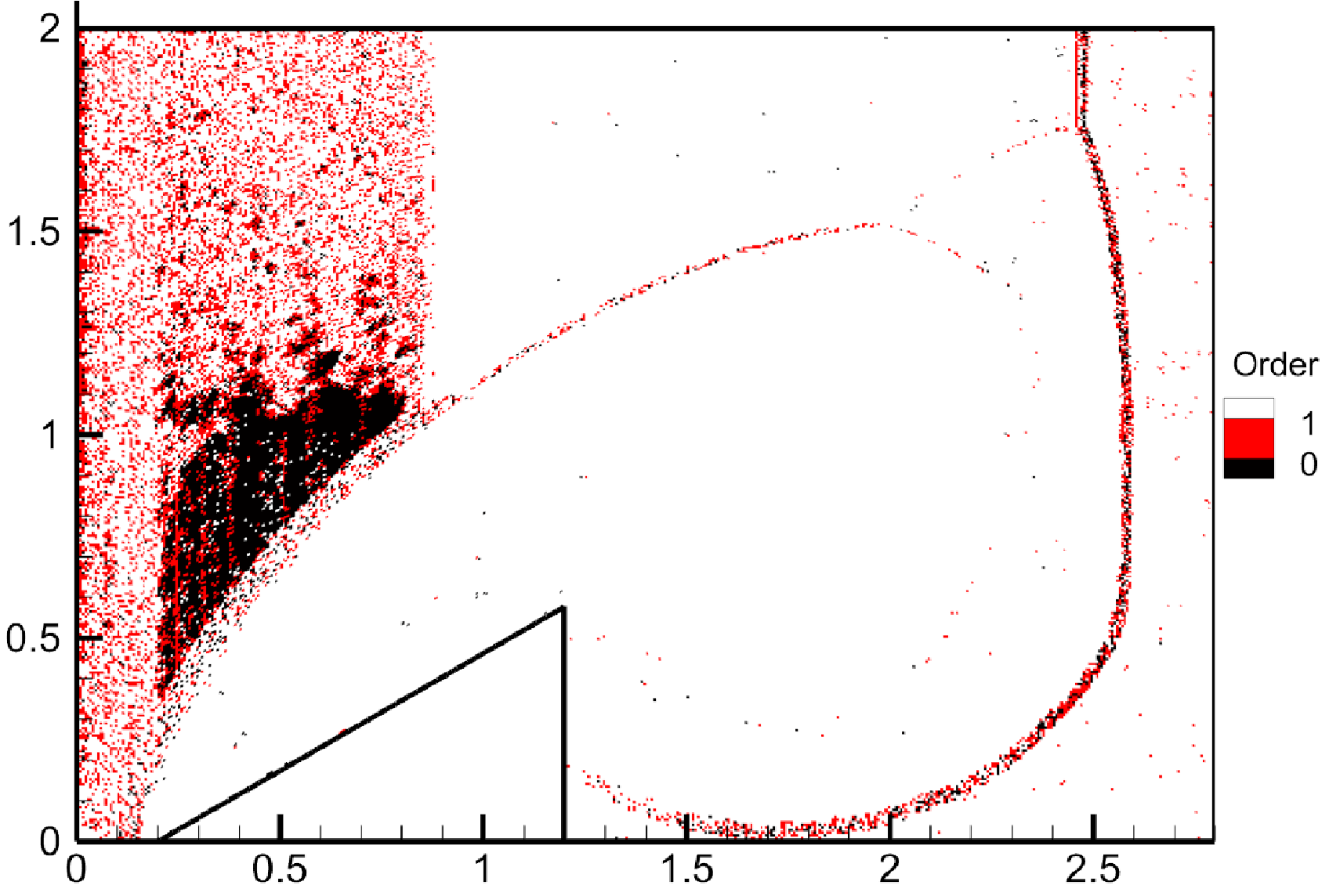}}
  \subfigure[$k=3$]{
  \includegraphics[height=5.5 cm]{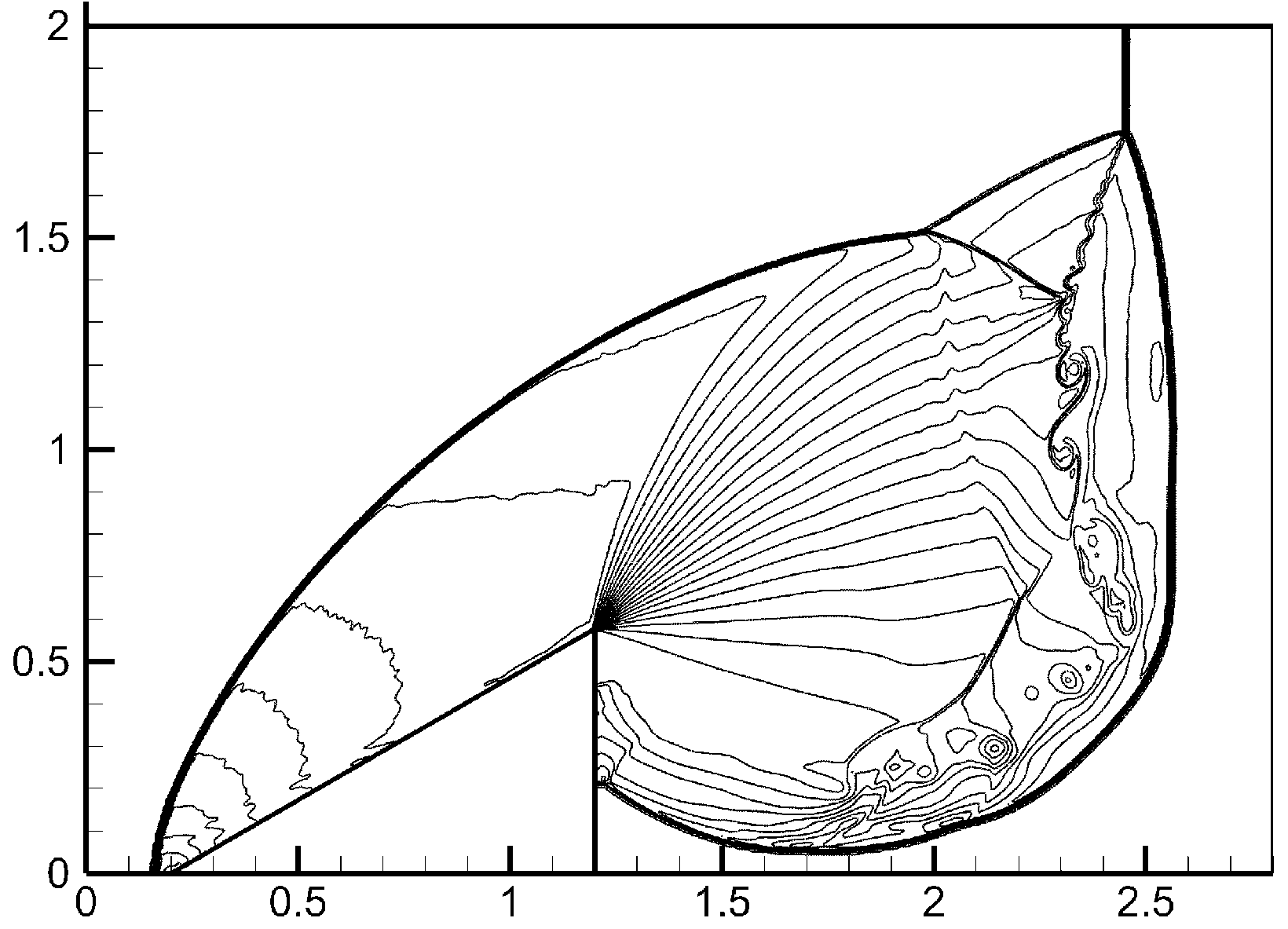}}
  \subfigure[$k=3$]{
  \includegraphics[height=5.5 cm]{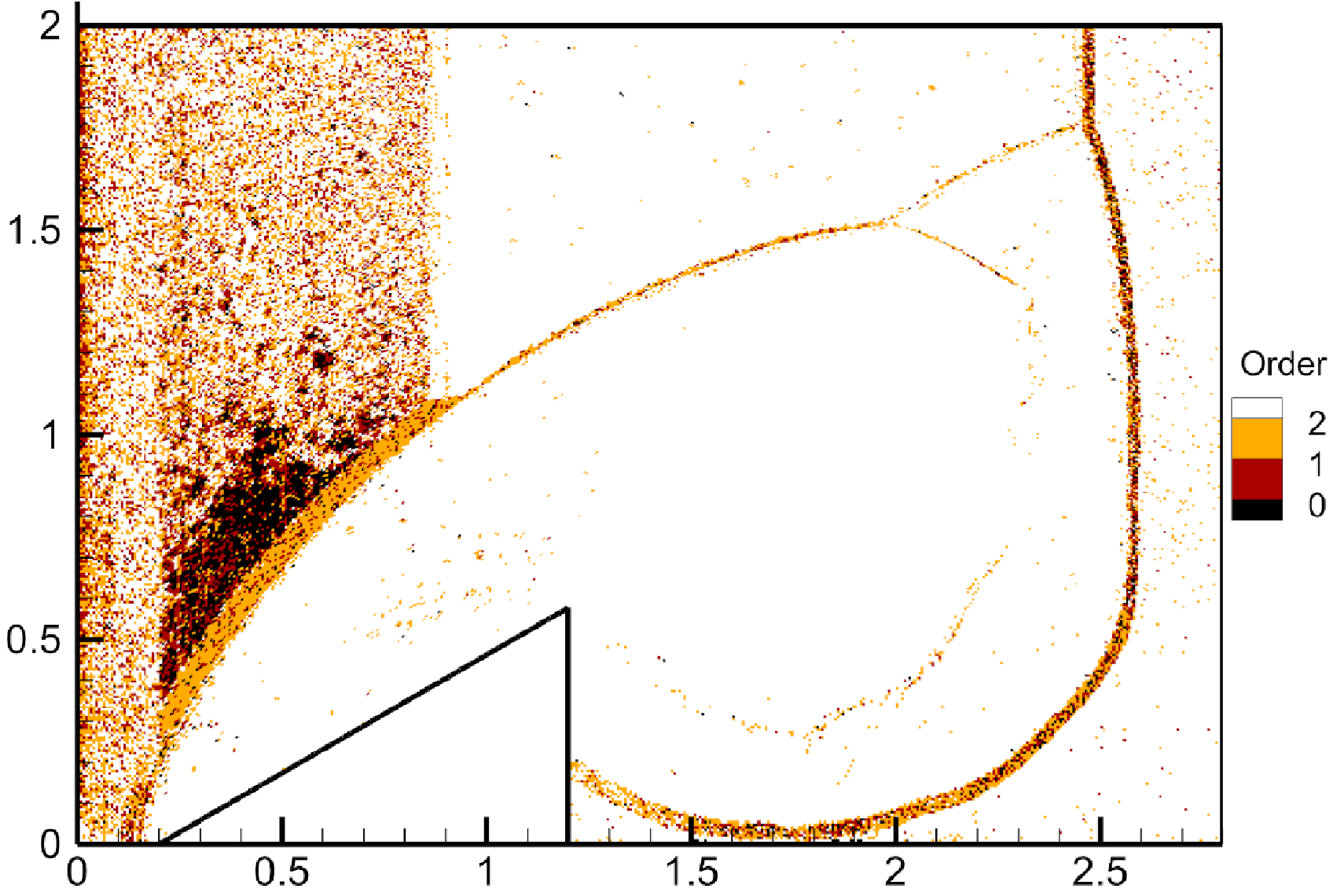}}
  \caption{Density contours (left column) and the polynomial order (right column) of Mach 10 shock passing a sharp corner at $t=0.24$ calculated 
  by RKDG schemes with the MR limiter using multi-block structured quadrilateral cells.
  The density contours contain 30 equidistant contours from 0.05 to 21.5.
  The white parts of the polynomial order plots represent the original $k$th-order DG polynomial.}
 \label{FIG:Shock_Diffraction_Qua}
 \end{figure}

\begin{figure}[htbp]
  \centering
  \subfigure[$k=1$]{
  \includegraphics[height=5.5 cm]{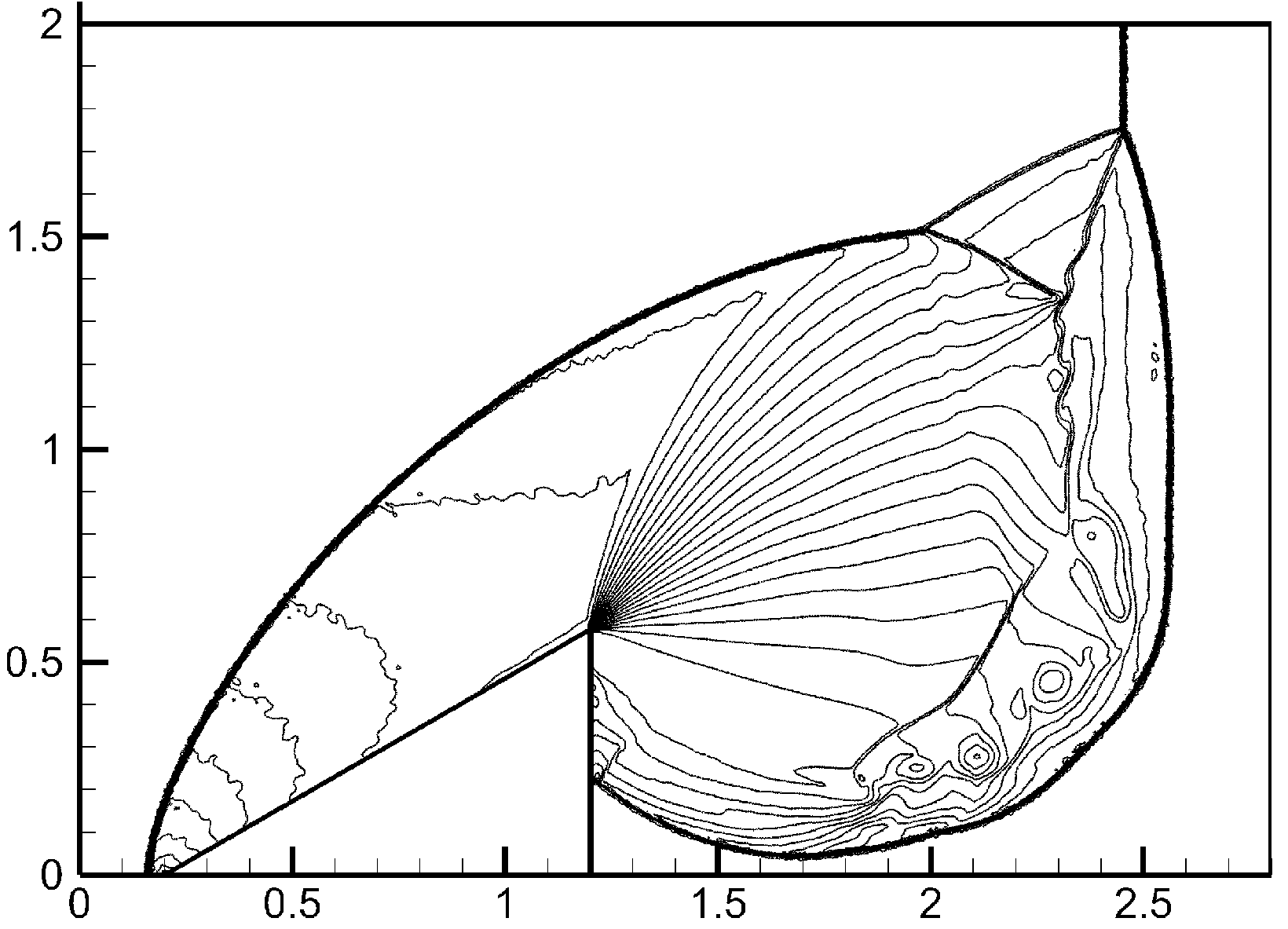}}
  \subfigure[$k=1$]{
  \includegraphics[height=5.5 cm]{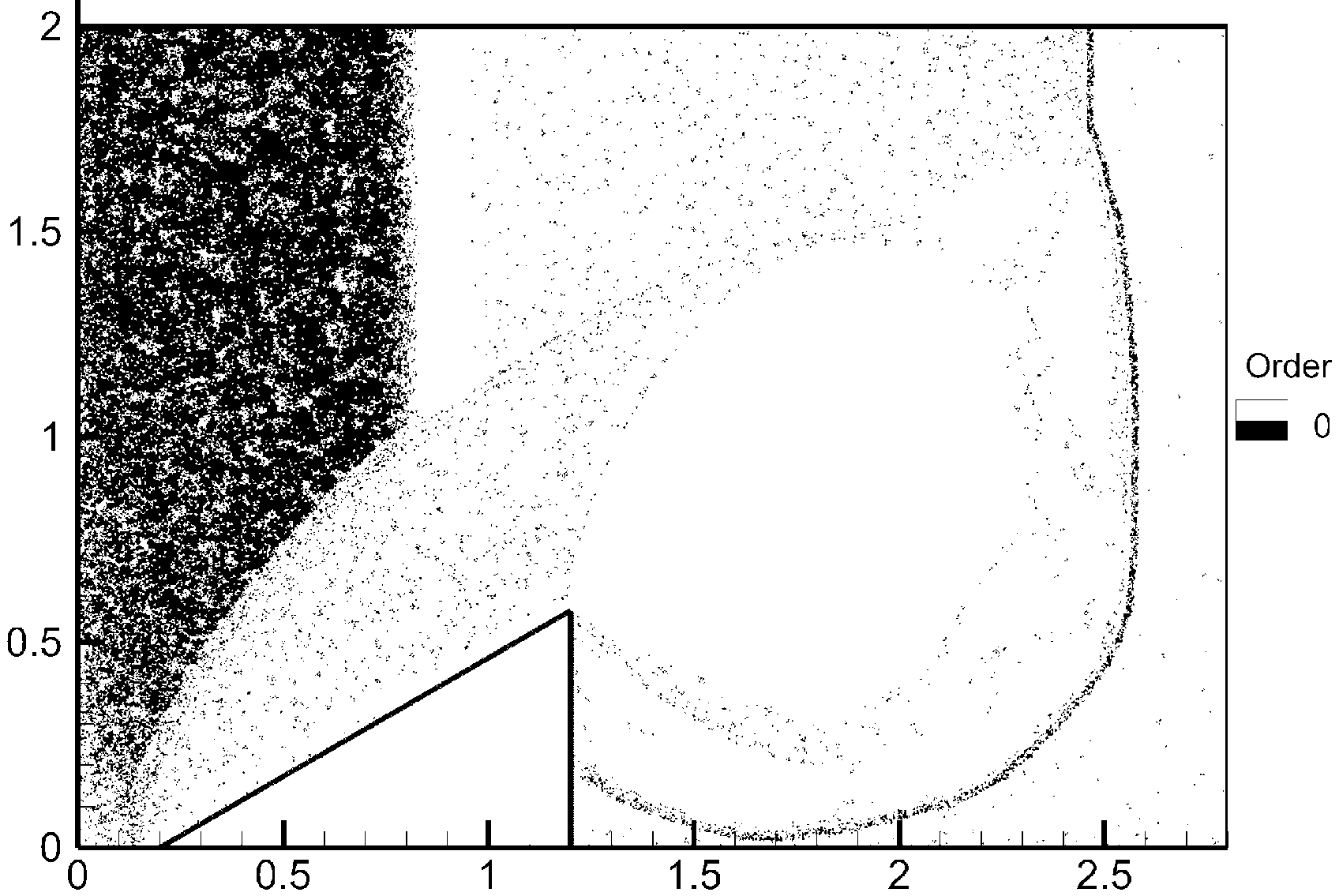}}
  \subfigure[$k=2$]{
  \includegraphics[height=5.5 cm]{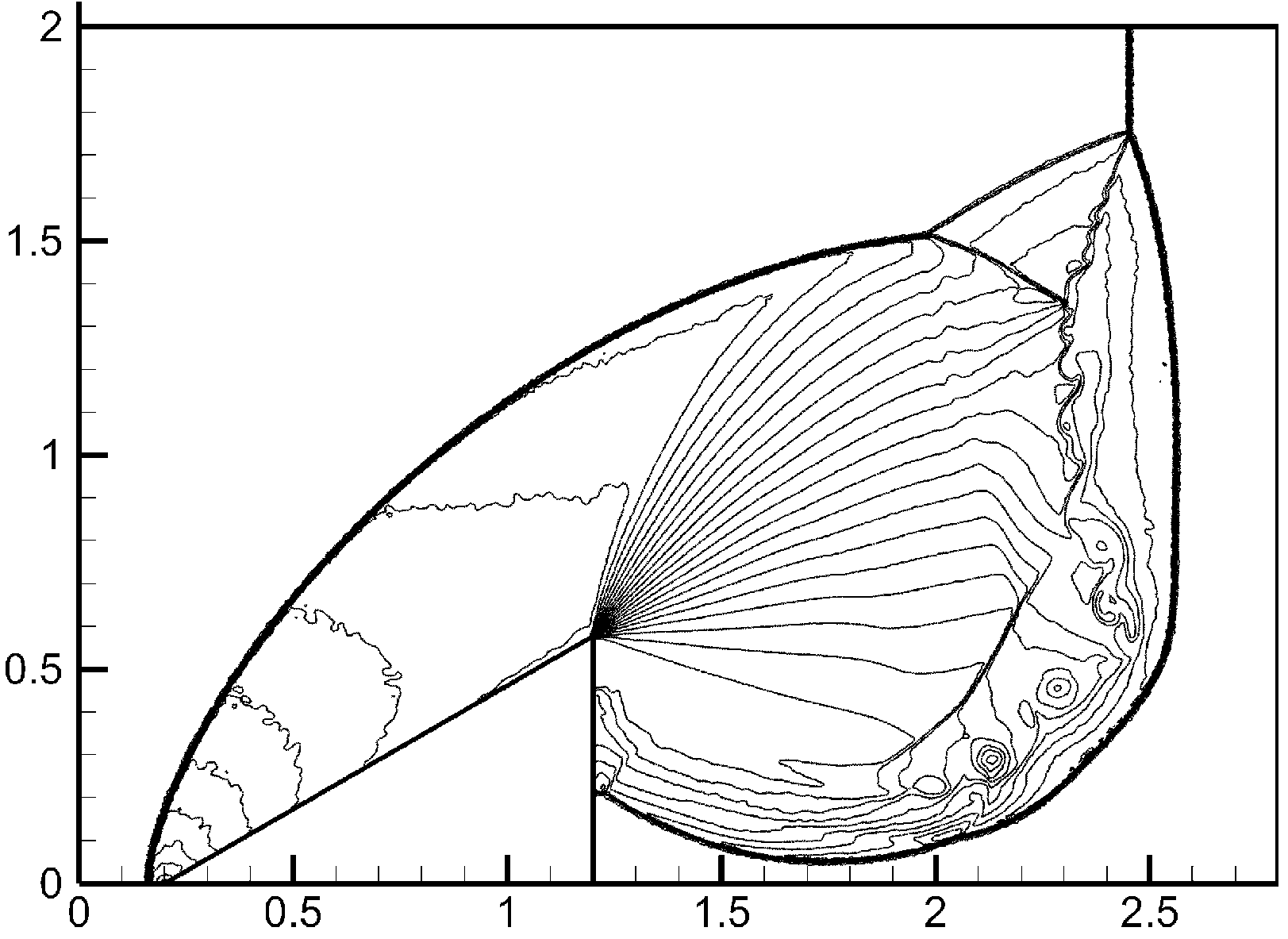}}
  \subfigure[$k=2$]{
  \includegraphics[height=5.5 cm]{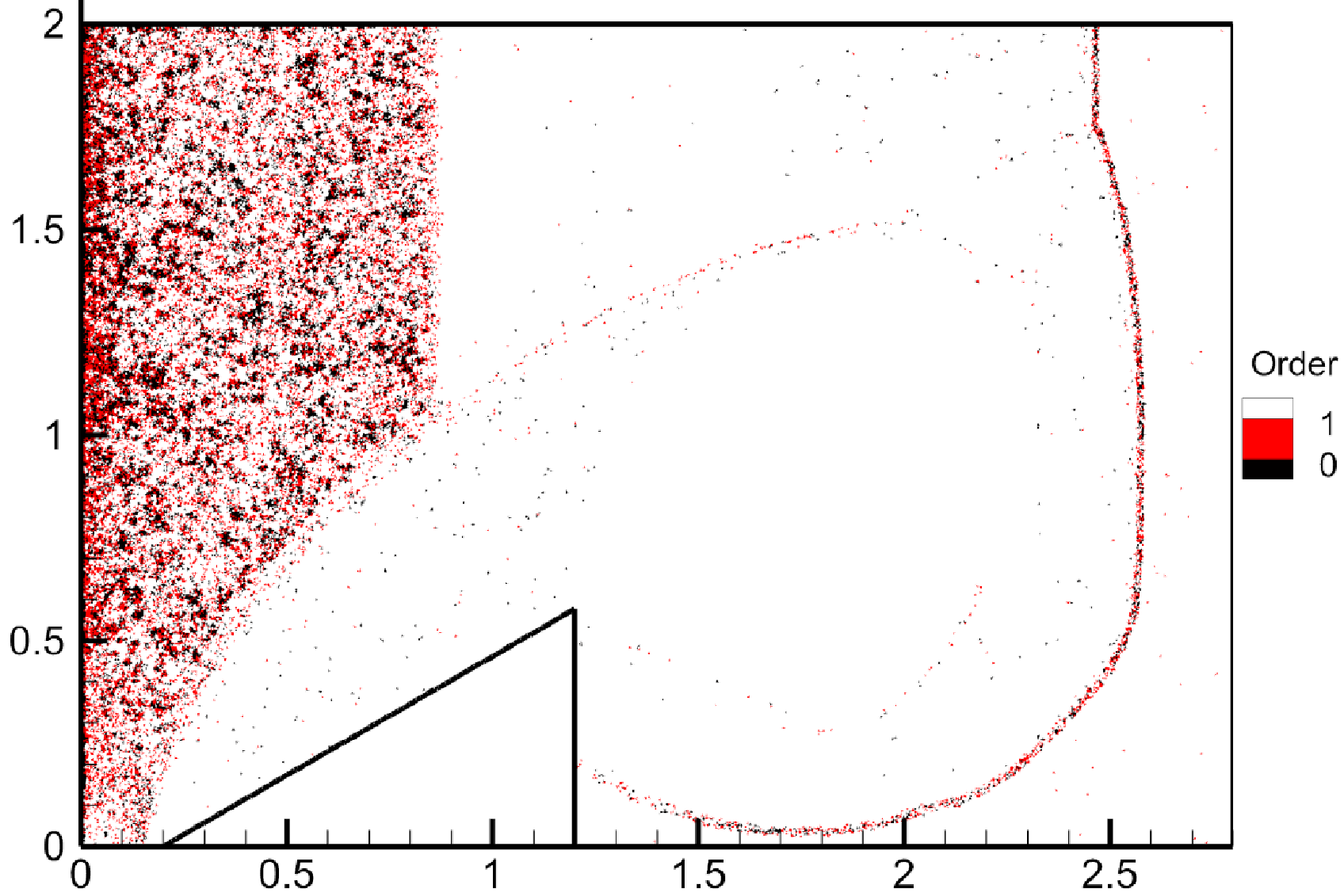}}
  \subfigure[$k=3$]{
  \includegraphics[height=5.5 cm]{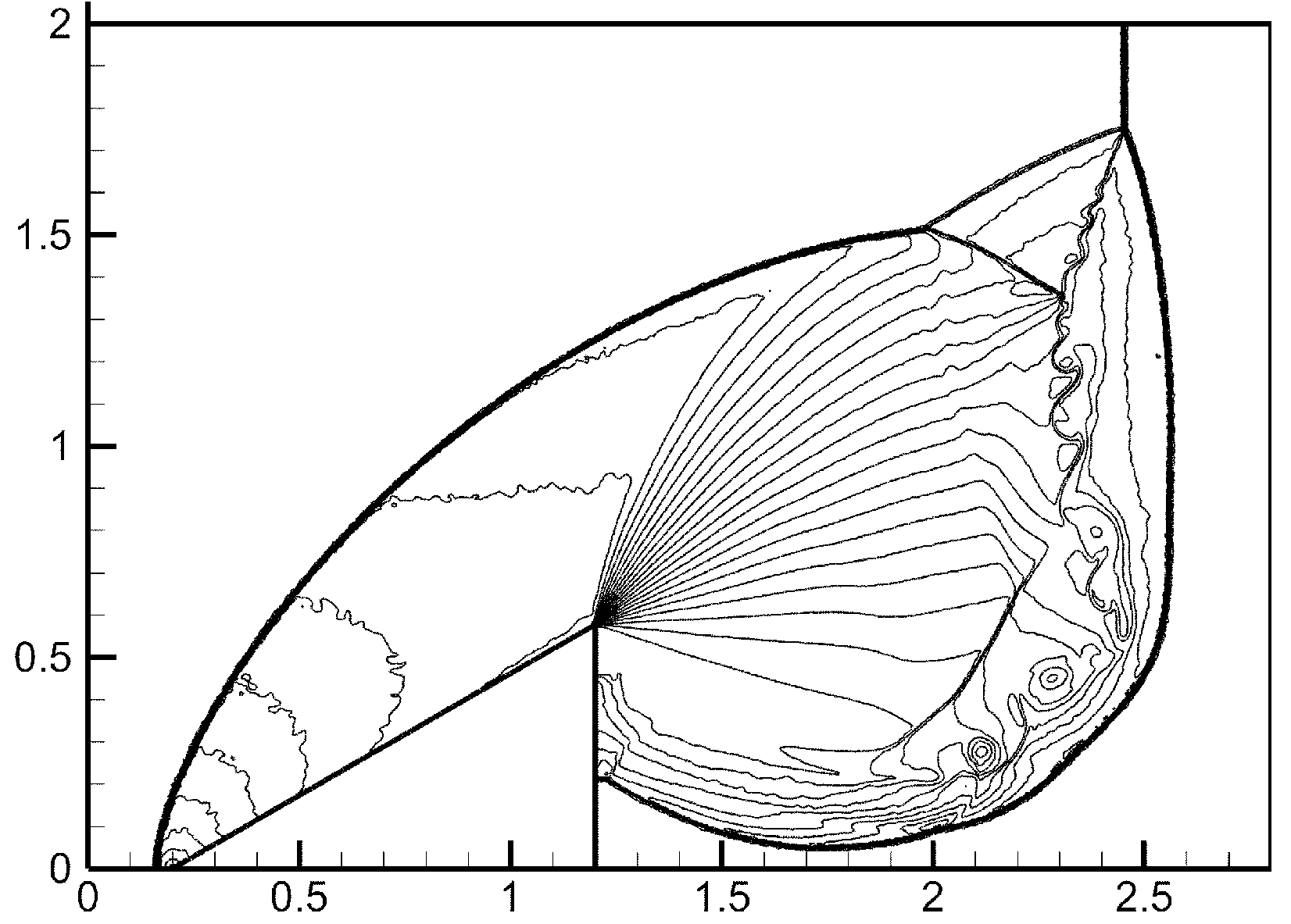}}
  \subfigure[$k=3$]{
  \includegraphics[height=5.5 cm]{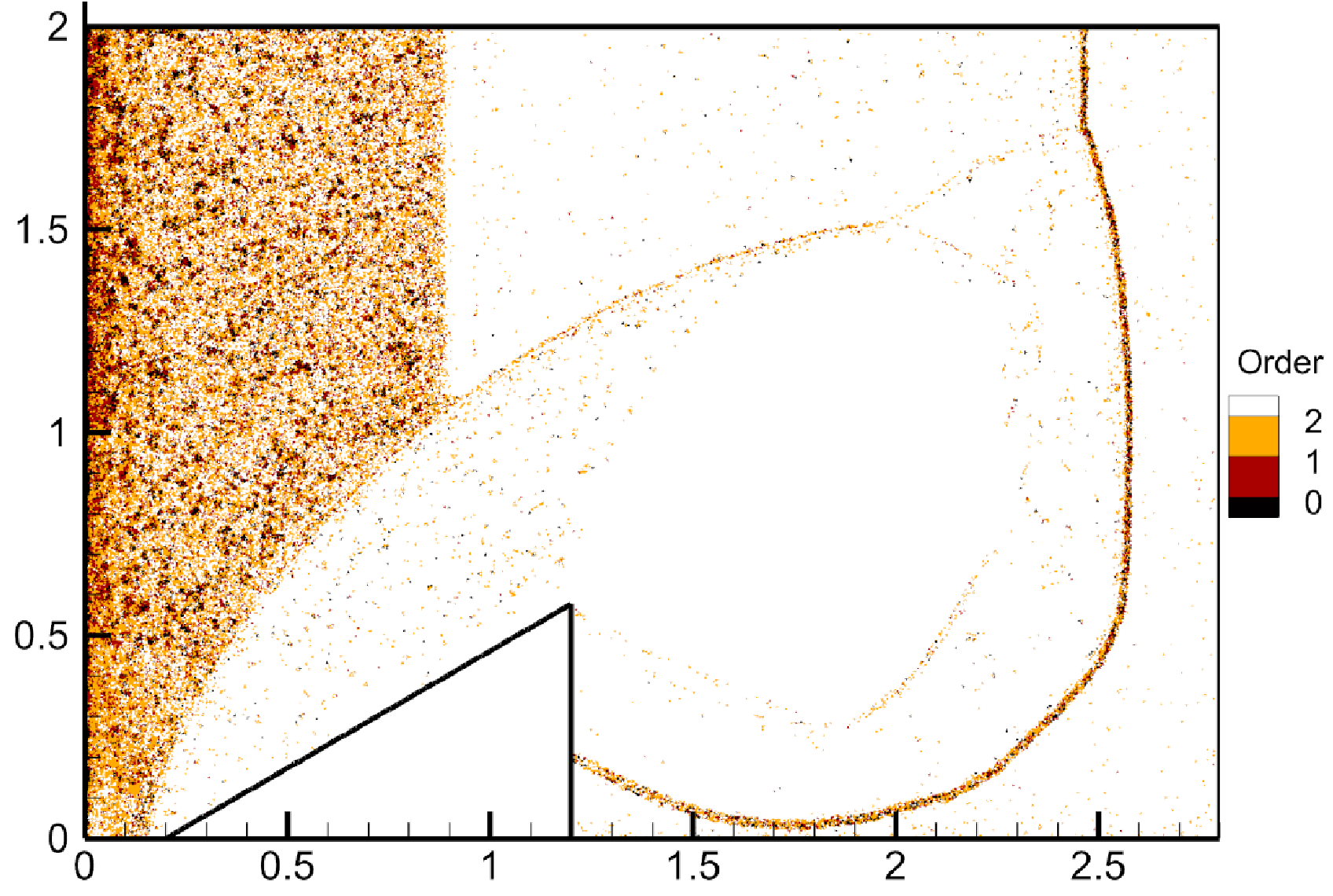}}
  \caption{Density contours (left column) and the polynomial order (right column) of Mach 10 shock passing a sharp corner at $t=0.24$ calculated 
  by RKDG schemes with the MR limiter using triangular cells.
  The density contours contain 30 equidistant contours from 0.05 to 21.5.
  The white parts of the polynomial order plots represent the original $k$th-order DG polynomial.}
 \label{FIG:Shock_Diffraction_Tri}
 \end{figure}

As shown in Fig. \ref{FIG:Shock_Diffraction_Domain}, we divide the computational domain into three blocks 
and then discretize each block with structured quadrilateral cells.
Fig. \ref{FIG:Shock_Diffraction_Qua} shows the density contours and the polynomial order at $t=0.24$
computed by the RKDG schemes with the MR limiter using the quadrilateral meshes.
In addition, we also compute this problem using triangular cells 
with the edge length equal to $1/160$ on the boundaries of the computational domain.
Fig. \ref{FIG:Shock_Diffraction_Tri} shows the corresponding results.
The numerical results demonstrate that the RKDG schemes with the MR limiter
can effectively detect the discontinuities appearing in the flow field
and can capture shock waves, strong expansion waves, and flow instabilities without numerical oscillations.
The RKDG schemes with only positivity-preserving limiter 
generate a lot of numerical oscillations (Fig. 15 in \cite{Zhang2017PositivityDGNS}).

 \paragraph{Example 4.3.3} Forward-facing step problem.

 The last numerical example is another classical test case proposed by Woodward and Colella \cite{Woodward1984JCP}.
 This problem depicts a Mach 3 supersonic flow passing through a forward-facing step in a wind tunnel.
 The dimensionless height and width of the wind tunnel are 1 and 3, respectively.
  The step with a height of 0.2 is placed at a distance of 0.6 from the entrance of the wind tunnel. 
  The initial conditions are set as $(\rho,u,v,p)=(1.4,3,0,1)$, which is equivalent to a Mach 3 uniform flow in the tunnel. 
  Supersonic inflow and outflow boundary conditions are respectively applied to the left and right boundaries, 
and reflective wall boundary conditions are applied to the other boundaries. 
To avoid an erroneous entropy layer at the downstream bottom wall caused by the singularity at the corner of the step,
we follow the modification of Woodward and Colella \cite{Woodward1984JCP}.

\begin{figure}[htbp]
  \centering
  \subfigure[$k=1$]{
  \includegraphics[height=3 cm]{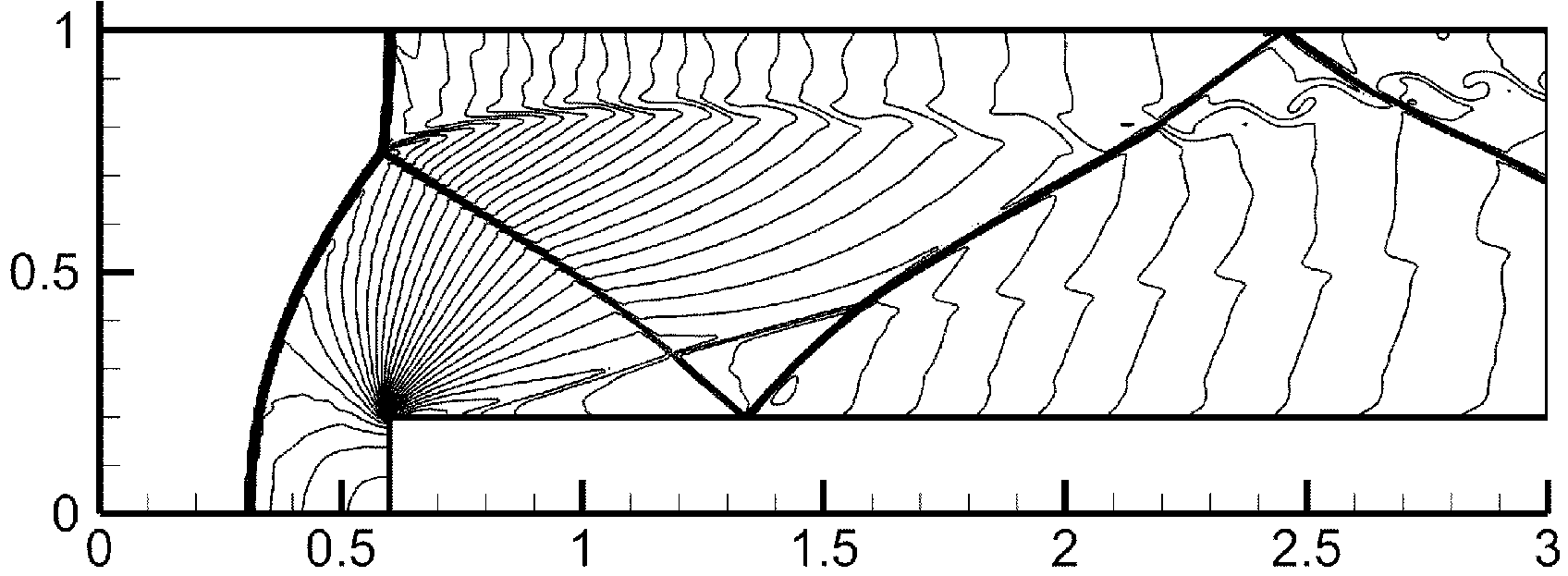}}
  \subfigure[$k=1$]{
  \includegraphics[height=3 cm]{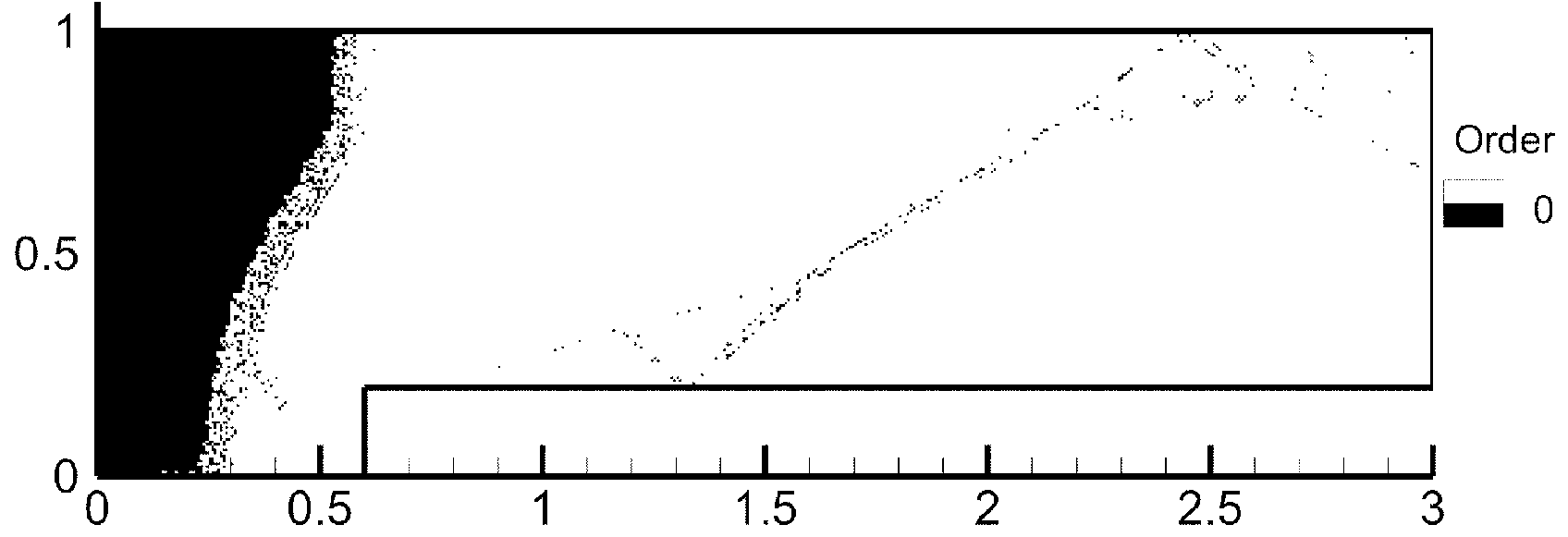}}
  \subfigure[$k=2$]{
  \includegraphics[height=3 cm]{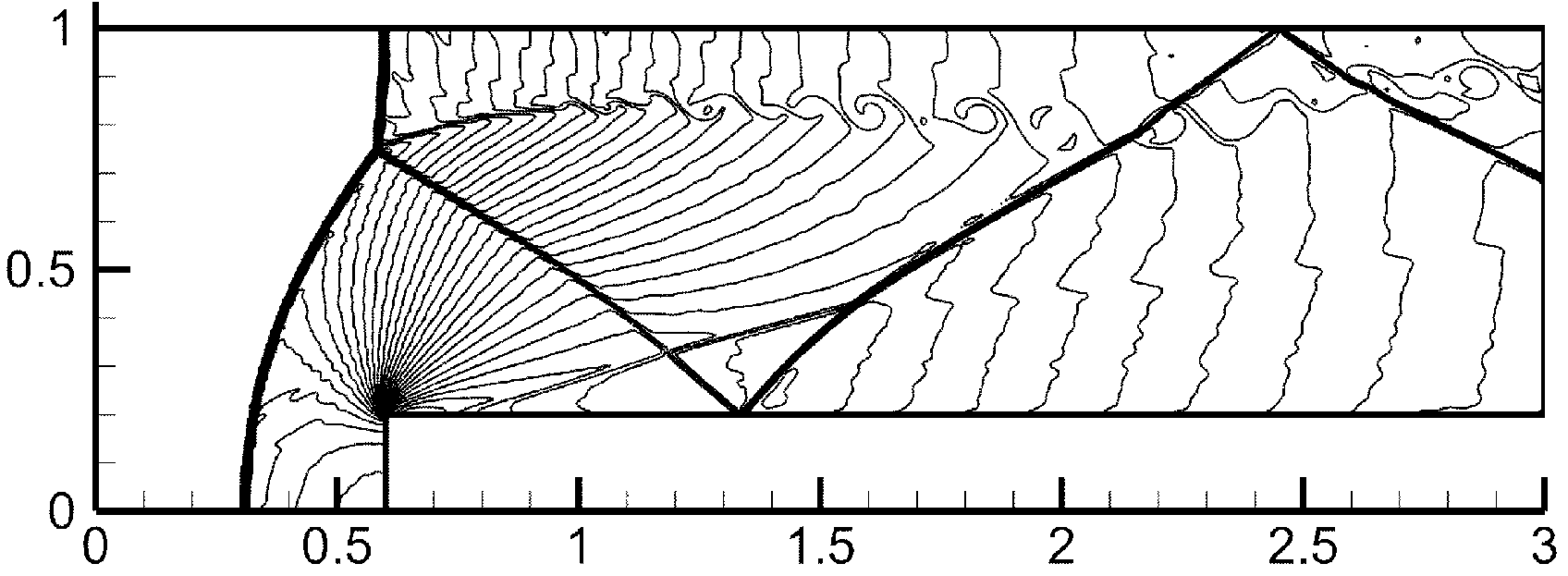}}
  \subfigure[$k=2$]{
  \includegraphics[height=3 cm]{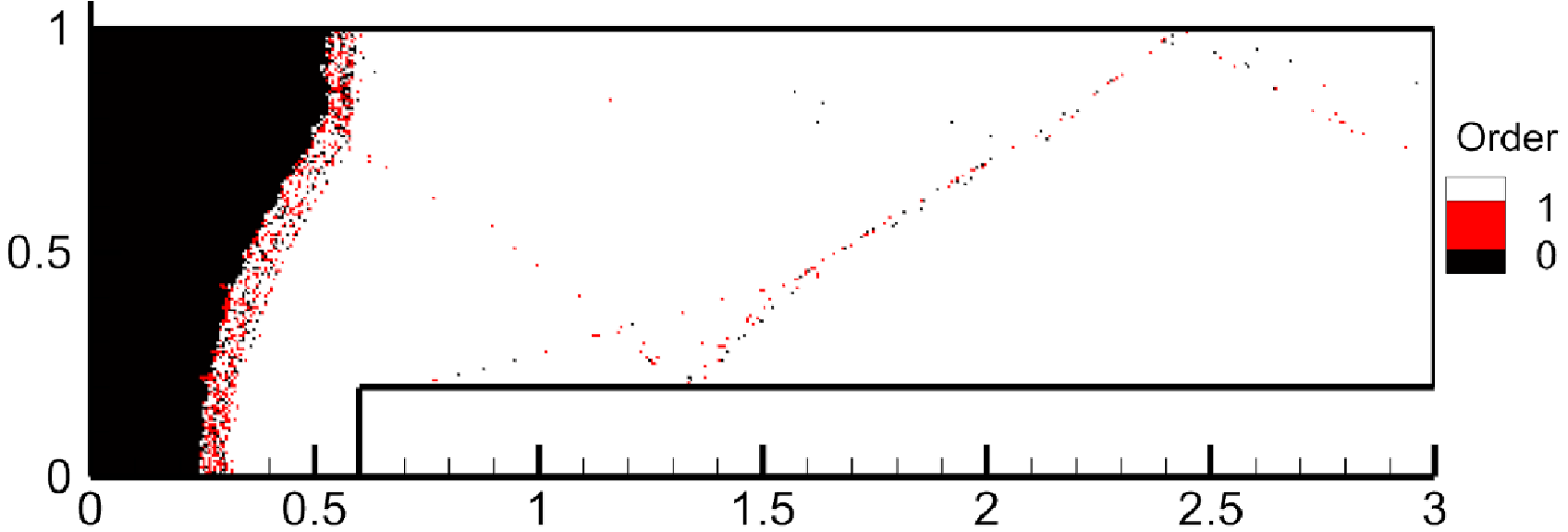}}
  \subfigure[$k=3$]{
  \includegraphics[height=3 cm]{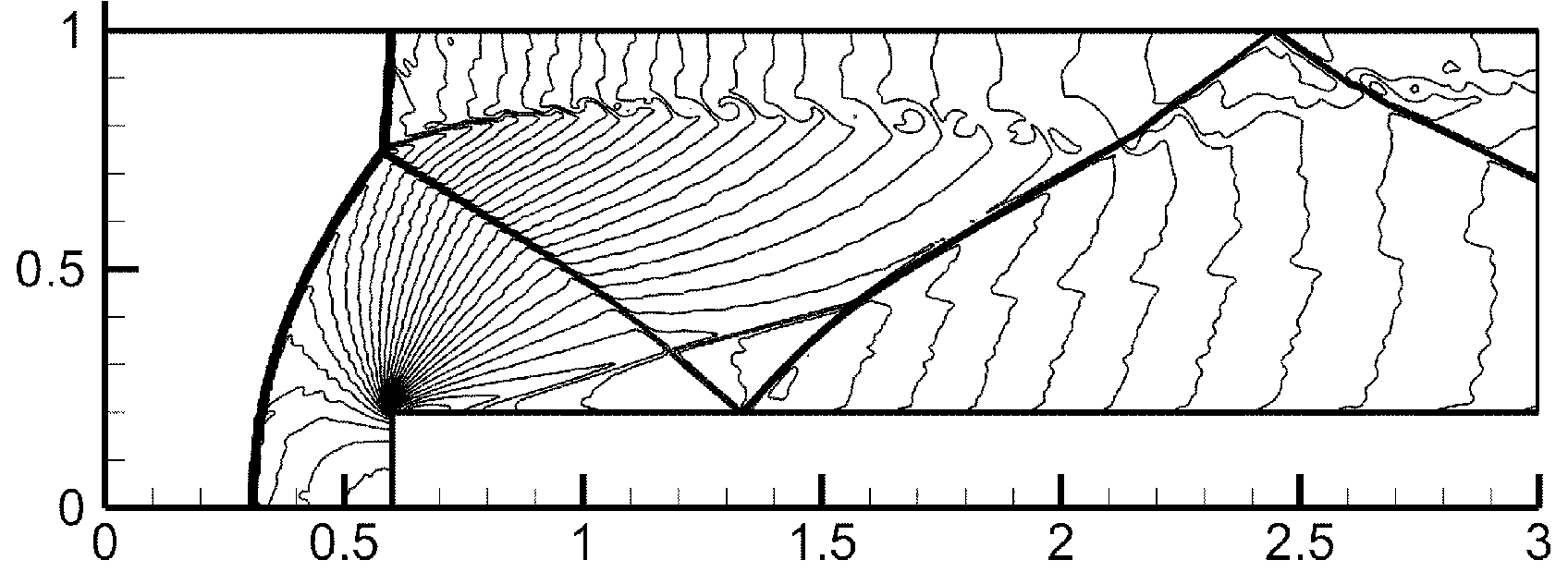}}
  \subfigure[$k=3$]{
  \includegraphics[height=3 cm]{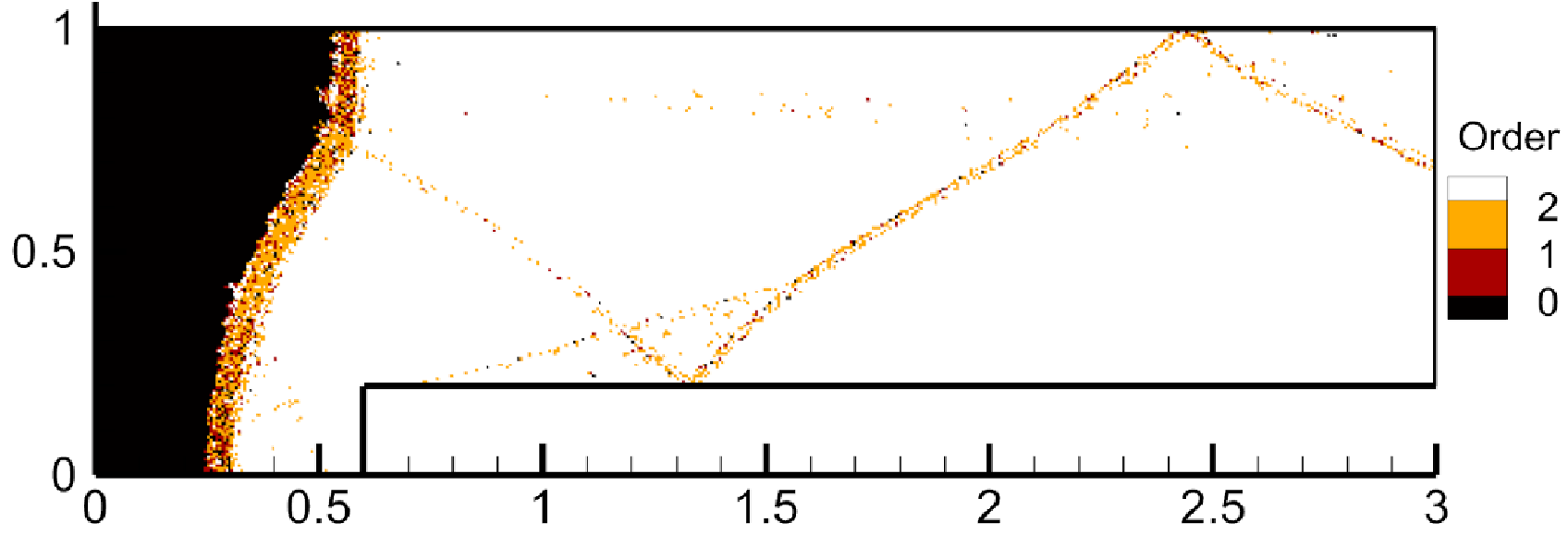}}
  \caption{Density contours (left column) and the polynomial order (right column) of the forward-facing step problem at $t=4$ calculated 
  by RKDG schemes with the MR limiter using uniform quadrilateral cells.
  The density contours contain 30 equidistant contours from 0.32 to 6.15.
  The white parts of the polynomial order plots represent the original $k$th-order DG polynomial.}
 \label{FIG:FST_Qua}
 \end{figure}

\begin{figure}[htbp]
  \centering
  \subfigure[$k=1$]{
  \includegraphics[height=3 cm]{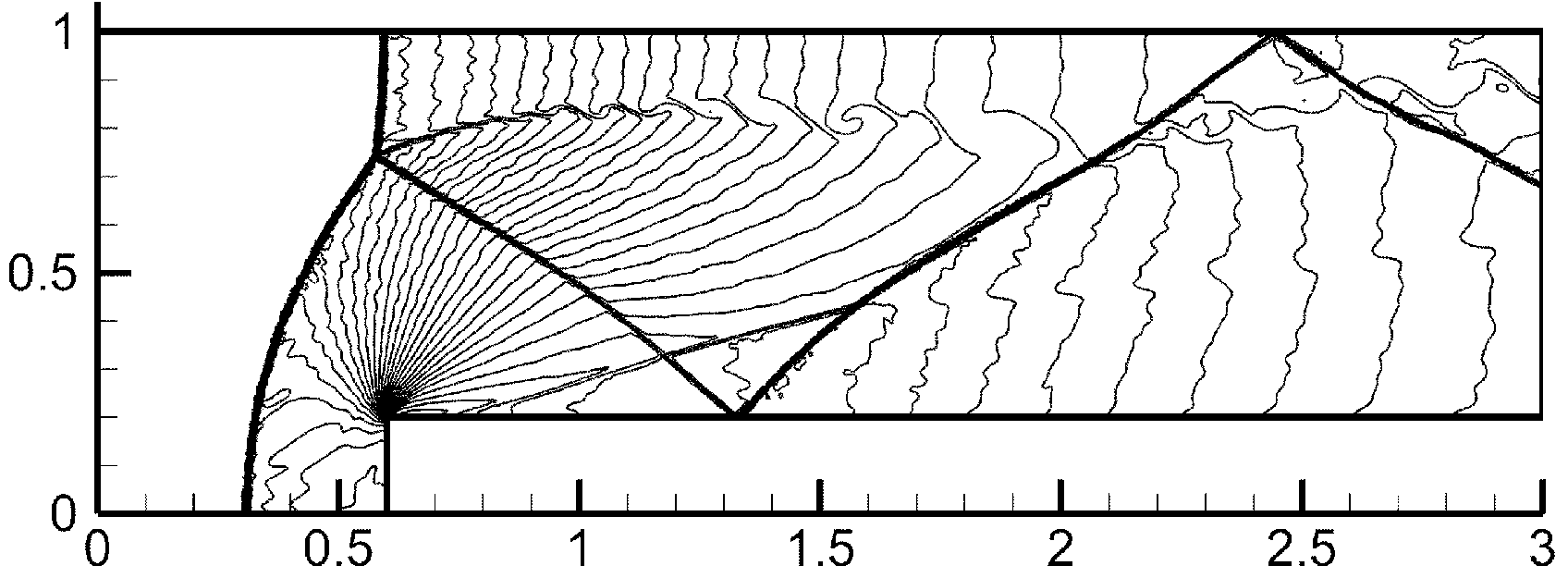}}
  \subfigure[$k=1$]{
  \includegraphics[height=3 cm]{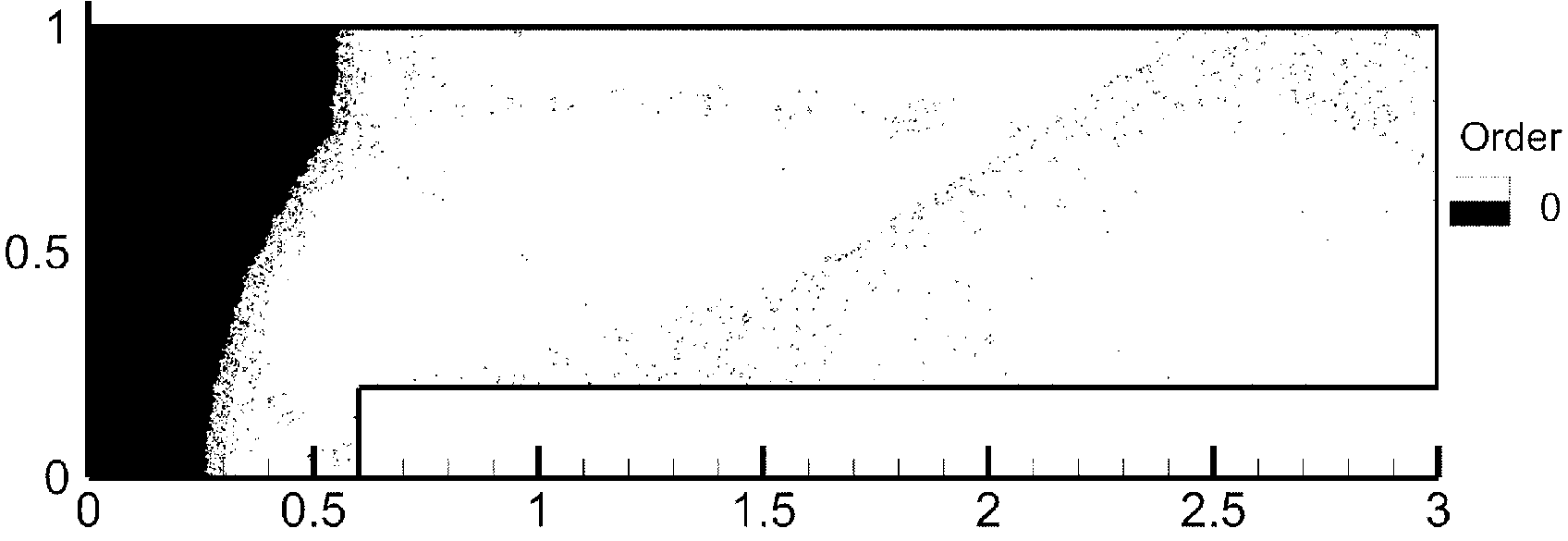}}
  \subfigure[$k=2$]{
  \includegraphics[height=3 cm]{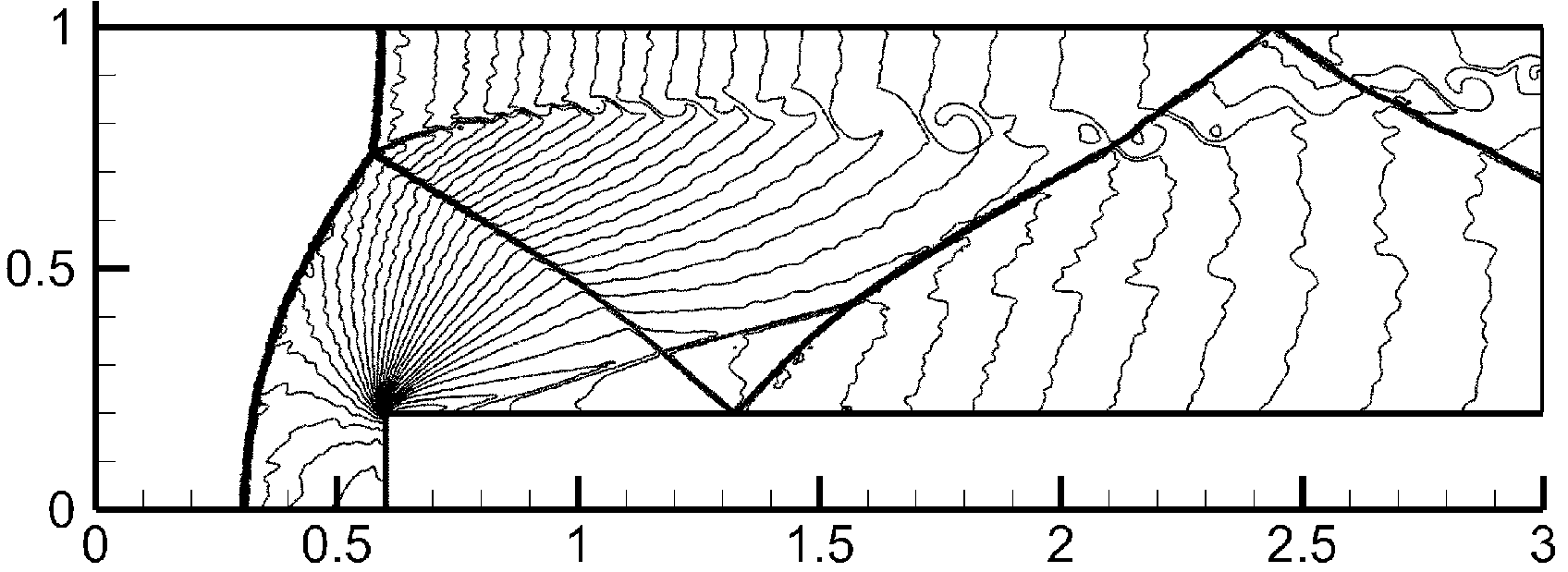}}
  \subfigure[$k=2$]{
  \includegraphics[height=3 cm]{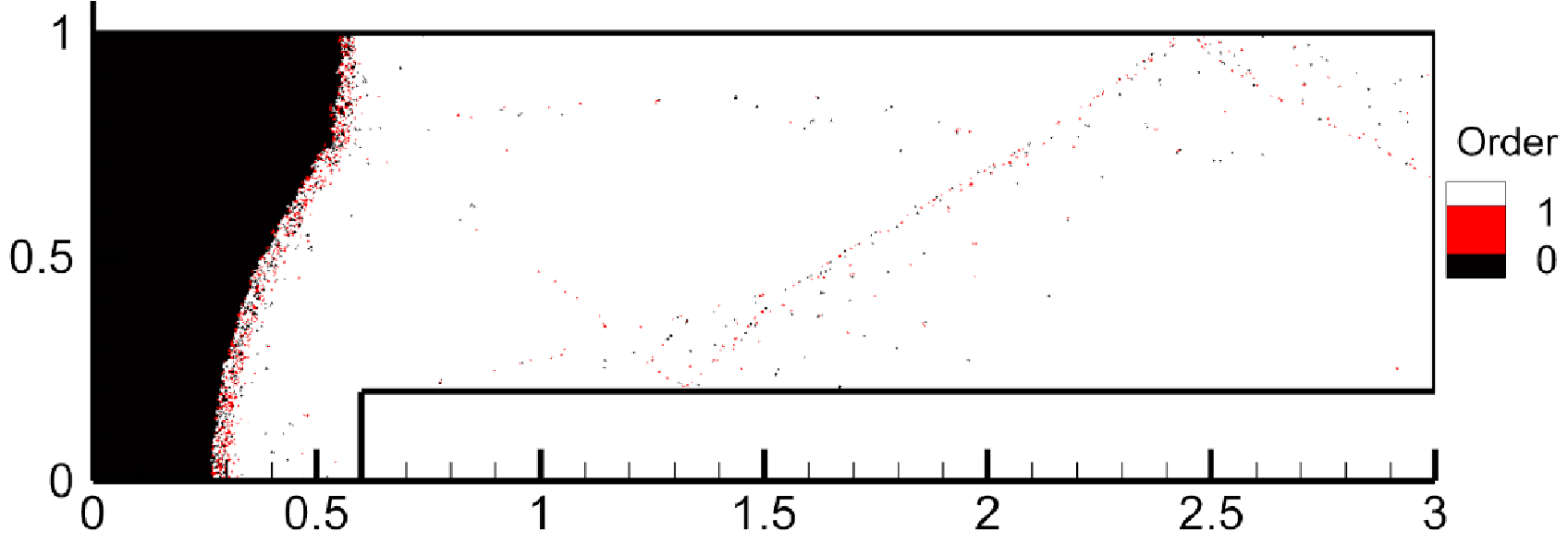}}
  \subfigure[$k=3$]{
  \includegraphics[height=3 cm]{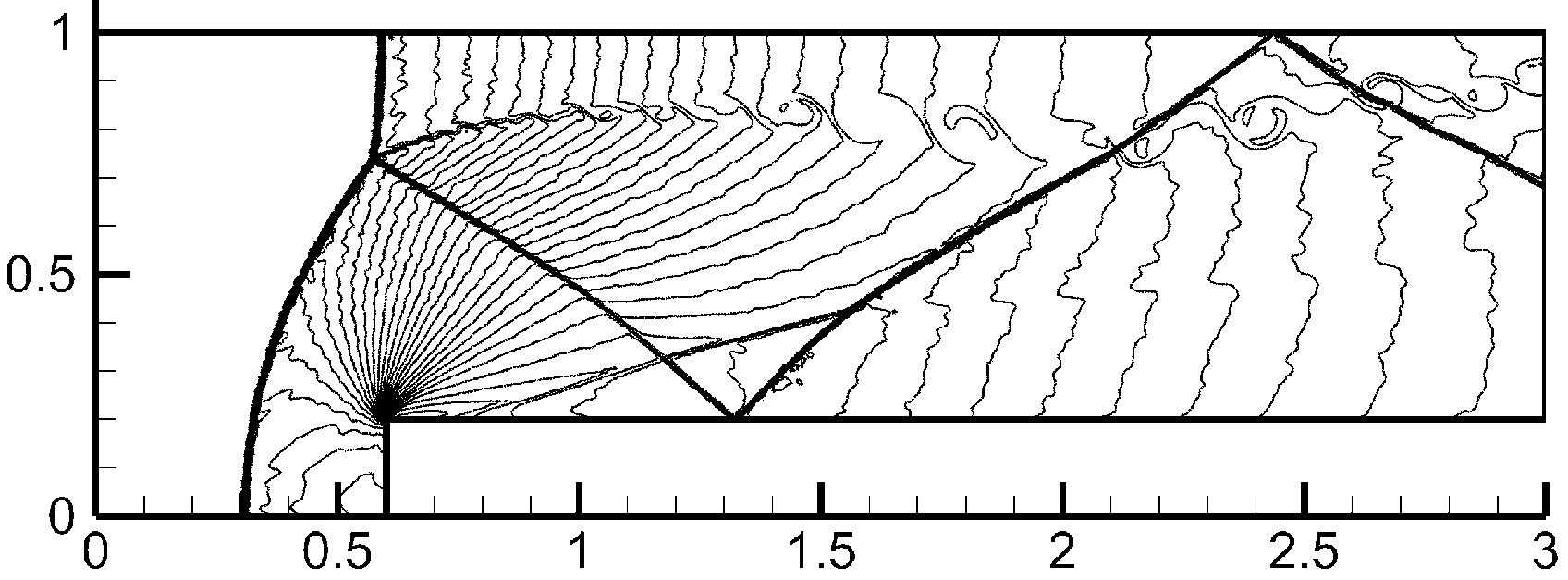}}
  \subfigure[$k=3$]{
  \includegraphics[height=3 cm]{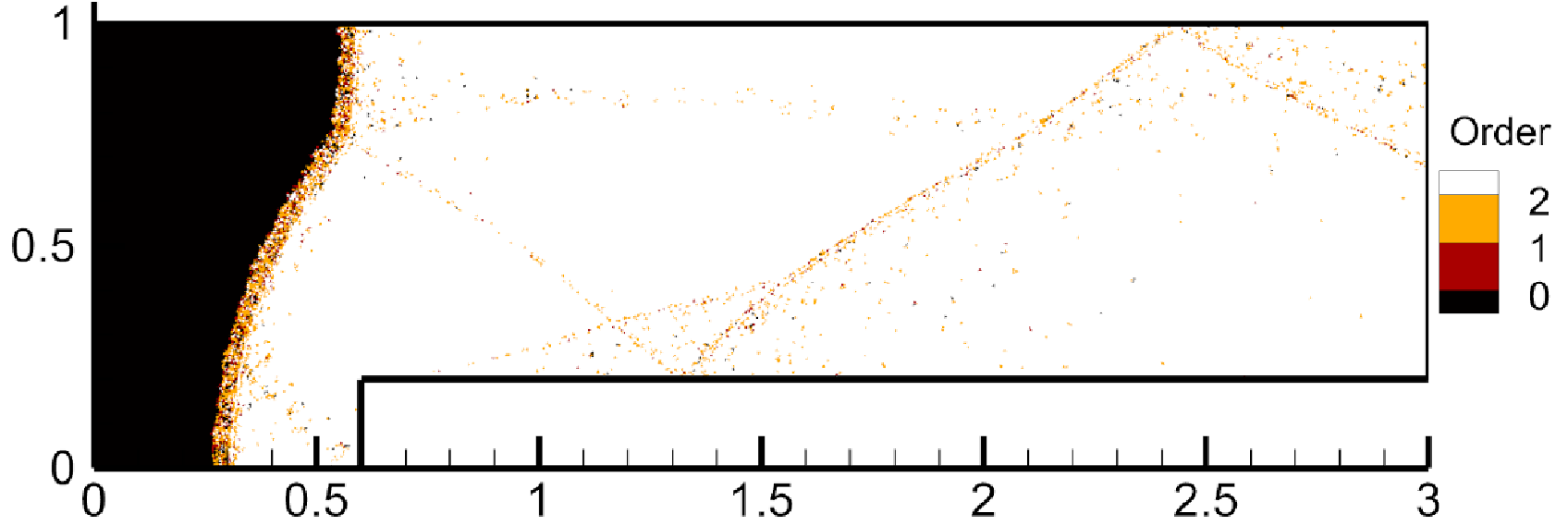}}
  \caption{Density contours (left column) and the polynomial order (right column) of the forward-facing step problem at $t=4$ calculated 
  by RKDG schemes with the MR limiter using quasi-uniform triangular cells.
  The density contours contain 30 equidistant contours from 0.32 to 6.15.
  The white parts of the polynomial order plots represent the original $k$th-order DG polynomial.}
 \label{FIG:FST_Tri}
 \end{figure}

Fig. \ref{FIG:FST_Qua} and Fig. \ref{FIG:FST_Tri} respectively show the density contours and the polynomial order
computed by the RKDG schemes with the MR limiter using uniform quadrilateral meshes and quasi-uniform triangular meshes.
The edge length of both the quadrilateral cells and triangular cells on the boundaries is $1/160$.
We can see that the RKDG schemes with the MR limiter can capture the shock waves without numerical oscillations.
When $k\ge2$, flow instabilities can be resolved by using the current mesh size.

\section{Conclusions}
We proposed a multi-resolution limiter for the RKDG schemes.
Extensive numerical examples show that the MR limiter can eliminate numerical oscillations near discontinuities 
while preserving the high-order accuracy of the RKDG schemes.
The MR limiter works well for different problems and different meshes without tuning the limiter. 
The search for the optimal threshold $C_k$ for different orders, rather than using a constant $C_k$,
may further improve the performance of the MR limiter.

\section*{Acknowledgments}
H. S. was partially supported by the National Natural Science Foundation of China (No. 62231016) and the Innovation Group Grant of Sichuan Province (No. 2025NSFTD0004)
B. S. was partially supported by the  National Natural Science Foundation of China (No. 12201437). 
\bibliography{mybibfile}

\begin{thebibliography}{10}
\expandafter\ifx\csname url\endcsname\relax
  \def\url#1{\texttt{#1}}\fi
\expandafter\ifx\csname urlprefix\endcsname\relax\def\urlprefix{URL }\fi
\expandafter\ifx\csname href\endcsname\relax
  \def\href#1#2{#2} \def\path#1{#1}\fi

\bibitem{Reed1973DG}
W.~H. Reed, T.~Hill, Triangular mesh methods for the neutron transport
  equation, Los Alamos Report LA-UR-73-479.

\bibitem{Cockburn1989RKDGII}
B.~Cockburn, C.-W. Shu, {TVB} {R}unge-{K}utta local projection discontinuous
  {G}alerkin finite element method for conservation laws. {II}. {G}eneral
  framework, Mathematics of Computation 52~(186) (1989) 411--435.

\bibitem{Cockburn1989RKDGIII}
B.~Cockburn, S.-Y. Lin, C.-W. Shu, {TVB} {R}unge-{K}utta local projection
  discontinuous {G}alerkin finite element method for conservation laws. {III}:
  {O}ne-dimensional systems, Journal of Computational Physics 84~(1) (1989)
  90--113.

\bibitem{Cockburn1990RKDGIV}
B.~Cockburn, S.~Hou, C.-W. Shu, The {R}unge-{K}utta local projection
  discontinuous {G}alerkin finite element method for conservation laws. {IV}.
  {T}he multidimensional case, Mathematics of Computation 54~(190) (1990)
  545--581.

\bibitem{Cockburn1998RKDGV}
B.~Cockburn, C.-W. Shu, The {R}unge-{K}utta local projection discontinuous
  {G}alerkin finite element method for conservation laws. {V}:
  {M}ultidimensional systems, Journal of Computational Physics 141~(2) (1998)
  199--224.

\bibitem{Shu1988EfficientENO}
C.-W. Shu, S.~Osher, Efficient implementation of essentially non-oscillatory
  shock-capturing schemes, Journal of Computational Physics 77 (1988) 439--471.

\bibitem{Shu1989EfficientENOII}
C.-W. Shu, S.~Osher, Efficient implementation of essentially non-oscillatory
  shock-capturing schemes, ii, Journal of Computational Physics 83 (1989)
  32--78.

\bibitem{Shu1987TVB}
C.-W. Shu, {TVB} uniformly high-order schemes for conservation laws,
  Mathematics of Computation 49~(179) (1987) 105--121.

\bibitem{Biswas1994}
R.~Biswas, K.~D. Devine, J.~E. Flaherty, Parallel, adaptive finite element
  methods for conservation laws, Applied Numerical Mathematics 14~(1-3) (1994)
  255--283.

\bibitem{Burbeau2001}
A.~Burbeau, P.~Sagaut, C.-H. Bruneau, A problem-independent limiter for
  high-order {R}unge--{K}utta discontinuous {G}alerkin methods, Journal of
  Computational Physics 169~(1) (2001) 111--150.

\bibitem{Qiu2004WENOLimiter}
J.~Qiu, C.-W. Shu, Hermite {WENO} schemes and their application as limiters for
  {R}unge--{K}utta discontinuous {G}alerkin method: one-dimensional case,
  Journal of Computational Physics 193~(1) (2004) 115--135.

\bibitem{Qiu2005WENOLimiter}
J.~Qiu, C.-W. Shu, Runge--{K}utta discontinuous {G}alerkin method using {WENO}
  limiters, SIAM Journal on Scientific Computing 26~(3) (2005) 907--929.

\bibitem{Qiu2005ComparisonIndicator}
J.~Qiu, C.-W. Shu, A comparison of troubled-cell indicators for
  {R}unge--{K}utta discontinuous {G}alerkin methods using weighted essentially
  nonoscillatory limiters, SIAM Journal on Scientific Computing 27~(3) (2005)
  995--1013.

\bibitem{Harten1989Subcell}
A.~Harten, Eno schemes with subcell resolution, Journal of Computational
  Physics 83~(1) (1989) 148--184.

\bibitem{Krivodonova2004KXRCF}
L.~Krivodonova, J.~Xin, J.-F. Remacle, N.~Chevaugeon, J.~E. Flaherty, Shock
  detection and limiting with discontinuous galerkin methods for hyperbolic
  conservation laws, Applied Numerical Mathematics 48~(3-4) (2004) 323--338.

\bibitem{Zhu2008WENOLimiterTri}
J.~Zhu, J.~Qiu, C.-W. Shu, M.~Dumbser, {R}unge--{K}utta discontinuous
  {G}alerkin method using {WENO} limiters ii: unstructured meshes, Journal of
  Computational Physics 227~(9) (2008) 4330--4353.

\bibitem{Zhu2012WENOLimiter3D}
J.~Zhu, J.~Qiu, {R}unge-{K}utta discontinuous {G}alerkin method using
  {WENO}-type limiters: three-dimensional unstructured meshes, Communications
  in Computational Physics 11~(3) (2012) 985--1005.

\bibitem{Zhong2013SimpleDG_limiter}
X.~Zhong, C.-W. Shu, A simple weighted essentially nonoscillatory limiter for
  {R}unge--{K}utta discontinuous {G}alerkin methods, Journal of Computational
  Physics 232~(1) (2013) 397--415.

\bibitem{Zhu2013SimpleWENOLimiterTri}
J.~Zhu, X.~Zhong, C.-W. Shu, J.~Qiu, {R}unge--{K}utta discontinuous {G}alerkin
  method using a new type of {WENO} limiters on unstructured meshes, Journal of
  Computational Physics 248 (2013) 200--220.

\bibitem{Dumbser2014subcell}
M.~Dumbser, O.~Zanotti, R.~Loub{\`e}re, S.~Diot, A posteriori subcell limiting
  of the discontinuous {G}alerkin finite element method for hyperbolic
  conservation laws, Journal of Computational Physics 278 (2014) 47--75.

\bibitem{Clain2011MOOD}
S.~Clain, S.~Diot, R.~Loub{\`e}re, A high-order finite volume method for
  systems of conservation laws—{M}ulti-dimensional {O}ptimal {O}rder
  {D}etection ({MOOD}), Journal of Computational Physics 230~(10) (2011)
  4028--4050.

\bibitem{Du2022ImprovedWENOLimiterTri}
J.~Du, C.-W. Shu, X.~Zhong, An improved simple {WENO} limiter for discontinuous
  {G}alerkin methods solving hyperbolic systems on unstructured meshes, Journal
  of Computational Physics 467 (2022) 111424.

\bibitem{Wei2024HybridLimiter}
L.~Wei, Y.~Xia, An indicator-based hybrid limiter in discontinuous galerkin
  methods for hyperbolic conservation laws, Journal of Computational Physics
  498 (2024) 112676.

\bibitem{Wei2025JumpFilter}
L.~Wei, L.~Zhou, Y.~Xia, The jump filter in the discontinuous {G}alerkin method
  for hyperbolic conservation laws, Journal of Computational Physics 520 (2025)
  113498.

\bibitem{Lu2021OFDG}
J.~Lu, Y.~Liu, C.-W. Shu, An oscillation-free discontinuous {G}alerkin method
  for scalar hyperbolic conservation laws, SIAM Journal on Numerical Analysis
  59~(3) (2021) 1299--1324.

\bibitem{Liu2022OFDG}
Y.~Liu, J.~Lu, C.-W. Shu, An essentially oscillation-free discontinuous
  {G}alerkin method for hyperbolic systems, SIAM Journal on Scientific
  Computing 44~(1) (2022) A230--A259.

\bibitem{Peng2025OEDG}
M.~Peng, Z.~Sun, K.~Wu, {OEDG}: {O}scillation-eliminating discontinuous
  {G}alerkin method for hyperbolic conservation laws, Mathematics of
  Computation 94~(353) (2025) 1147--1198.

\bibitem{Ding2025OEDG}
S.~Ding, S.~Cui, K.~Wu, Robust {DG} schemes on unstructured triangular meshes:
  {O}scillation elimination and bound preservation via optimal convex
  decomposition, Journal of Computational Physics 526 (2025) 113769.

\bibitem{Fu2017NewLimiter}
G.~Fu, C.-W. Shu, A new troubled-cell indicator for discontinuous galerkin
  methods for hyperbolic conservation laws, Journal of Computational Physics
  347 (2017) 305--327.

\bibitem{Vuik2016automated}
M.~J. Vuik, J.~K. Ryan, Automated parameters for troubled-cell indicators using
  outlier detection, SIAM Journal on Scientific Computing 38~(1) (2016)
  A84--A104.

\bibitem{Ray2018ANN}
D.~Ray, J.~S. Hesthaven, An artificial neural network as a troubled-cell
  indicator, Journal of Computational Physics 367 (2018) 166--191.

\bibitem{Ray2019ANN}
D.~Ray, J.~S. Hesthaven, Detecting troubled-cells on two-dimensional
  unstructured grids using a neural network, Journal of Computational Physics
  397 (2019) 108845.

\bibitem{Zhu2021K_means}
H.~Zhu, H.~Wang, Z.~Gao, A new troubled-cell indicator for discontinuous
  {G}alerkin methods using {K}-means clustering, SIAM Journal on Scientific
  Computing 43~(4) (2021) A3009--A3031.

\bibitem{Zhang2010MaximumDG}
X.~Zhang, C.-W. Shu, On maximum-principle-satisfying high order schemes for
  scalar conservation laws, Journal of Computational Physics 229~(9) (2010)
  3091--3120.

\bibitem{Zhang2010PositivityDG}
X.~Zhang, C.-W. Shu, On positivity-preserving high order discontinuous
  {G}alerkin schemes for compressible {E}uler equations on rectangular meshes,
  Journal of Computational Physics 229~(23) (2010) 8918--8934.

\bibitem{Zhang2011MaximumReview}
X.~Zhang, C.-W. Shu, Maximum-principle-satisfying and positivity-preserving
  high-order schemes for conservation laws: survey and new developments,
  Proceedings of the Royal Society A: Mathematical, Physical and Engineering
  Sciences 467~(2134) (2011) 2752--2776.

\bibitem{Zhang2012MaximumDGTri}
X.~Zhang, Y.~Xia, C.-W. Shu, Maximum-principle-satisfying and
  positivity-preserving high order discontinuous {G}alerkin schemes for
  conservation laws on triangular meshes, Journal of Scientific Computing
  50~(1) (2012) 29--62.

\bibitem{Zhang2017PositivityDGNS}
X.~Zhang, On positivity-preserving high order discontinuous {G}alerkin schemes
  for compressible {N}avier--{S}tokes equations, Journal of Computational
  Physics 328 (2017) 301--343.

\bibitem{Shen2025ENO_MR}
H.~Shen, An efficient class of increasingly high-order {ENO} schemes with
  multi-resolution, Computers \& Fluids 291 (2025) 106589.

\bibitem{Shen2020RotatedDecompositon}
H.~Shen, M.~Parsani, A rotated characteristic decomposition technique for
  high-order reconstructions in multi-dimensions, Journal of Scientific
  Computing 88~(87).

\bibitem{Jiang_Shu1996WENO}
G.-S. Jiang, C.-W. Shu, Efficient implementation of weighted {ENO} schemes,
  Journal of Computational Physics 126~(1) (1996) 202--228.

\bibitem{Woodward1984JCP}
P.~Woodward, P.~Colella, The numerical simulation of two-dimensional fluid flow
  with strong shocks, Journal of computational physics 54~(1) (1984) 115--173.

\end{thebibliography}

\end{document}